\documentclass[final,12pt]{elsarticle}

\usepackage{mathptmx}      

\usepackage{color}

\usepackage{setspace} 

\usepackage{lscape}
\usepackage{rotating}
\usepackage{geometry}
\usepackage{longtable}

\usepackage{amsmath}
\usepackage{epstopdf}
\usepackage[T1]{fontenc}
\usepackage[latin1]{inputenc}
\usepackage{alltt}   %
\usepackage{color}   %
\usepackage{xcolor}  %
\usepackage{lmodern} %
\usepackage[footnotesize]{caption} %
\usepackage{verbatim}    
\usepackage{bm} 
\usepackage{float}      %
\usepackage{morefloats} %

\usepackage{amsfonts}
 \usepackage{amssymb,amsmath}

\DeclareMathAlphabet{\mathbfsf}{\encodingdefault}{\sfdefault}{bx}{sl}

\makeatletter
\newcommand*{\centerfloat}{%
  \parindent \z@
  \leftskip \z@ \@plus 1fil \@minus \textwidth
  \rightskip\leftskip
  \parfillskip \z@skip}
\makeatother

\journal{}

\begin{document}

\begin{frontmatter}

\title{A continuous/discontinuous Galerkin solution of the shallow water equations with dynamic viscosity, high-order wetting and drying, and implicit time integration}


\author[1]{Simone Marras\corref{cor1}}\ead{smarras@stanford.edu}
\author[2]{Michal A. Kopera}
\author[3,4]{Emil M. Constantinescu}
\author[1]{Jenny Suckale}
\author[2]{Francis X. Giraldo}

\cortext[cor1]{Corresponding author} 

\address[1]{Stanford University, Dept. of Geophysics, Stanford, CA 94305, U.S.A.}
\address[2]{Naval Postgraduate School, Dept. of Applied Mathematics, Monterey, CA 93943, U.S.A.}
\address[3]{Argonne National Laboratory, Mathematics and Computer Science, Argonne, IL 60439, U.S.A.}
\address[4]{The University of Chicago, Computation Institute, Chicago, IL 60637, U.S.A.}

\begin{abstract}
The high-order numerical solution of the non-linear shallow water equations (and of hyperbolic systems in general) 
is susceptible to unphysical Gibbs oscillations that form in the proximity of strong gradients.
The solution to this problem is still an active field of research as no general cure has been found yet.
In this paper, we tackle this issue by presenting a dynamically adaptive viscosity based on a residual-based sub-grid scale model
that has the following properties:
$(i)$ it removes the spurious oscillations in the proximity of strong wave fronts while preserving 
the overall accuracy and sharpness of the solution. This is possible because of the residual-based definition of the 
dynamic diffusion coefficient.
$(ii)$ For coarse grids, it prevents energy from building up at small wave-numbers. $(iii)$ The model has no tunable parameter.

Our interest in the shallow water equations is tied to the simulation of coastal inundation, where a careful handling of the transition between dry and wet surfaces is particularly challenging for high-order Galerkin approximations. In this paper, we extend to a unified continuous/discontinuous Galerkin (CG/DG) framework a very simple, yet effective wetting and drying algorithm originally designed for DG [Xing, Zhang, Shu (2010) {\it Positivity-preserving high order well-balanced discontinuous Galerkin methods for the shallow water equations} Adv. Water Res. 33:1476-1493]. 
We show its effectiveness for problems in one and two dimensions 
on domains of increasing characteristic lengths varying from centimeters to kilometers.
 

Finally, to overcome the time-step restriction incurred by the high-order Galerkin approximation, 
we advance the equations forward in time via a three stage, second
order explicit-first-stage, singly diagonally implicit Runge-Kutta 
(ESDIRK) time integration scheme. Via ESDIRK, we are able to preserve numerical stability for an advective CFL number 10 times larger than 
its explicit counterpart.
\end{abstract}
\end{frontmatter}

\section{Introduction}
The shallow water equations (Saint-Venant \cite{Saint-Venant1871}) are a common $(d - 1)$ approximation to the $d$--dimensional Navier-Stokes equations to model incompressible, free surface flows.
Due to the ability of high-order Galerkin methods to keep dissipation and dispersion errors low (Ainsworth et al.\ \cite{ainsworthEtAl2006})
and their flexibility with arbitrary geometries and {\it hp}-adaptivity, 
these methods are 
proving their mettle for solving the shallow water equations (SW) in the modeling of non-linear waves in different geophysical flows (Ma \cite{ma1993}, Iskandarani \cite{iskandaraniEtAl1995}, Taylor et al.\ \cite{taylorTribbia1997}, Giraldo \cite{giraldo2001}, Giraldo et al.\ \cite{giraldoHesthavenWarburton2002}), Dawson and Aizinger \cite{dawsonAizinger2005}, Kubatko et al.\ \cite{kubatkoEtAl2006}, Nair et al.\ \cite{nairThomasLoft2007}, Giraldo and Restelli \cite{giraldoRestelli2010},  Xing et al.\ \cite{xingZhangShu2010}, K\"arn\"a et al.\ \cite{karnaEtAl2011}, Hendricks et al.\ \cite{hendricksKoperaGiraldo2015}, Hood \cite{karolineTHESIS2016}).
Although it is typically assumed that high-order Galerkin methods are not strictly necessary, they do offer many advantages over their low-order counterparts. Examples include their ability to resolve fine scale structures and to do so with fewer degrees of freedom, as well as their strong scaling properties on massively parallel computers (M\"uller et al.\ \cite{mullerEtAl2016}, Abdi et al.\ \cite{abdiEtAl2016}, Gandhem et al.\ \cite{gandhamEtAl2015}).
High-order methods are often attributed with some disadvantages as well. For example, they are constrained by 
small time-steps. To overcome this restriction, 
we follow K\"arn\"a et al.\ \cite{karnaEtAl2011} and implement an implicit Runge-Kutta scheme based on Giraldo et al.\ \cite{giraldoEtAl2013}. 
Furthermore, the numerical approximation of non-linear hyperbolic systems via high-order methods is susceptible to unphysical Gibbs oscillations that form in the proximity of strong gradients such as those of sharp wave fronts (e.g. bores).
Filters such as Vandeven's and Boyd's \cite{vandeven1991, boyd1996} are the most common tools to handle this problem, as well as artificial diffusion of some sort. However, we noticed that filtering may not be sufficient as the flow strengthens and the 
wave sharpness intensifies. We have recently shown this issue with the solution of the non-linear Burgers' equation by high order spectral elements in \cite[\S 5]{marrasEtAl2015b} and will show how a proper dynamic viscosity does indeed improve on filters for the shallow water system as well. 
To preserve numerical stability without compromising the overall quality of the solution, Pham Van et al.\ \cite{phanVanEtAl2014} and Rakowsky et al.\ \cite{rakowskyEtAl2013} utilized a Lilly-Smagorinsky model \cite{lilly1962, smagorinsky1963}. 
To account for sub-grid scale effects, viscosity was also utilized in the DG model described by 
Gourgue et al.\ \cite{gourgueEtal2009} 
to improve their inviscid simulations.
Recently, Pasquetti et al.\ \cite{pasquettiGuermondPopoICOSAHOM2014} stabilized the high-order spectral element solution of the one-dimensional Saint-Venant equations via the entropy viscosity method.
Michoski et al.\ \cite{michoskiEtAl2016} compared artificial viscosity, limiters, and filters for the (modal) DG solution of SW, 
concluding that a dynamically adaptive diffusion may be the most effective means of regularization at higher orders. 
The diffusion model that we present in this work stems from the sub-grid scale (SGS) model that 
was proposed by Nazarov and Hoffman \cite{nazarovHoffman2013} 
to stabilize the linear finite element solution of compressible flows with shock waves. 
This was later applied to stabilize high-order continuous and discontinuous Galerkin (CG/DG)
in the context of stratified, low Mach number atmospheric flows by Marras et al.\ in \cite{marrasEtAl2015b}.
This approach is based on the ideas of scales splitting in large eddy simulation, where
the scales of physical importance are split into 
grid resolvable and unresolvable. The unresolved scales are then parameterized via the SGS model at hand ({\it Dyn-SGS}). 
{\it Dyn-SGS} is applicable to any numerical method and is not limited to Galerkin approximations; 
its implementation in existing numerical models is straightforward. 
Stabilization becomes more important as the grid resolution is decreased. The metric that we utilize to 
evaluate its importance is the spectrum of kinetic energy across wave numbers.
To understand the impact of diffusion 
with respect to the energetic behavior of the numerical solution, we compare the spectra obtained for the 
stabilized continuous and discontinuous Galerkin solutions against their inviscid counterparts. 
As viscosity is accounted for, the solution improves at coarser resolutions for both CG and DG. 
In the case of DG, however, dissipation is already built-in by the flux computation across elements.
For this reason, the amount of artificial diffusion that DG requires for stabilization is smaller than the one required by CG.
This property of DG has been exploited to build implicit Large Eddy Simulation (ILES) models (e.g.,  Uranga et al.\ \cite{urangaEtAl2011}). ILES relies on the numerical dissipation of the discretization method to 
model the effects on the unresolved, sub-grid scale eddies rather than using an explicit sub-grid scale (SGS) model, under the requirement of 
high-resolution \cite{rider2006, drikakisEtAl2007}. The MPDATA method of Smolarkiewicz \cite{smolarkiewicz1983,smolarkiewicz1984} is an example of such an ILES model used for modeling geophysical flows.

Finally, there is a known difficulty in including wetting and drying algorithms (typically designed for low-order methods) while preserving high-order accuracy. 
The application of wetting/drying with discontinuous Galerkin using low-order Lagrange polynomials  can be found in, e.g., Bunya et al.\ \cite{bunyaEtAl2009}, Vater et al.\ \cite{vaterEtAl2015}, K\"arn\"a et al.\ \cite{karnaEtAl2011}, 
or Gourgue et al . cite{gourgueEtal2009}, and using Bernstein polynomials up to order three in 
Beisiegel and Behrens \cite{beisiegelBehrens2015}.
The positivity preserving limiter of Xing et al.\ \cite{xingZhangShu2010} (Xing-Zhang-Shu limiter from now on) 
was designed for high-order discontinuous Galerkin to solve this problem in particular. Because it is mass conservative and preserves the global high order accuracy of the solution, we implemented it in our 
model and extended it to continuous Galerkin as well. 
This limiter is guaranteed to work if and only if the mean water height is positive. 
Its designers recommend using the total variation bounding (TVB) limiter of Shu \cite{shu1987, cockburnShu1989} to make the solution positive definite. Because the TVB limiter requires a bound on the derivative of the solution, which we do not know a priori, we rely on the artificial dissipation introduced above instead.

Looking at the larger picture of things, the ever increasing interest in the high-order accuracy of inundation models stems from a long list of 
coastal disasters in the last 10 years alone. From the 2004 tsunami in the Indian Ocean that caused 280,000 deaths, to the 2011 T\-ohoku tsunami, followed by Superstorm Sandy in 2012, the devastating 2013 typhoons in the Philippines and India, and the 2014 hurricane Odile; the list seems to be getting longer rapidly as the global climate patterns are changing.
To study the impact of similar events in the future, coastal planners around the world are more and more relying 
on numerical tools such as the one that we present in this paper.


The rest of the paper is organized as follows.
The shallow water equations are defined in \S \ref{equationsSct}. Their numerical solution is described in \S \ref{discreteSCT}. 
The inundation model is described in \S \ref{wetDrySECT} and is followed by the derivation of the SGS model in \S \ref{stabilizationSCT}. Numerical results are reported in \S \ref{testsSCT}. The conclusions are drawn in \S \ref{conclusionsSCT}.

\section{Equations}
\label{equationsSct}
Let $\Omega\in\mathbb{R}^d$ be a fixed domain of space dimension $d$ with boundary $\Gamma$ and Cartesian coordinates ${\bf x}=x$ in 1D and ${\bf x}=[x,y]$ in 2D 
and let $t\in\mathbb{R}^+$ identify time. Let $z$ always identify the direction of gravity.
The absolute free surface level and bathymetry are identified by the symbols $H_s(t,{\bf x})$ and, respectively, $H_b({\bf x})$, so that $H=H_s + H_b$.
The viscous shallow water equations are then given as:

\begin{subequations}
\label{SWE}
\begin{equation}
\label{massEqn}
\frac{\partial H}{\partial t} + \nabla\cdot \left(H{\bf u}\right) = \delta\nabla\cdot(\mu_{SGS}\nabla H),
\end{equation}
\begin{equation}
\label{momeEqn}
\frac{\partial H{\bf u}}{\partial t} + \nabla\cdot\left(H{\bf u}\otimes {\bf u} + \frac{g}{2}(H^2 - H^2_{b}) {\bf I} \right) + g H_s\nabla\cdot\left(H_{b} {\bf I} \right) =  \nabla\cdot \left(H_s\mu_{SGS}\nabla{\bf u}\right), 
\end{equation}
\end{subequations}
where $g=9.81\,{\rm m\,s^{-2}}$ is the acceleration of gravity, ${\bf I}$ is the $d\times d$ identity matrix, and $\mu_{SGS}$ is the dynamic dissipation coefficient to be defined shortly. In (\ref{massEqn}), the $\delta$ coefficient defines whether viscosity is turned on ($\delta=1$) or off ($\delta=0$)  in the continuity equation. We will return to this topic shortly.

%

\section{Space and time discretization}
\label{discreteSCT}
High-order spectral element (SEM or CG, for continuous Galerkin) and discontinuous Galerkin (DG) approximations on quadrilateral elements are used to discretize Eq.\  (\ref{SWE}). 

The numerical model used in this paper is in fact the two-dimensional version of the NUMA model described in Giraldo et al.\ \cite{giraldoHesthavenWarburton2002} and in Abdi and Giraldo \cite{abdiGiraldo2016}. Furthermore, the model is derived from the model described in Kopera and Giraldo \cite{koperaGiraldo2013a, koperaGiraldo2013b} which is here used for the shallow water equations.
The solution is advanced in time using a fully implicit Runge-Kutta
scheme (see \S \ref{sec:time:integration}).

\subsection{Spectral element and discontinuous Galerkin approximations}
We leave the details of the discretization to the work of Giraldo et al. \cite{giraldoHesthavenWarburton2002} and 
Kopera and Giraldo \cite{koperaGiraldo2013a, koperaGiraldo2013b}. 
Nonetheless, we introduce some notation that we are going to use later in the paper. 
To solve the shallow water equations by element-based Galerkin methods on a domain $\Omega$, we proceed by defining the weak form of (\ref{SWE}) that we first recast in compact notation as 

\begin{equation}
\label{compactSW}
\frac{\partial {\bf q}}{\partial t} + \nabla\cdot{\bf F}({\bf q}) = {\bf S}({\bf q}),
\end{equation}
where ${\bf q} = [H, H{\bf u}]^{\rm T}$ is the transposed array of the solution variables and ${\bf F}$ and ${\bf S}$ are the flux and source terms.

In the case of spectral elements, the space discretization yields the semi-discrete matrix problem
\begin{equation}
\label{matrixprobCG}
\frac{\partial {\bf q}}{\partial t} = \widehat{\bf D}^{\rm T}{\bf F}({\bf q}) + {\bf S}({\bf q})\,,
\end{equation}
where, for the global mass and differentiation matrices, ${\bf M}$ and ${\bf D}$, 
we have that $\widehat{\bf D}={\bf M}^{-1}{\bf D}$. The global matrices are obtained from their local counterparts, ${\bf M}^e$ and ${\bf D}^e$, by direct stiffness summation, which maps the local degrees of
freedom of an element $\Omega^h_e$ to the corresponding global degrees of freedom in $\Omega^h$, and 
adds the element values in the global system.
By construction, ${\bf M}$ is diagonal (assuming inexact integration), with an obvious advantage if explicit time integration is used. 

In the discontinuous Galerkin approximation, the problem at hand is solved only locally, and
unlike the case of spectral elements, the flux integral that stems from the integration by-parts must be discretized as well. 
The element-wise counterpart of the matrix problem (\ref{matrixprobCG}) is then written as:
\begin{equation}
\label{matrixprobDG}
\frac{\partial {\bf q}^e}{\partial t} = - (\widehat{\bf M}^{\Gamma,e})^{\rm T}\breve{\bf F}({\bf q}^e) + (\widehat{\bf D}^e)^{\rm T}{\bf F}({\bf q}^e) + {\bf S}({\bf q}^e),
\end{equation}
where we obtain $\widehat{\bf M}^{\Gamma,e}=({\bf M}^e)^{-1}{\bf M}^{\Gamma,e}$ from the element boundary matrix, ${\bf M}^{\Gamma,e}$, 
and the element mass matrix, ${\bf M}^e$.
Out of various possible choices for the definition of the numerical flux $\breve{\bf F}({\bf q})$ in Eq.\ (\ref{matrixprobDG}), 
we adopted the Rusanov flux.
The Laplace operator of viscosity is approximated using the Symmetric Interior Penalty method (SIP; the reader is referred to Arnold \cite{arnold1982} for details on its definition).

%
\subsection{Time integration\label{sec:time:integration}}
%

Equation \eqref{matrixprobCG} is integrated
in time by an implicit Runge-Kutta (RK) scheme that corresponds to
the implicit part of the implicit-explicit scheme used in
\cite{giraldoEtAl2013} (see also \cite{Butcher_1999}). The method
coefficients in standard ($A=a_{ij},\,b,\,c$) tableaux form are the
following  
\newcommand\ST{\rule[-0.75em]{0pt}{2em}}
\begin{align}
\label{eq:ARK:s3:p2:q2:LSTABLE}
& \begin{array}{c|ccc}
\ST 0&0\\
\ST 2-\sqrt{2}&1-\frac{1}{\sqrt{2}}&1-\frac{1}{\sqrt{2}}\\
\ST 1&\frac{1}{2\sqrt{2}}&\frac{1}{2\sqrt{2}}&1-\frac{1}{\sqrt{2}}\\
\cline{1-4}
\ST  &\frac{1}{2\sqrt{2}}&\frac{1}{2\sqrt{2}}&1-\frac{1}{\sqrt{2}}
\end{array}\,,\qquad \begin{array}{c|c}c&A\\ \hline &b\end{array}.
\end{align}
Scheme \eqref{eq:ARK:s3:p2:q2:LSTABLE} is a three stage second order
explicit-first-stage singly diagonally implicit RK (ESDIRK) scheme. This
scheme has desirable accuracy and stability properties: $(i)$ all
stages are second order accurate, $(ii)$ it is stiffly accurate and
$L$-stable, and $(iii)$ it is strong-stability-preserving (SSP)
\cite{gottlieb2001strong} with SSP coefficient of 2. These
properties allow us to take large time-steps with high
accuracy as well as alleviate the instability issues associated with sharp
solution gradients \cite{gottlieb2001strong}. 
The two-dimensional tests presented later in this paper showed to be 
the most demanding in terms of stability constraints.
Method \eqref{eq:ARK:s3:p2:q2:LSTABLE} allows us to gain up to one
order of magnitude in terms of maximum admissible advective CFL when
compared to an explicit method (explicit part of ARK3, \cite{Kennedy_2001}).
In particular, the explicit four-stage Runge-Kutta solution of the solitary wave 
against one isolated obstacle described in \S \ref{OneislandTest} 
preserved stability for up to CFL=0.21 using both CG and DG approximations. Although 
we were not able to use arbitrarily large time-steps with the ESDIRK with the current implementation (we will address this issue in a future work), 
we resolved the same problem at CFL = 1.8. 
Schemes with a subset of these properties are
employed by K\"{a}rn\"{a} et al.\ \cite{karnaEtAl2011} and shown to be robust in this context.
Method \eqref{eq:ARK:s3:p2:q2:LSTABLE} used in this study satisfies all properties ($i$-$iii$).   

Computationally, at each of the two implicit stages we have to solve a
nonlinear equation ${\bf G}({\bf Q}^{(i)})=0$, where ${\bf Q}^{(i)}$
are the stage values, $i=1,2$. We do so by using Newton iterations
with a stopping criterion based on the relative decrease in the
residual; that is, stop at iteration $k$ if $||{\bf G}({\bf
  Q}^{(i)}_k)||/||{\bf G}({\bf Q}^{(i)}_0)||<Tol_N$. 
At each Newton iteration we have to solve a linear system ${\bf J}({\bf
  Q}^{(i)}_k - {\bf
  Q}^{(i)}_{k-1}) = - {\bf G}({\bf Q}^{(i)}_{k-1})$, where ${\bf J}$
is the Jacobian matrix of ${\bf G}({\bf Q}^{(i)})$.
%
%
We approximate the Jacobian using directional
finite differences and iterate with the generalized minimal residual (GMRES)
method, which is effectively a
Jacobian-free Newton--Krylov method \cite{knoll2004jacobian}. The
GMRES stopping criterion is also based on the relative residual. 
The
first stage is explicit and equal to the last stage of the previous
step, effectively making it a two-stage method, which saves some computational time.


\section{The Dyn-SGS model for the shallow water equations}
\label{stabilizationSCT}
There are different ways to derive the viscous model described by Eq.\ (\ref{SWE}) from the inviscid Saint-Venant equations. Analogous to our previous work on the large eddy simulation of stratified atmospheric flows \cite{marrasEtAl2015b}, the current model
builds on the separation between grid resolved (indicated as $\overline{f}({\bf x})$ for any quantity $f({\bf x})$) and unresolved (sub-grid) scales.
The unresolved scales are modeled via the eddy viscosity scheme described in this paper.
Given an element $\Omega_e$ of order $N$ and with side lengths $\Delta x, \Delta y$ of comparable orders of magnitude,
we define the following characteristic length (and hence filter width \cite{sagautBook}):
\[
\overline\Delta = \min{\left(\Delta x, \Delta y\right)}/(N+1).
\]
The value of $\Delta$ sets the size of the smallest resolvable scales in the same way as cut-off filters do in large eddy simulation 
models.

The application of this filter to the continuity equation (\ref{massEqn}) 
results in the presence of an additional sub-grid term on the right-hand side. 
It is often debated whether diffusion should be applied to the continuity equation \cite{pasquettiGuermondPopoICOSAHOM2014, gerbeauEtAl2001,Guermond_Popov_2014}; 
clearly, should the discrete viscous operator not be conservative, mass dissipation should not be used.
However, by relying on spectral elements with integration by parts of the second-order diffusion operator, 
the discrete viscous operator is conservative, as shown in \cite{gubaEtAl2014}. To get a sense of how necessary a stabilized continuity equation may be, we will show a few results for both conditions in \S \ref{OneislandTest}.

Scale separation in the momentum equation yields a new equation that includes the gradient of the quantity 
\begin{equation}
\label{tauDef}
{\bm \tau}^{SGS} \approx \overline{H} \mu_{SGS} \nabla\widetilde{\bf u},
\end{equation}
where the $\widetilde{\bf u}$ indicates sub-grid velocity.
The coefficient $\mu_{SGS}$ is defined element-wise and is given as:
\begin{equation}
\label{mun}
\mu_{SGS} = \max\left(0.0, \min (\mu_{\rm max}|_{\Omega_e}, \mu_{{\rm res}}|_{\Omega_e} ) \right),
\end{equation}
where
\begin{equation}
\label{mu1}
\mu_{{\rm res}}|_{\Omega_e} = \overline{\Delta}^2 \, \max \left(\frac{\| R(H)\|_{\infty,\Omega_e}}{\| H - \widehat{H}\|_{\infty,\Omega}} , \frac{\| R(H{\bf u})\|_{\infty,\Omega_e}}{\| H{\bf u} - \widehat{H \bf u} \|_{\infty,\Omega}} \right)
\end{equation}
and
\begin{equation}
\label{mumax}
\mu_{\rm max}|_{\Omega_e} = 0.5\overline{\Delta}\left\Vert  |{\bf u}| + \sqrt{g H_s}  \right\Vert_{\infty,\Omega_e}.
\end{equation}
In (\ref{mu1}, \ref{mumax}), $\ \widehat{\cdot}\ $ indicates the spatially averaged value of the quantity at hand over the global domain $\Omega$,  the norms $\|\cdot \|_{\infty,\Omega}$ at the denominator are used to preserve the physical dimension of the resulting equation, and $R(H, H{\bf u})$ are the residuals of the inviscid governing equations. At each time-step, the residuals are known. 
The presence of $R$ makes the artificial diffusion mathematically consistent\footnote{In the context of stabilization, {\it mathematical consistency} is the property of artificial diffusion such that it vanishes when the residual is zero.}.
The quantity $|{\bf u}| + \sqrt{g H_s}$ in (\ref{mu1}, \ref{mumax}) is the maximum wave speed.

\paragraph{\bf Remark 1: physical dimensions} We would like to emphasize the necessity for the physically correct dimensions of diffusion. This is an important issue that is often underestimated and not accounted for in the design of artificial diffusion methods for stabilization purposes.\\

\section{Wetting and drying interface}
\label{wetDrySECT}
The dry regions are handled by adding an infinitesimal layer of water on the dry surfaces. 
We utilized a threshold value of 1e-3 m\footnote{One millimeter of water is physically negligible.} for all the tests that we ran.
Furthermore, the limiter of Xing et al.\ \cite{xingZhangShu2010} is applied on the velocity and water depth.
Even at high order, this approach is suitable for wetting and drying problems solved via both continuous and discontinuous Galerkin methods.

Albeit simple, the Xing-Zhang-Shu limiter works well. 
However, it is guaranteed to work if and only if the mean water height is positive. 
Xing et al.\ \cite{xingZhangShu2010} recommend using the total variation bounding (TVB) limiter of Shu \cite{shu1987, cockburnShu1989} to make the solution positive definite. 
We rely on artificial dissipation to achieve this as the TVB limiter requires a bound on the derivative of the solution, which we do not know a priori.

\section{Numerical tests}
\label{testsSCT}
We verify the correctness of our numerical model through five standard benchmarks in both one- and two-dimensions.
In 1D, given the N-wave by Carrier et al.\ \cite{carrierEtAl2003} that mimics a wave generated by 
an offshore submarine landslide, we compute the solution of the tsunami approaching a sloping beach. 
See \S \ref{1drunupCarrier}. The numerical solution is evaluated against a set of tabulated surface height and momentum data provided by \cite{tsunamiWorkshopData}. 
The next 1D problem consists of the oscillation of a flat lake in a parabolic bowl \cite{gallardoEtAl2007,xingZhangShu2010}. 
The existence of an analytic solution to this problem allows us to measure the accuracy of the numerical solutions. See \S \ref{gridConvergence}. In addition, various 1D wetting and drying test cases were analyzed in the thesis of Hood \cite{karolineTHESIS2016} to assess our model.

In 2D, we analyzed the results for the following tests. The vertical oscillation of a column of water in an axisymmetric paraboloid
with analytic solution \cite{thacker1981}. We describe this problem in \S \ref{2dparaboloidTest} to test the inundation against the two-dimensional effects of a varying topography.
The second two-dimensional test involves the flooding in a closed channel as presented by Gallardo et al.\ \cite{gallardoEtAl2007} and by Xing and Zhang \cite{xingZhang2013}. No exact solution is presented for this test, but the results can be compared against previous studies. This is shown in  \S \ref{3islands}.
We conclude the set of 2D problems with the simulation of the flow in a two-dimensional channel with an isolated obstacle \cite{synolakis1987}. We use this test to gain insight into the practical -- other than theoretical as stated in \cite{gerbeauEtAl2001} -- need to include viscosity by analyzing the energetics of the solutions. The test is run during a sufficiently long time to assess the robustness of the inundation and viscosity schemes to handle fast motion with interacting waves that impinge against a steep obstacle. See \S \ref{OneislandTest}. 
The description and analysis of more tests can be found in the collection recently compiled by Delestre et al.\ \cite{delestreEtAl2014}.

\subsection{1D tsunami runup over a sloping beach}
\label{1drunupCarrier}
\begin{figure}
\centering
\includegraphics[width=0.7\textwidth]{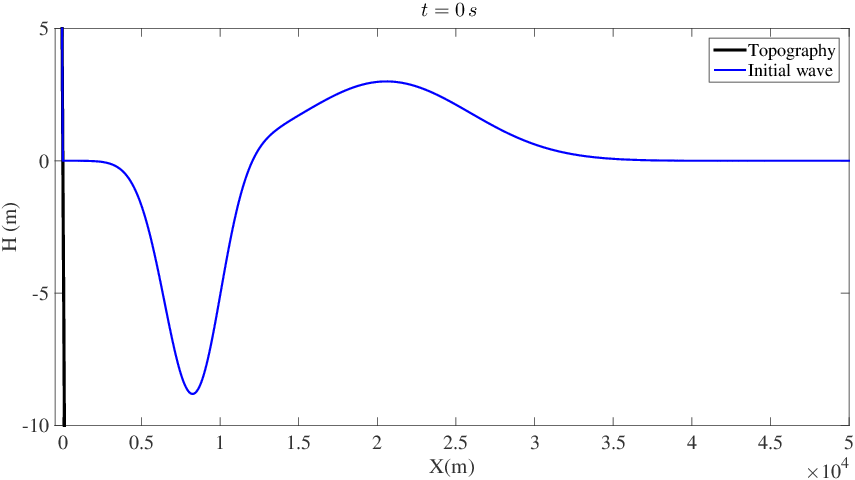}
\caption{1D tsunami runup over a sloping beach: initial N-wave.}
\label{runuptinitialwave}
\end{figure}

The runup of a long wave on a uniformly sloping beach was originally proposed at the third international workshop on long-wave runup models \cite{tsunamiWorkshopData}. 
The one-dimensional computational domain is defined as $\Omega = x = [-500, 50000]$ m. The dry initial beach is 500 m long.
The initial waveform was defined by Carrier et al.\ in \cite{carrierEtAl2003} for an $L=8$ m domain as:

\begin{equation}
\label{carrierWaveEQN}
\eta = a_1 \exp\{ -\hat{k}_1(x - \hat{x}_1)^2\} - a_2\exp\{ \hat{k}_2(x - \hat{x}_2)^2 \},
\end{equation}
with constants $(a_1,a_2,\hat{k}_1,\hat{k}_2, \hat{x}_1, \hat{x}_2) = (0.006, 0.018, 0.4444, 4.0, 4.1209, 1.6384)$.
To scale the wave to the $L=50000$ m long domain used for the current test, the scaling factor 
$\delta = L/8$ is introduced and Eq.\ (\ref{carrierWaveEQN}) is re-expressed with respect to $x_{1,2} = \hat{x}_{1,2} \delta$ and $k_{1,2} = \hat{k}_{1,2}/\delta^2$, with larger amplitudes $(a_1, a_2)=(3.0, -8.8)$. The initial wave is plotted in Fig.\  \ref{runuptinitialwave}.

\begin{figure}
\centering

\includegraphics[width=0.495\textwidth]{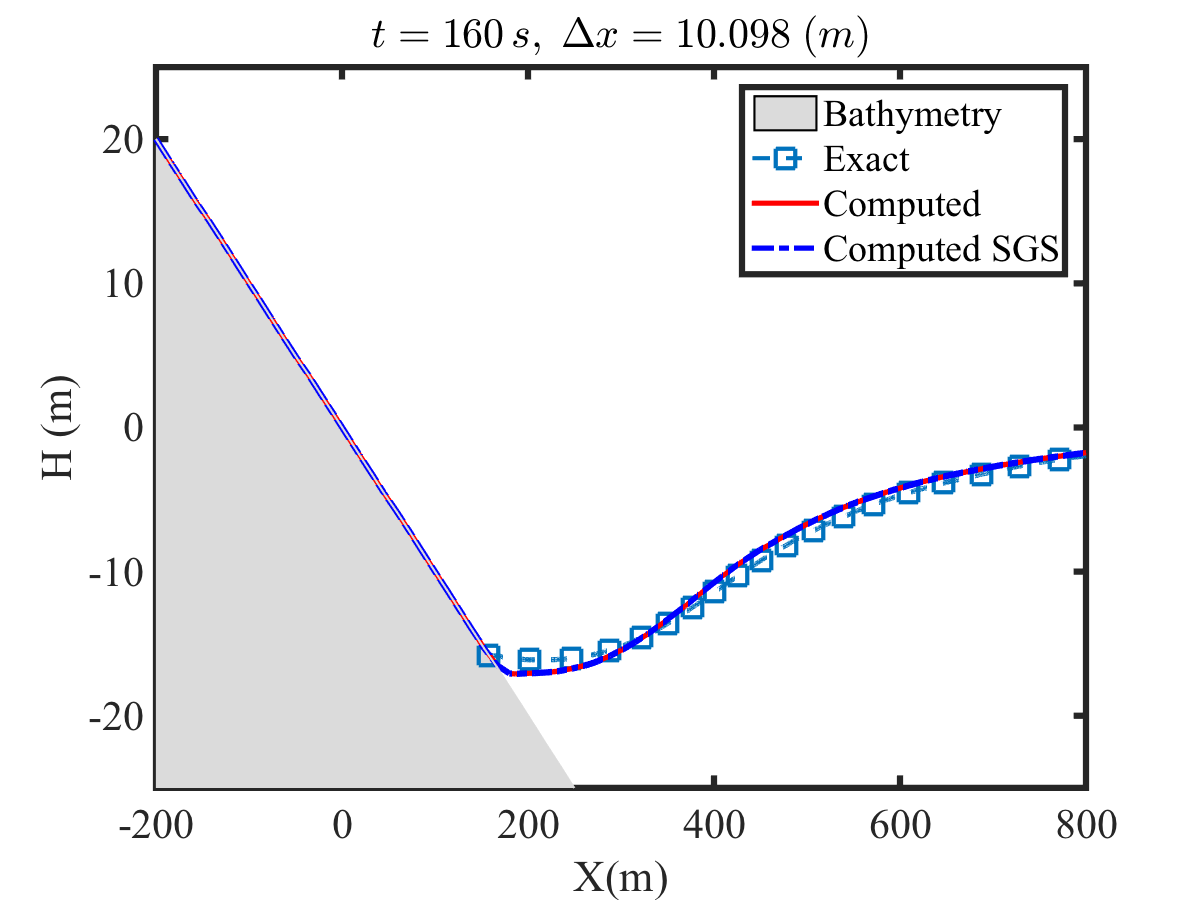}
\includegraphics[width=0.495\textwidth]{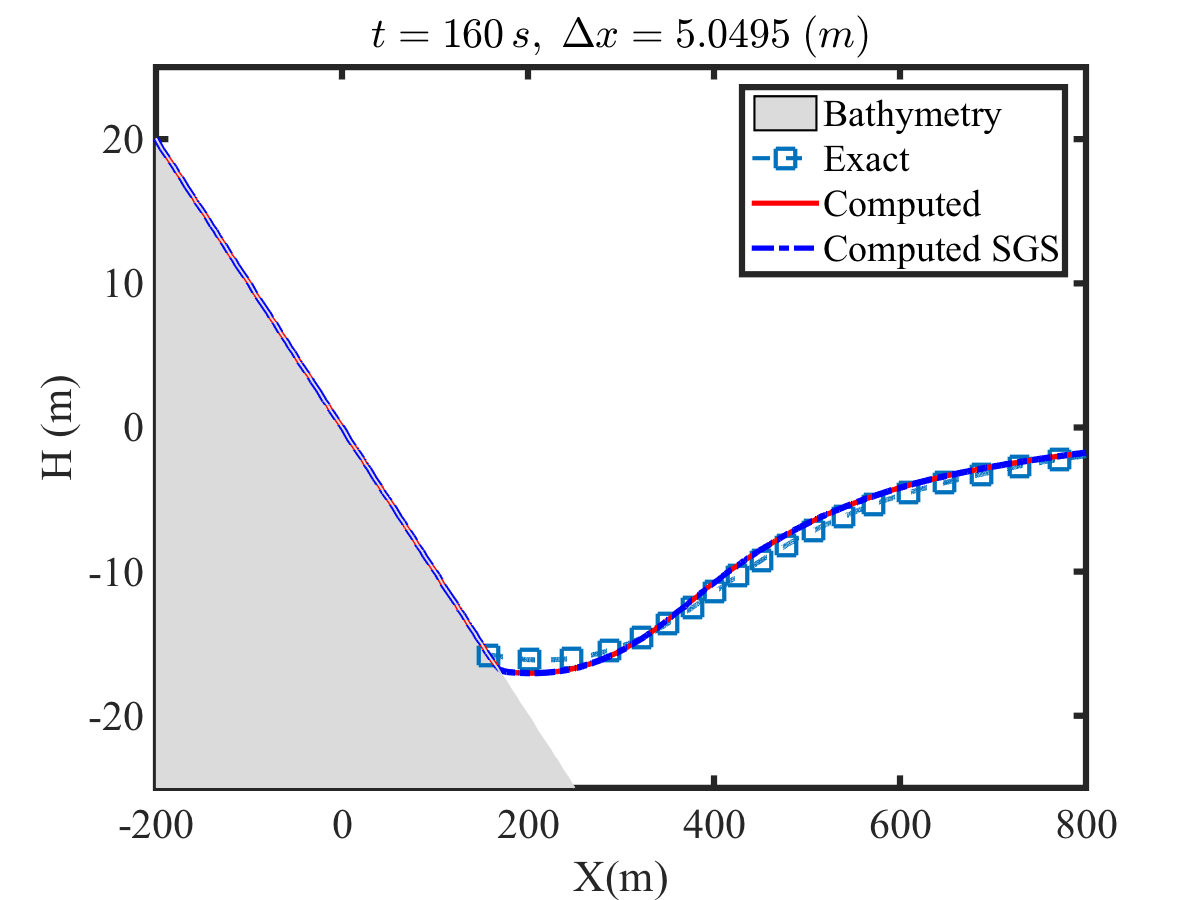}\\

\includegraphics[width=0.495\textwidth]{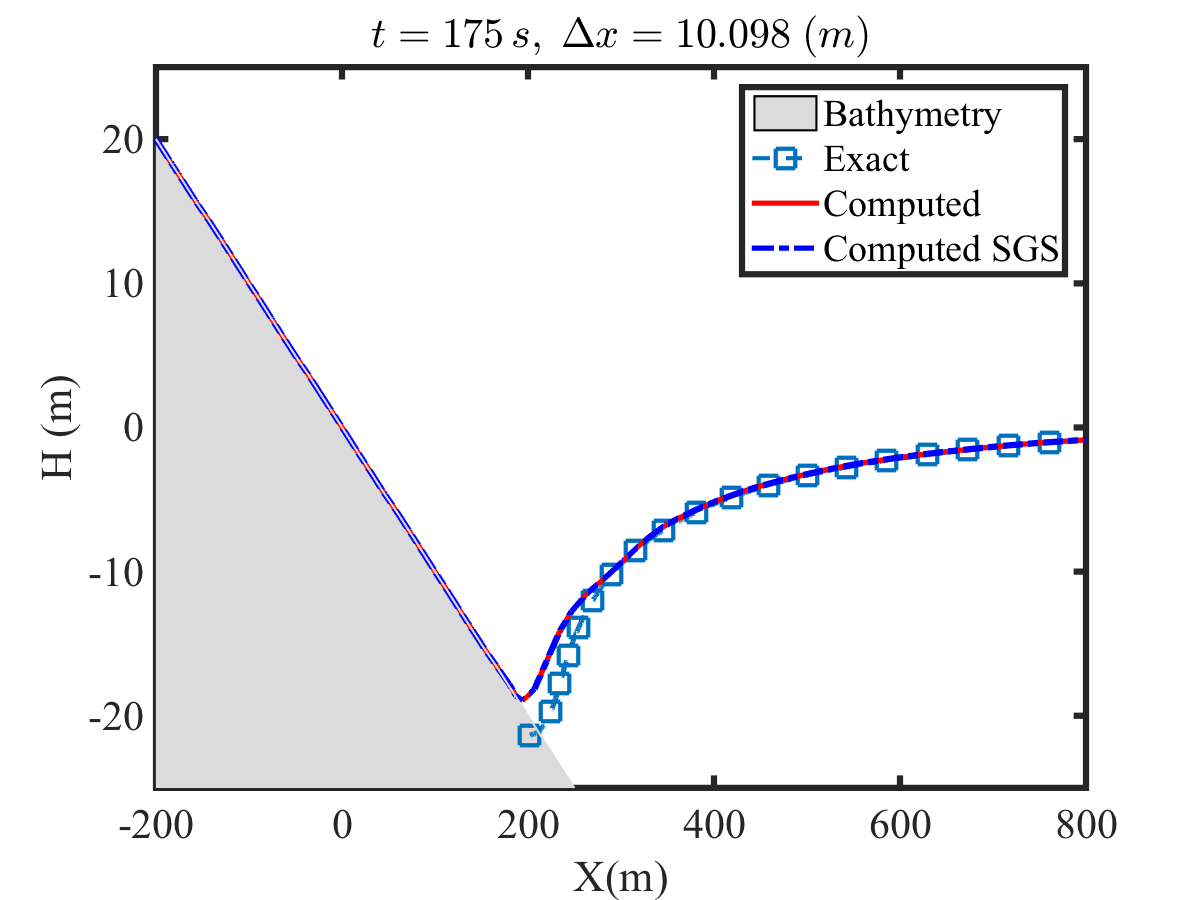}
\includegraphics[width=0.495\textwidth]{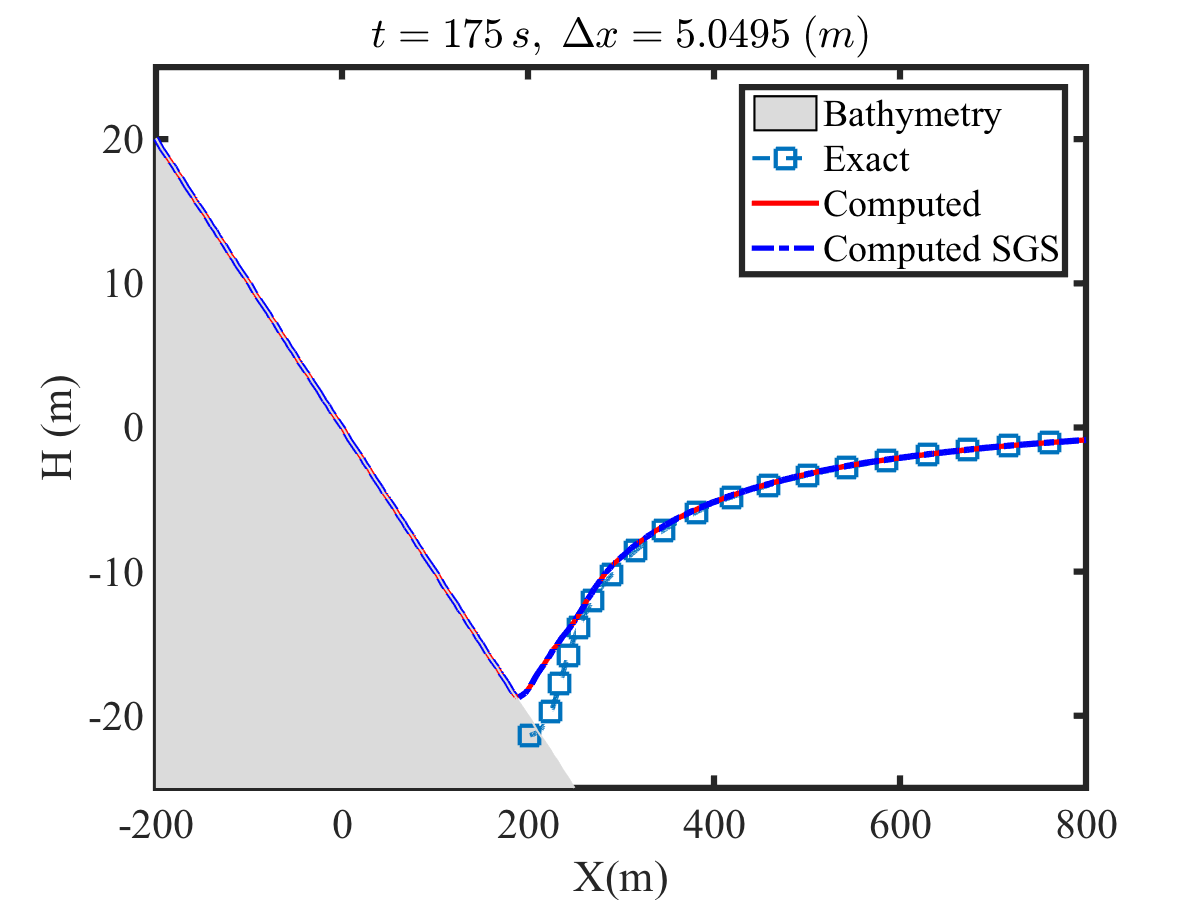}\\

\includegraphics[width=0.495\textwidth]{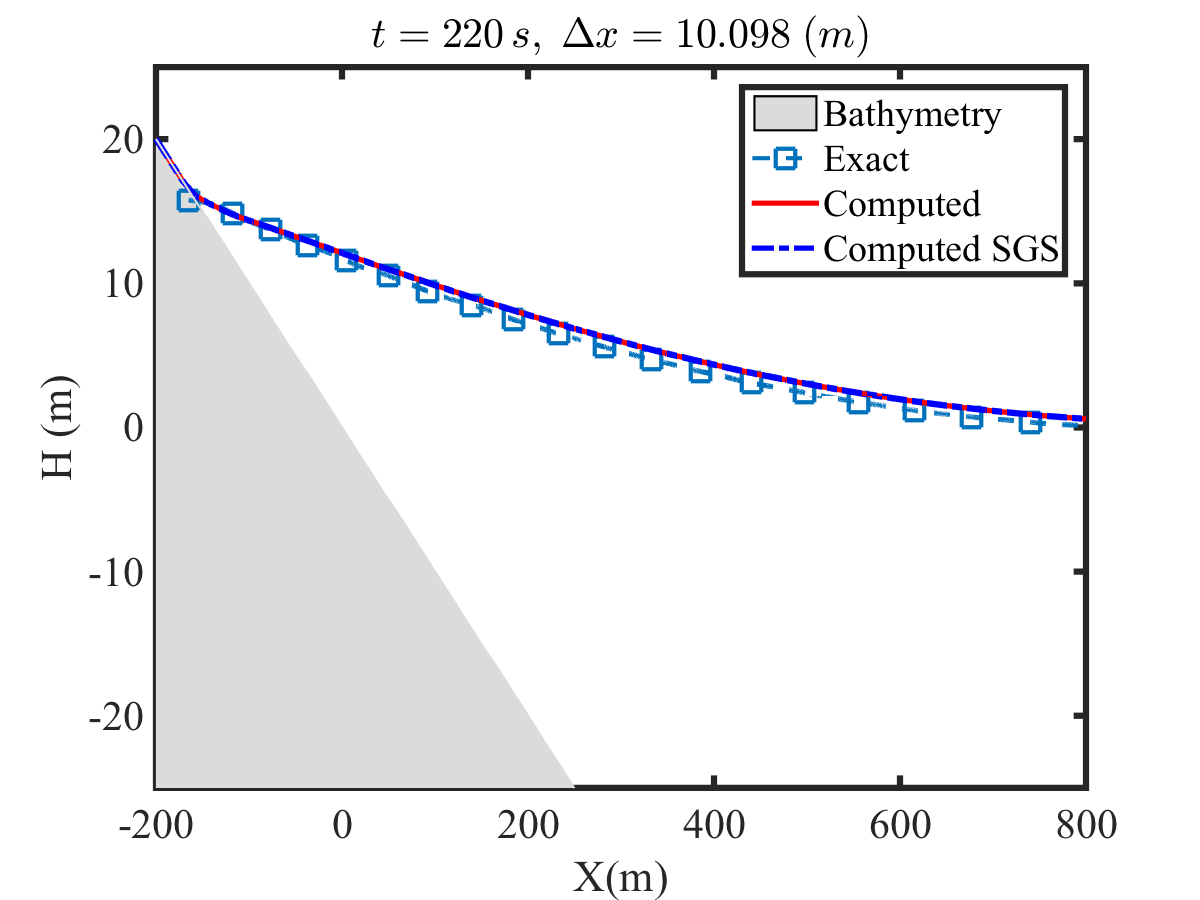}
\includegraphics[width=0.495\textwidth]{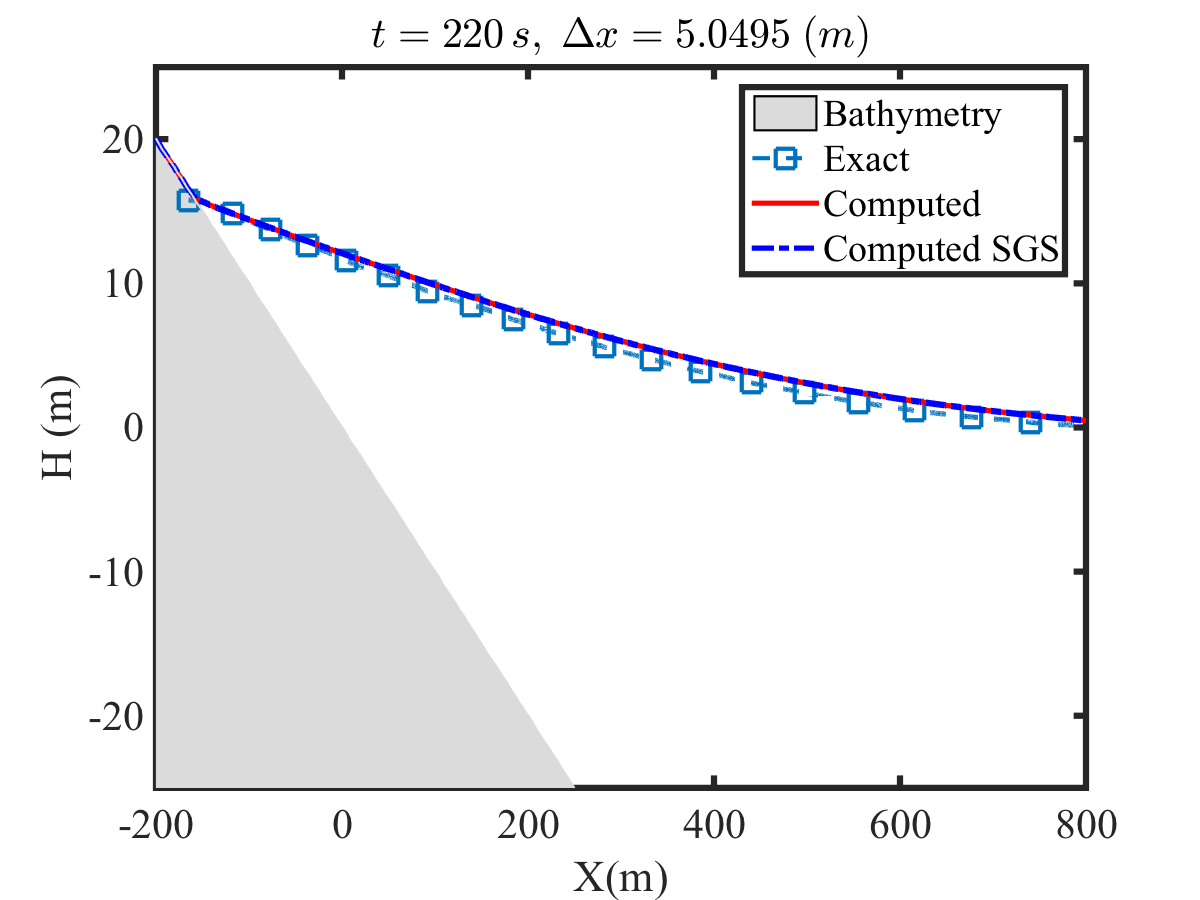}

\caption{1D tsunami runup over a sloping beach. 
Computed and exact solutions at $t=[160,175,220]$ s. Left column: $\Delta x \approx 10$ m (1250 elements of order $4$). 
Right column: $\Delta x \approx 5$ m (2500 elements of order $4$). 
The computed inviscid and viscous (SGS) solutions appear almost perfectly superimposed. 
The problem is smooth almost everywhere so that the intensity of the dynamic dissipation is minimal. 
This translates into an almost imperceptible effect of diffusion on the solution. The structure of the dissipation coefficient is plotted in Fig.\ \ref{runup1D_SGS}.}
\label{runupSolution220s_2500el}
\end{figure}

\begin{figure}
\centering
\includegraphics[width=0.7\textwidth]{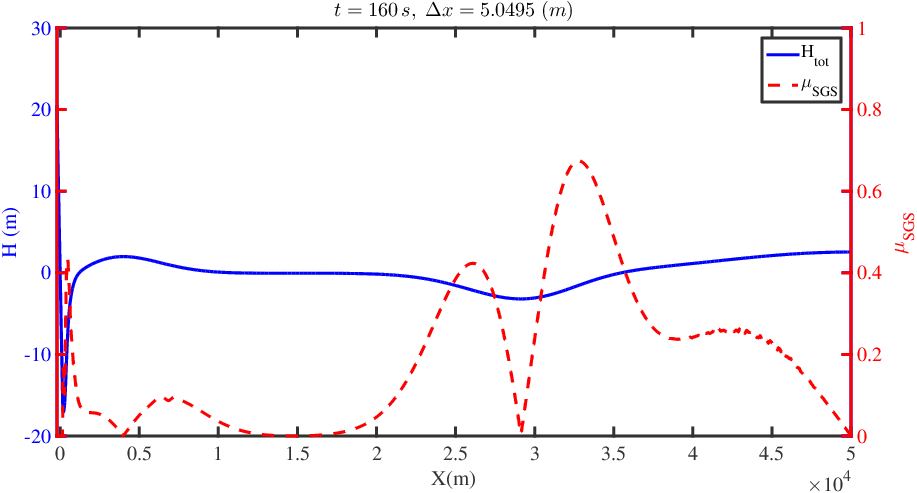} 
\includegraphics[width=0.7\textwidth]{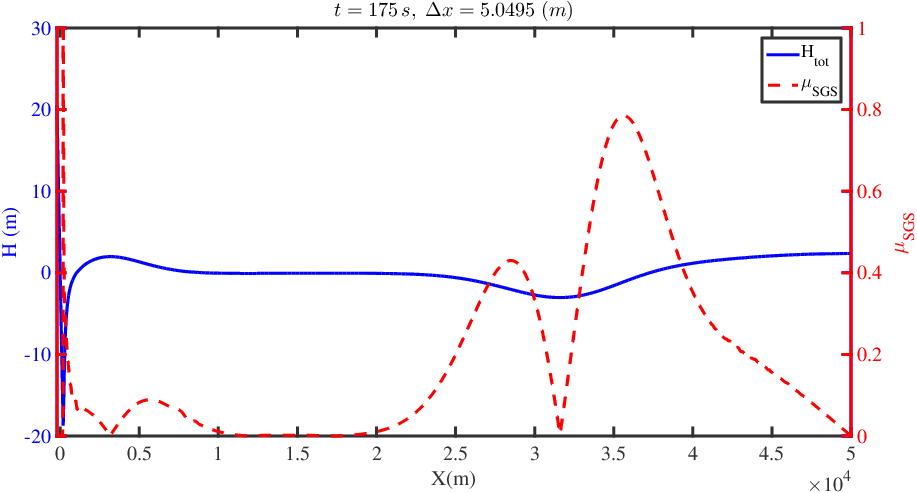} 
\includegraphics[width=0.7\textwidth]{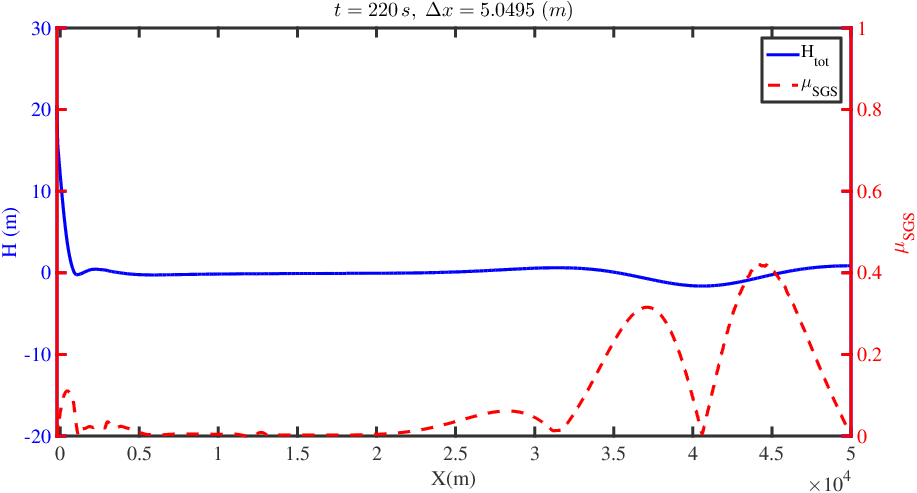} 
\caption{1D tsunami runup over a sloping beach. Dynamic $\mu_{SGS}$ (red, dashed line) and water surface (blue, solid line) in the full domain. The effect of diffusion on the solution of Fig.\ \ref{runupSolution220s_2500el} is minimal as the value of the coefficient is indeed very small. The solution is smooth almost everywhere, which is the reason for the very small values of the dynamic diffusion coefficient.}
\label{runup1D_SGS}
\end{figure}

The solutions at times $t=[0,\, 110,\, 220]$ s are plotted in Fig.\ \ref{runupSolution220s_2500el} and are 
compared against the tabulated data available in \cite{tsunamiWorkshopData}.
Fig.\ \ref{runupSolution220s_2500el} shows that the effect of diffusion on the water surface solution is clearly negligible. This can be explained by looking at 
the structure and values of $\mu_{SGS}$ in Fig.\  \ref{runup1D_SGS}. With a water surface that is smooth almost everywhere, 
the dynamic diffusion coefficient is so small that its effect becomes minimal. 
We will see later that this will not be the case in problems with a greater degree of irregularity of the surface.

\begin{figure}
\centering
\includegraphics[width=0.49\textwidth]{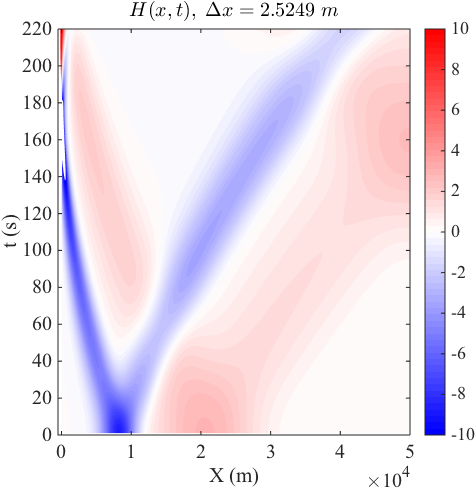} 
\includegraphics[width=0.49\textwidth]{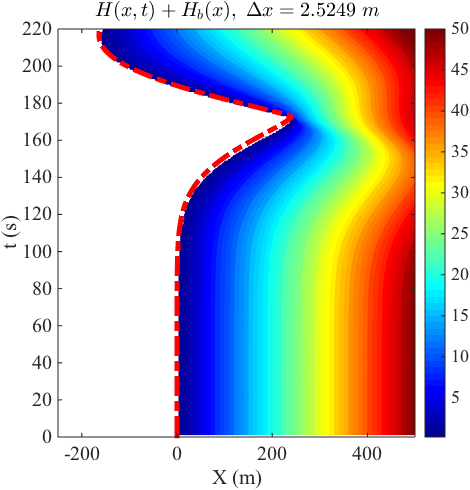} 
\caption{1D tsunami runup over a sloping beach. Left: wave trajectory in the full 50 km long domain. Right: $x-t$ variation of the water depth in the proximity of the coast.
The shoreline is at the interface between the white area (dry shore) and the color shading (water surface). The dashed red curve is the exact shoreline.}
\label{runupXTspaceFullDomain}
\end{figure}

We plot the space-time evolution of the shoreline in Fig.\ \ref{runupXTspaceFullDomain} where the wave elevation and total depth are plotted in the proximity of the shore. 
The dashed red curve represents the tabulated shoreline \cite{tsunamiWorkshopData}. 
By direct comparison with Carrier's results \cite{carrierEtAl2003}, the patterns of the water surface elevation ($\eta({x,t})$) and total water depth 
are in full agreement.

\subsection{Grid convergence rate}
\label{gridConvergence}
\begin{figure}
\centering
\includegraphics[width=0.49\textwidth]{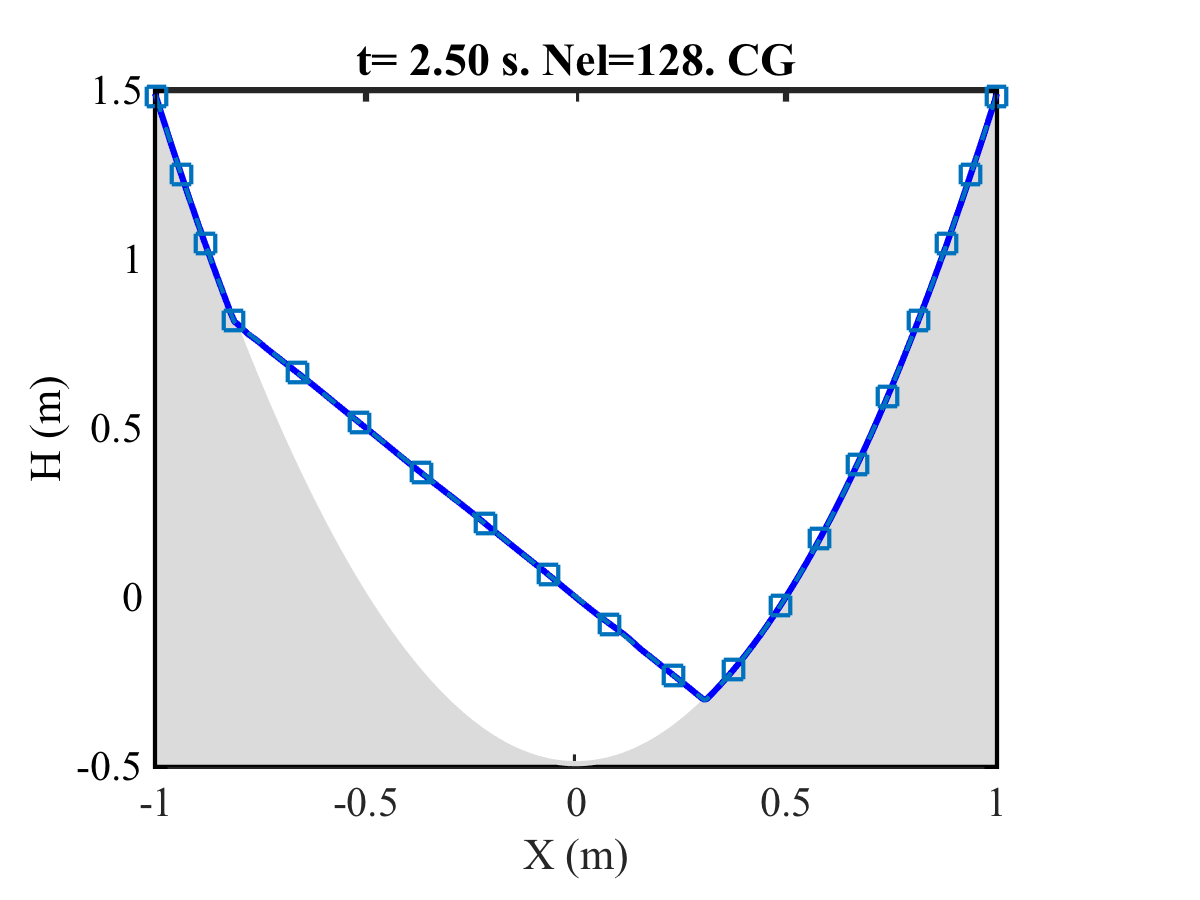}
\includegraphics[width=0.49\textwidth]{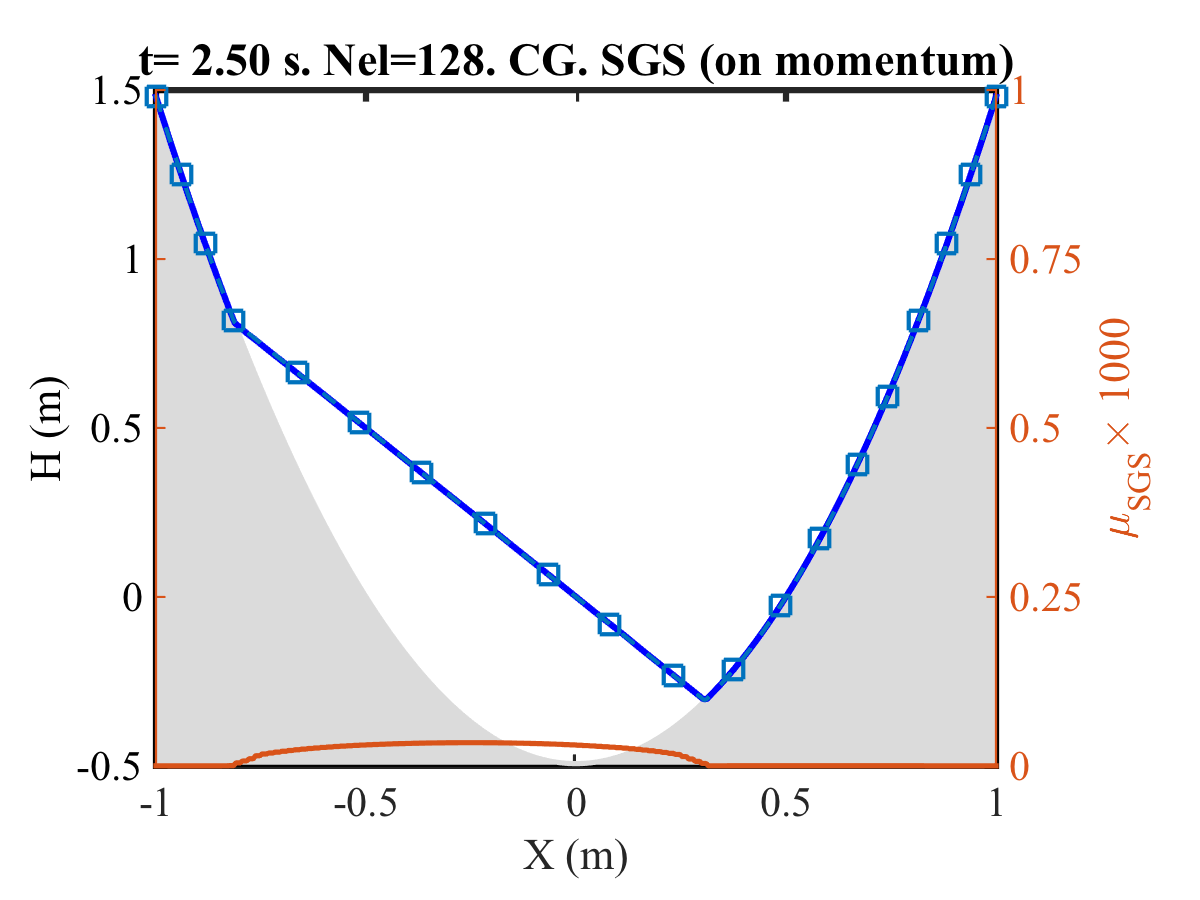}\\

\includegraphics[width=0.49\textwidth]{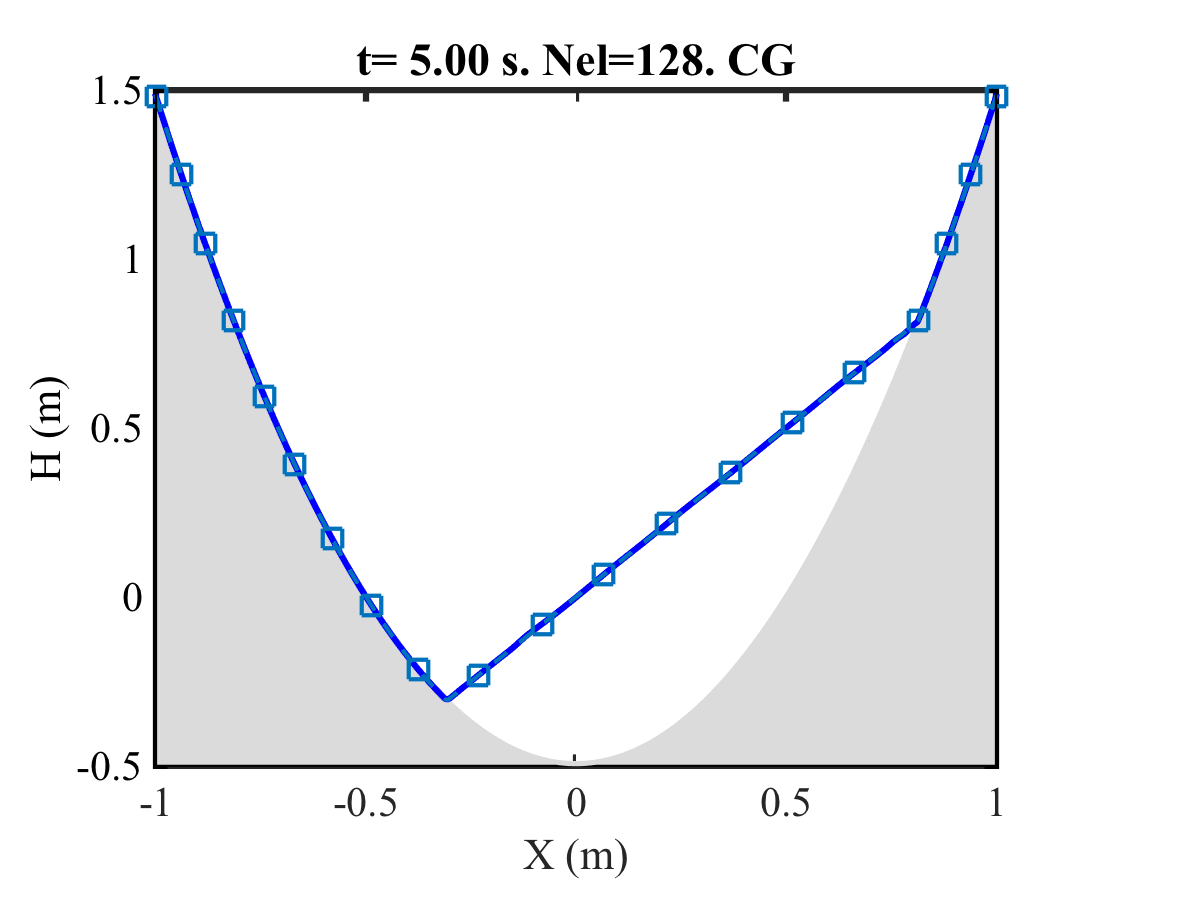}
\includegraphics[width=0.49\textwidth]{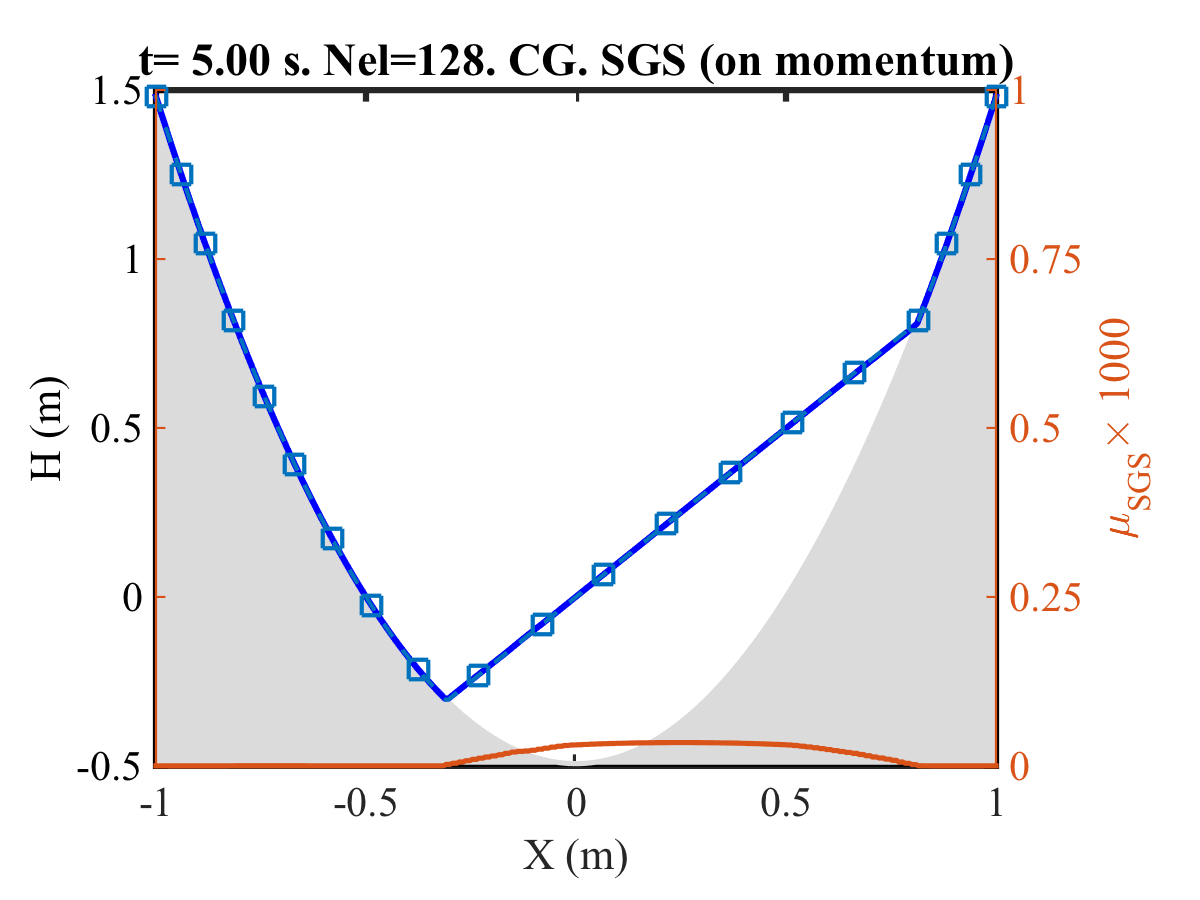}\\

\includegraphics[width=0.49\textwidth]{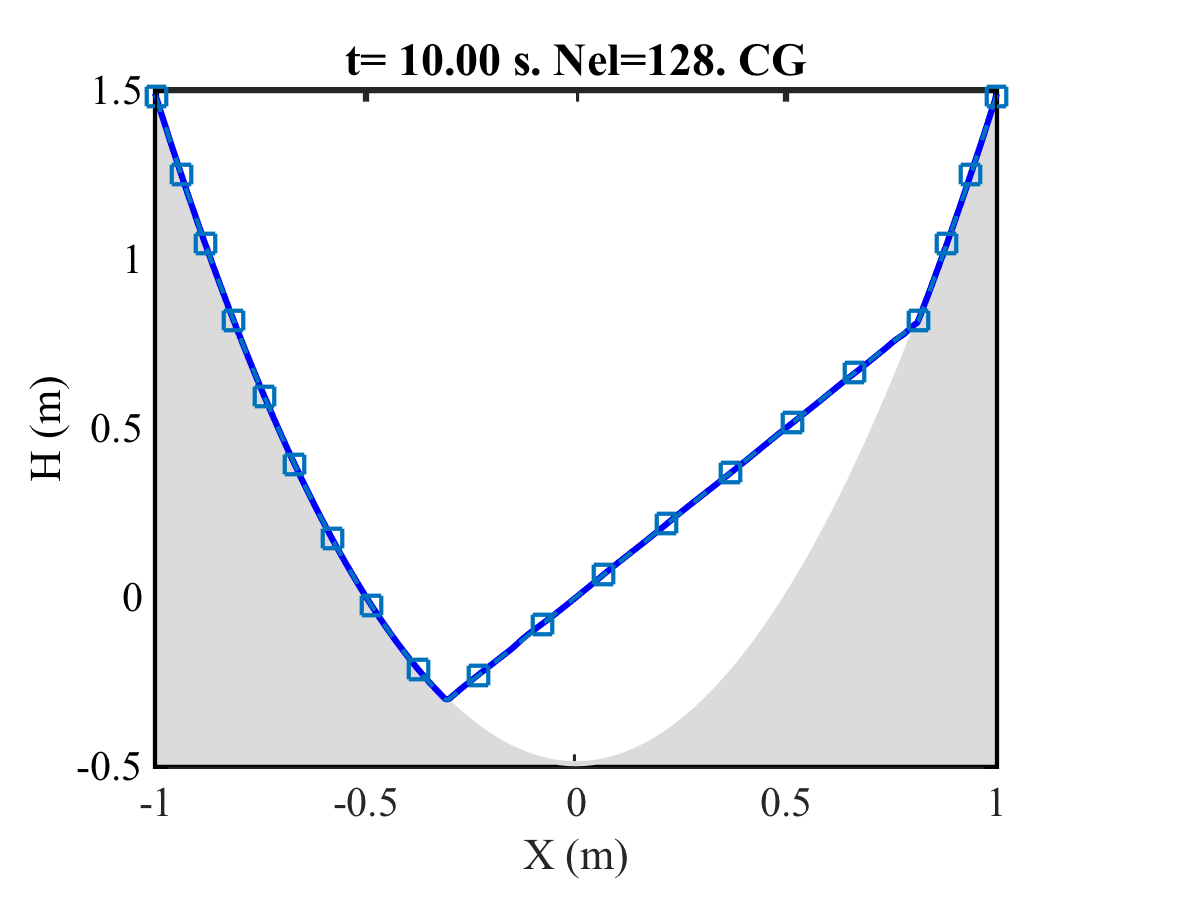}
\includegraphics[width=0.49\textwidth]{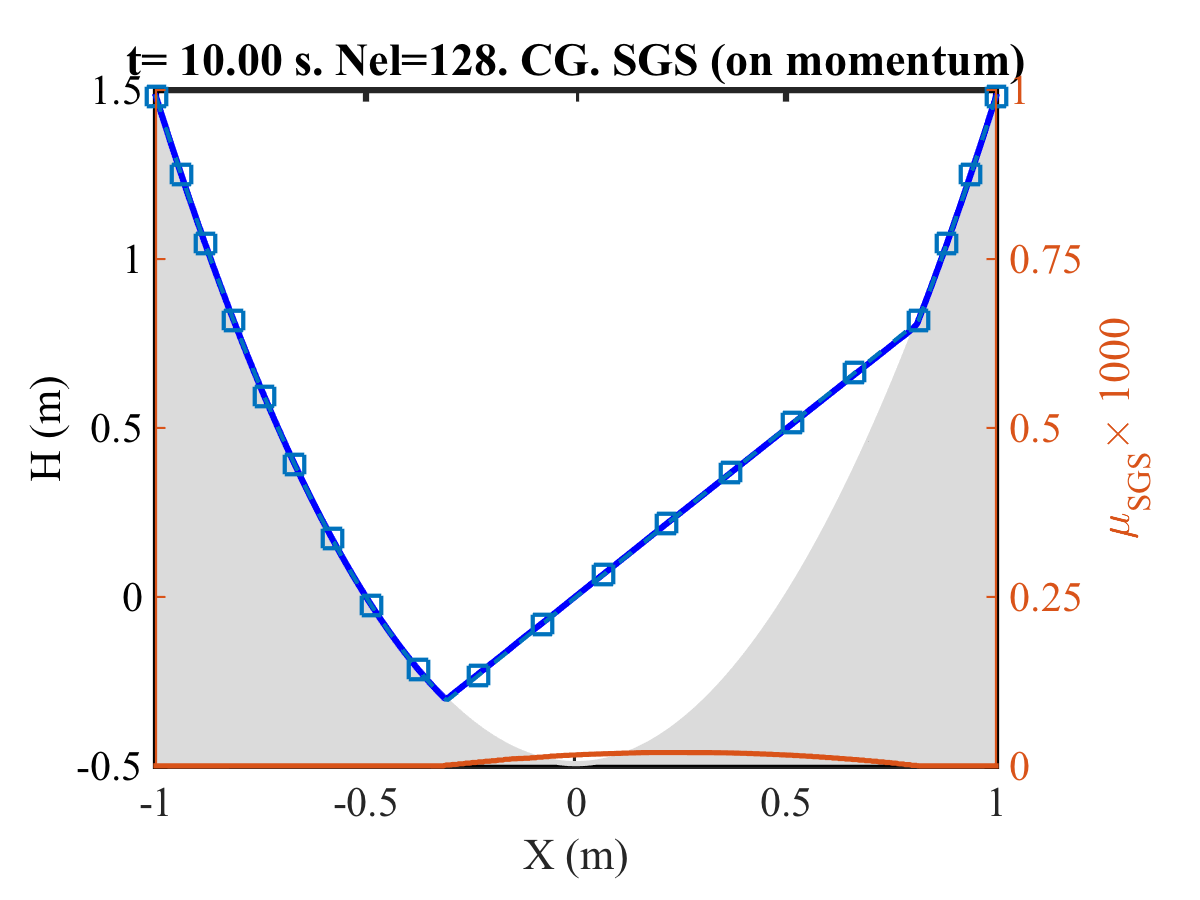}\\
\caption{1D flow in a parabolic bowl. {\bf CG} solution. Left column: inviscid. Right column: viscous. Computed water level (blue, solid line), exact solution (blue, dashed line with open squares), and $\mu_{SGS}$ (red, solid line) at different times using 128 elements of order 4. For visualization,  $\mu_{SGS}$ is scaled by a factor of 1000. The solution is sufficiently smooth that the dynamic viscosity is negligible. On the other hand, the limiter$+$wetting/drying scheme is responsible for the positivity of the water surface in the proximity of the moving boundary.}
\label{1DparabolicBowl128CG}
\end{figure}

\begin{figure}
\centering

\includegraphics[width=0.49\textwidth]{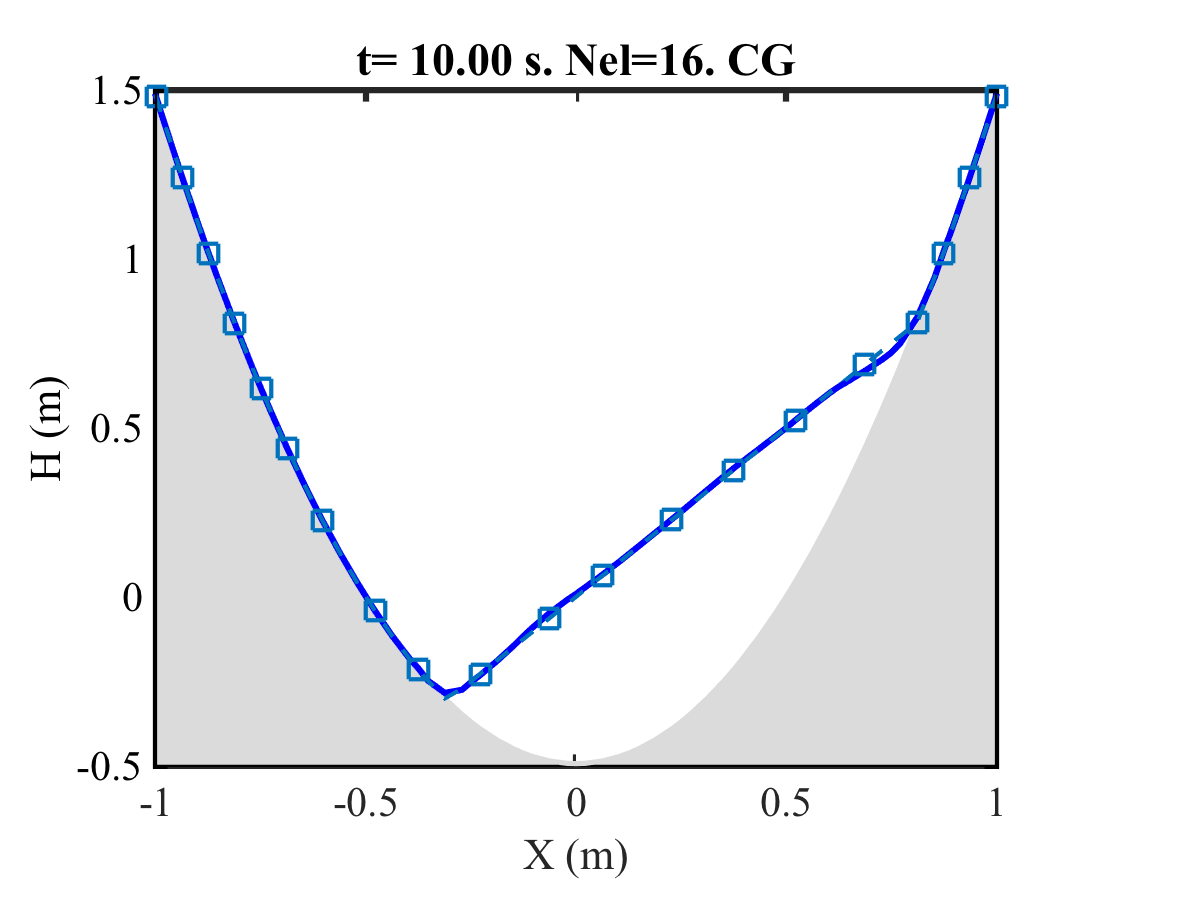}
\includegraphics[width=0.49\textwidth]{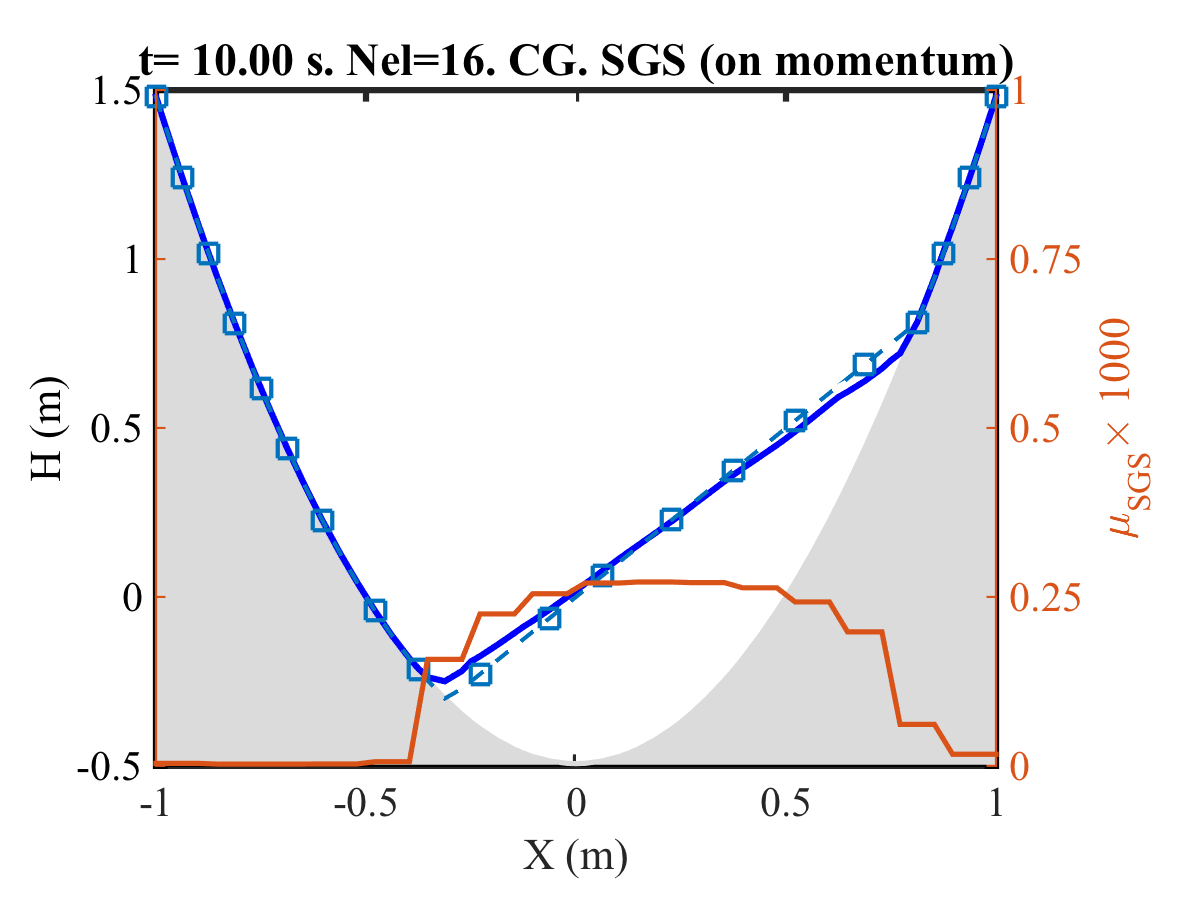}

\caption{1D flow in a parabolic bowl. Low resolution {\bf CG} solutions using 16 elements of order 4. 
Left column: inviscid. Right column: viscous. Computed water level (blue, solid line), exact solution (blue, dashed line with open squares), and $\mu_{SGS}$ (red, solid line) at $t=10$ s. For visualization, $\mu_{SGS}$ is scaled by a factor of 1000.}
\label{1DparabolicBowlCGlowRes}
\end{figure}

\begin{figure}
\centering
\includegraphics[width=0.49\textwidth]{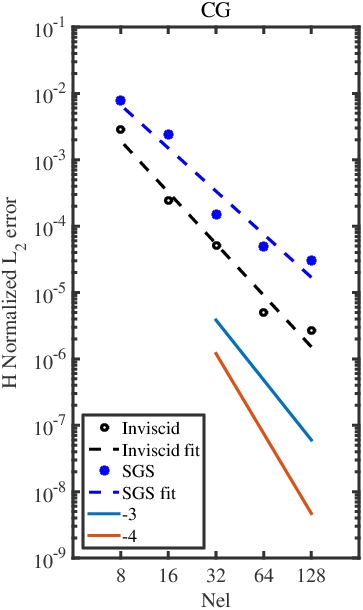}
\includegraphics[width=0.49\textwidth]{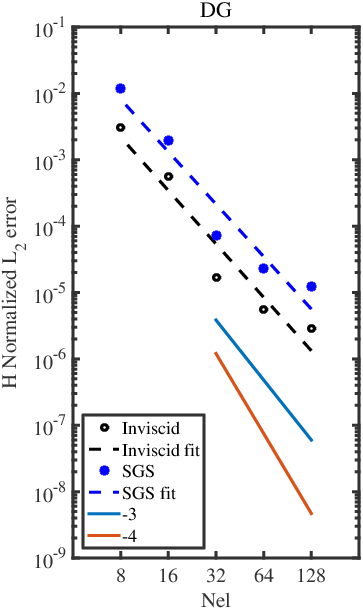}
\caption{Normalized $L_2$ error of water surface at $t=10$ s for {\bf CG} (left) and {\bf DG} (right). 
The $-3$ and $-4$ curves indicate the reference rates.
}
\label{L2error1DparabolicBowl}
\end{figure}

\begin{figure}
\centering
\includegraphics[width=0.49\textwidth]{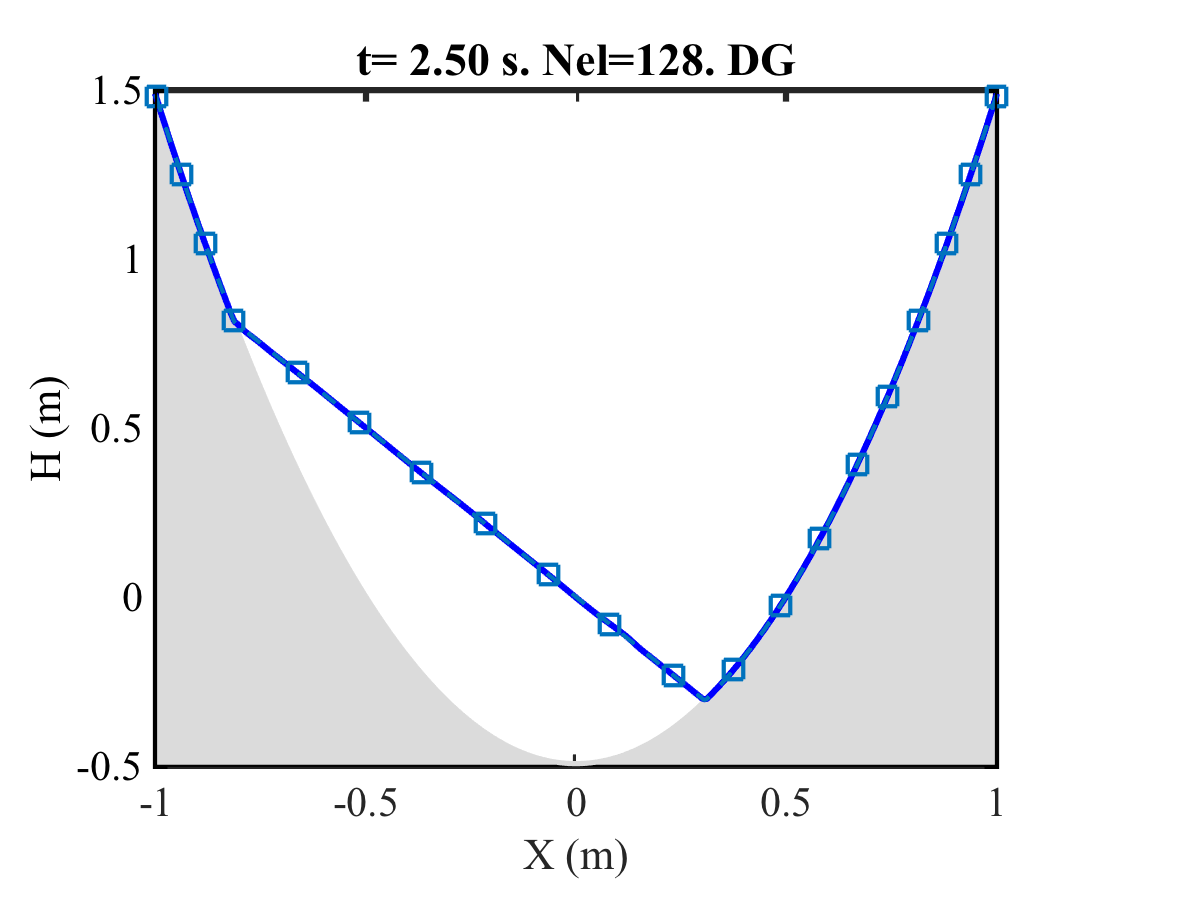}
\includegraphics[width=0.49\textwidth]{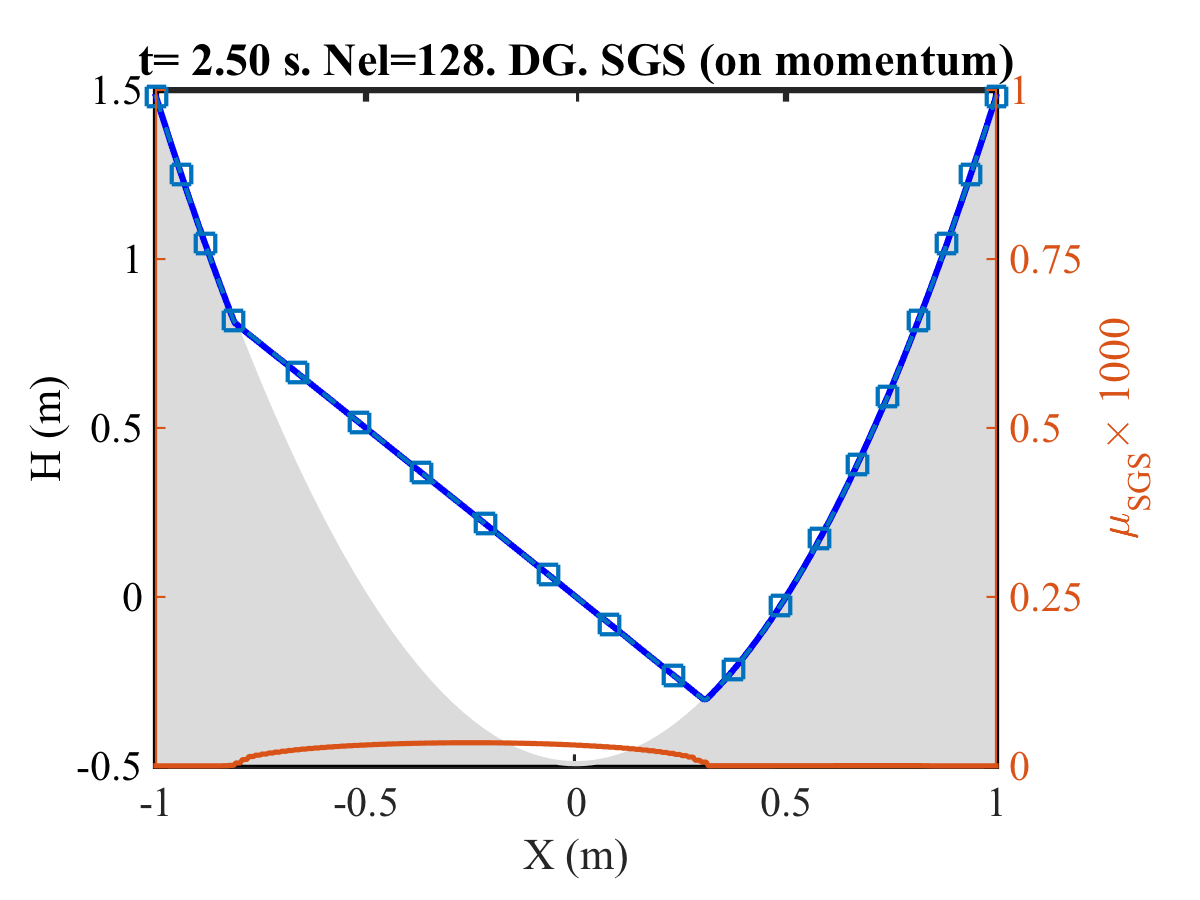}\\

\includegraphics[width=0.49\textwidth]{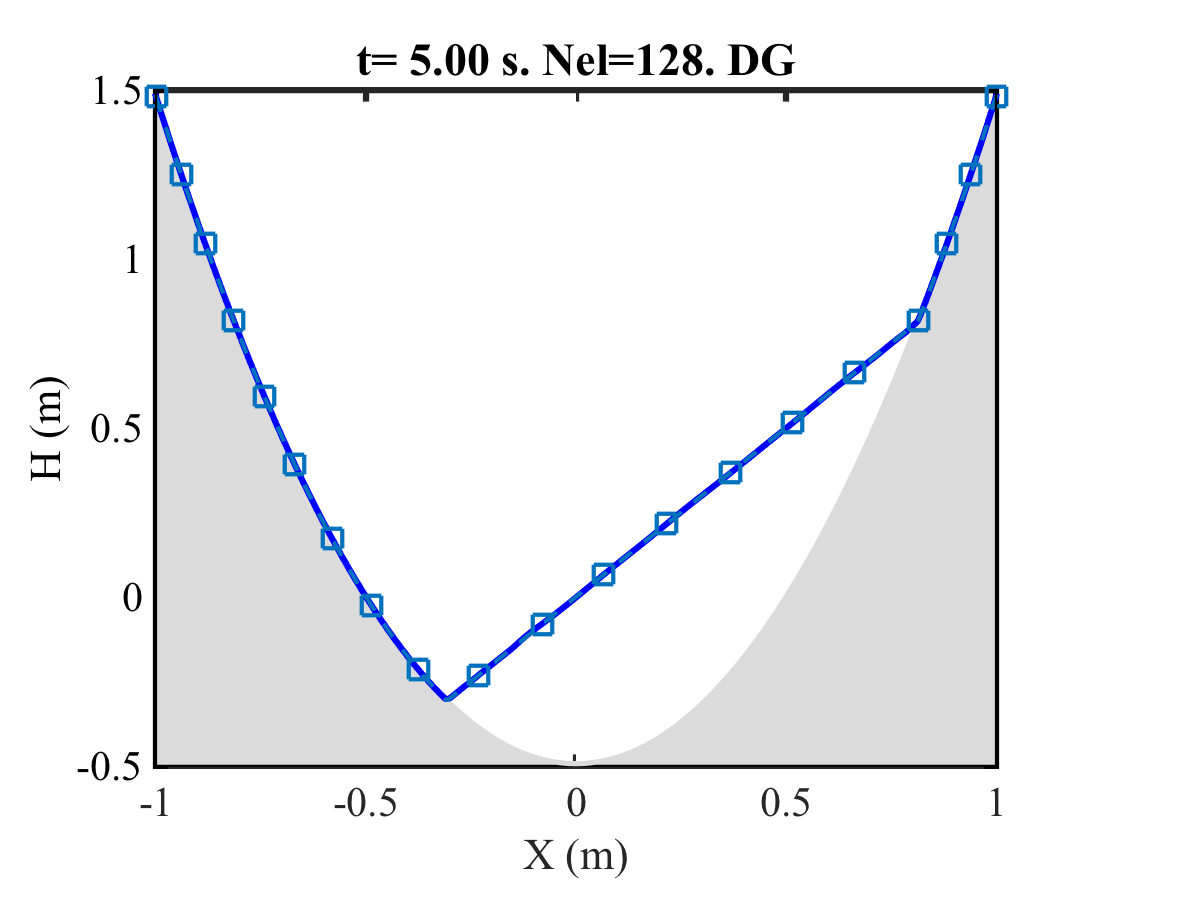}
\includegraphics[width=0.49\textwidth]{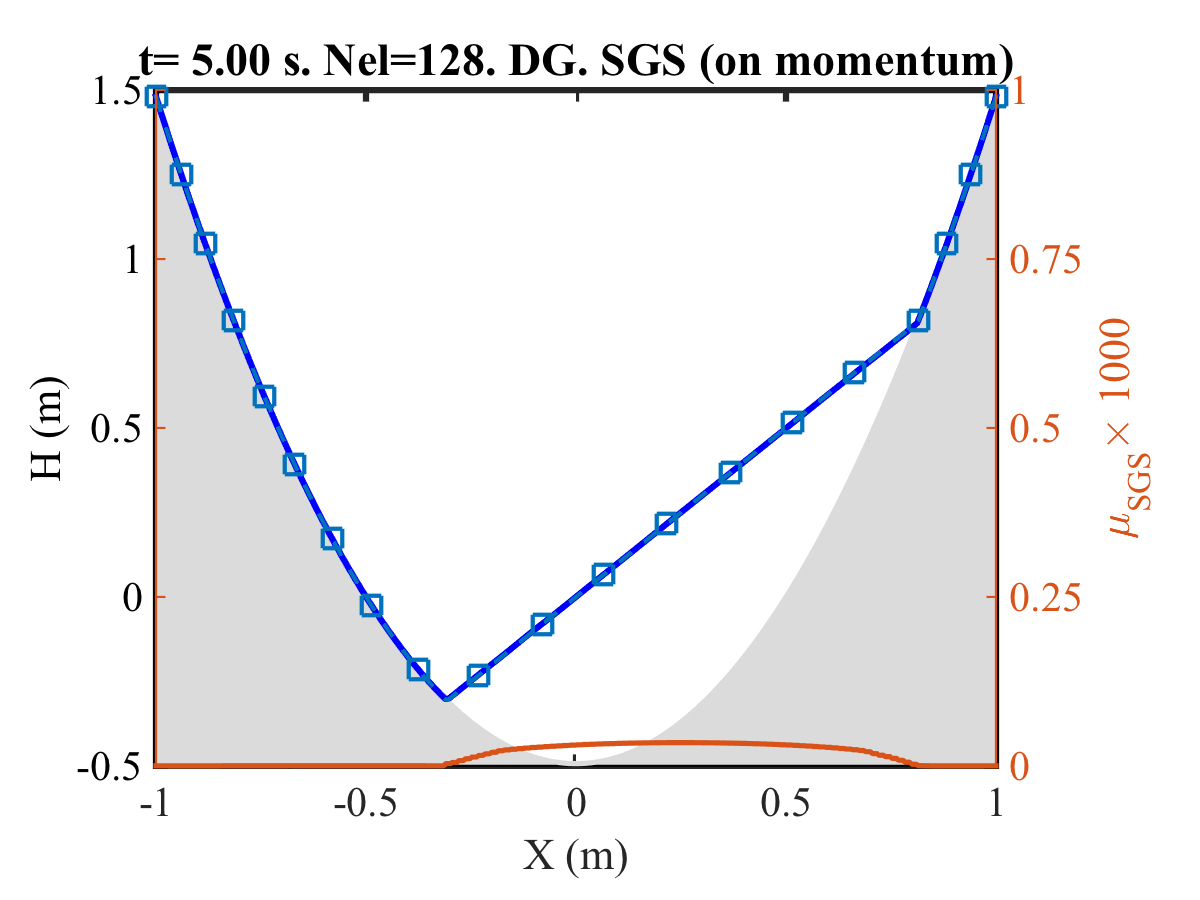}\\

\includegraphics[width=0.49\textwidth]{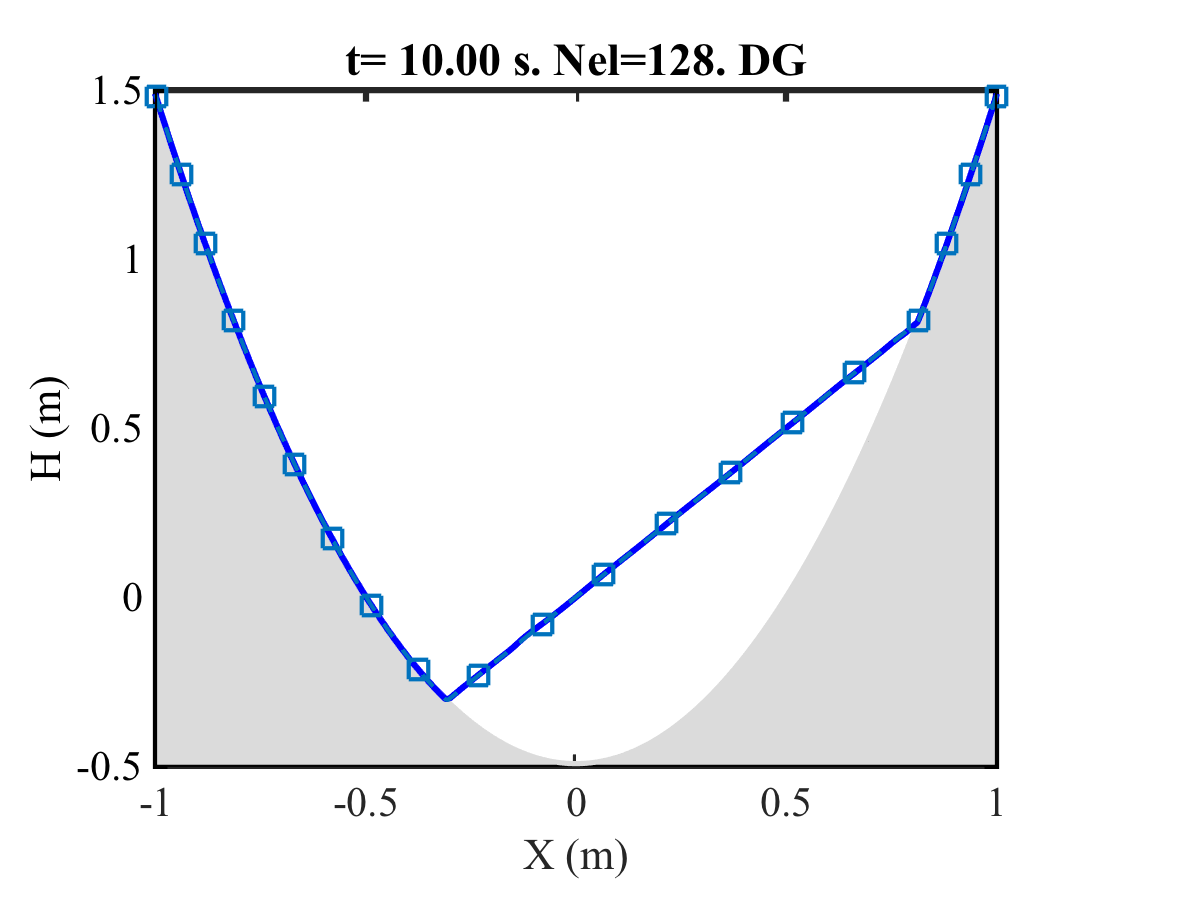}
\includegraphics[width=0.49\textwidth]{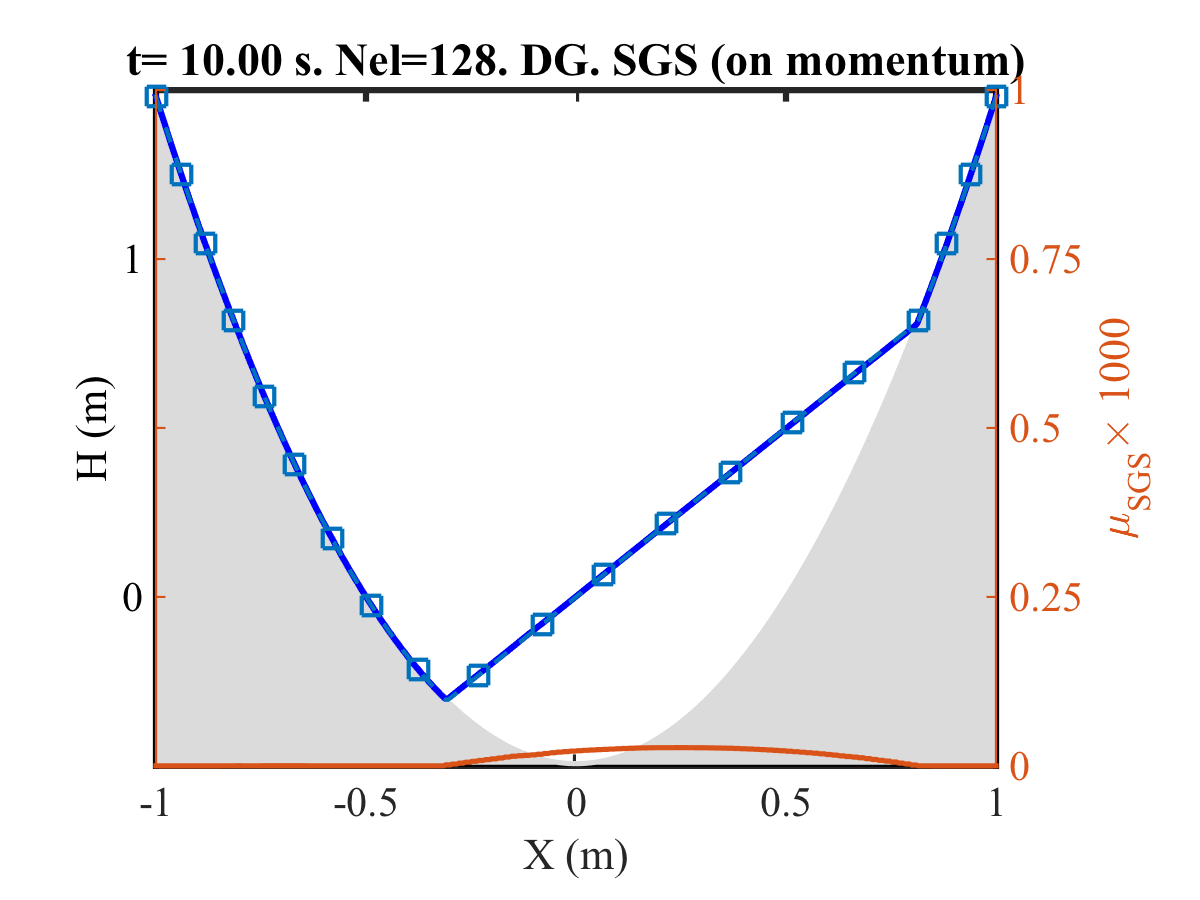}\\
\caption{1D flow in a parabolic bowl. Like Fig.\ \ref{1DparabolicBowl128CG} but using {\bf DG}.}
\label{1DparabolicBowl128DG}
\end{figure}

\begin{figure}
\centering

\includegraphics[width=0.49\textwidth]{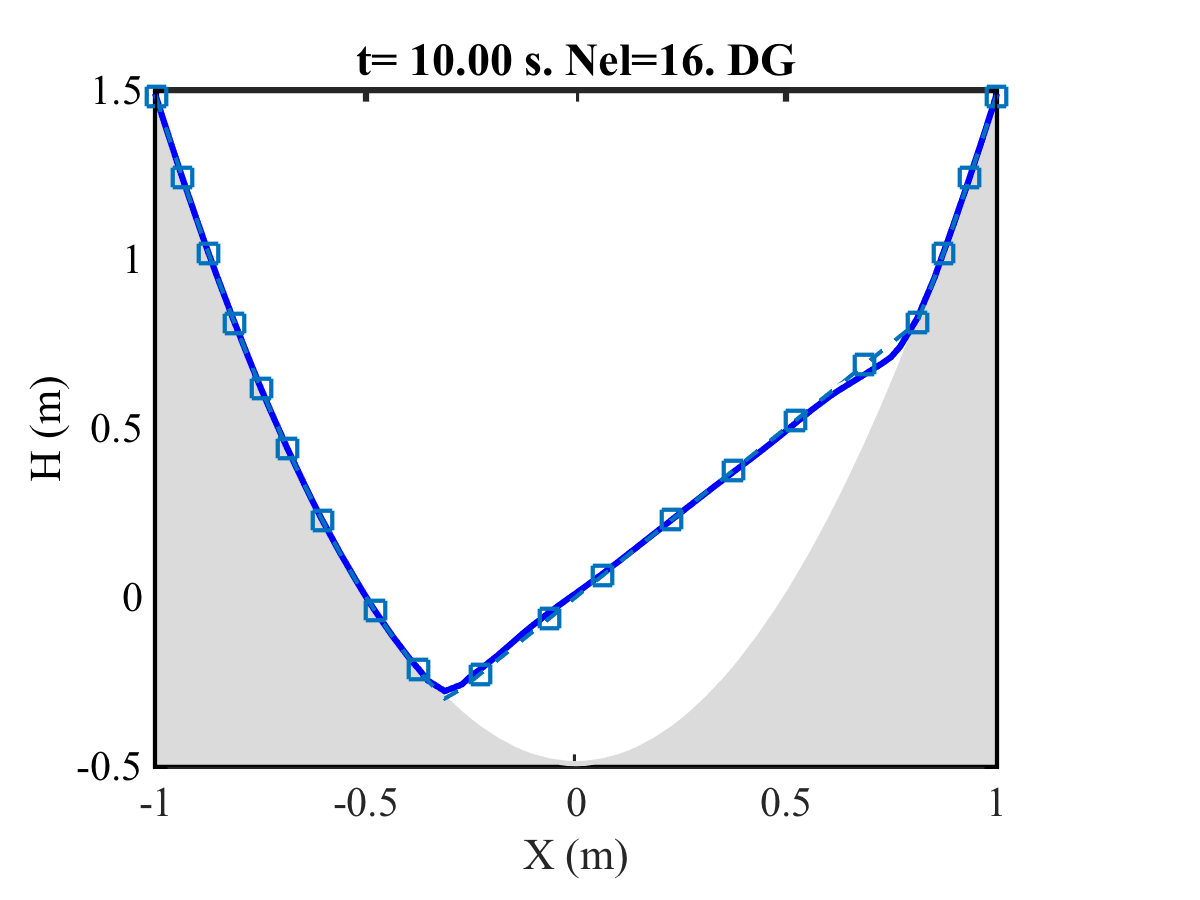}
\includegraphics[width=0.49\textwidth]{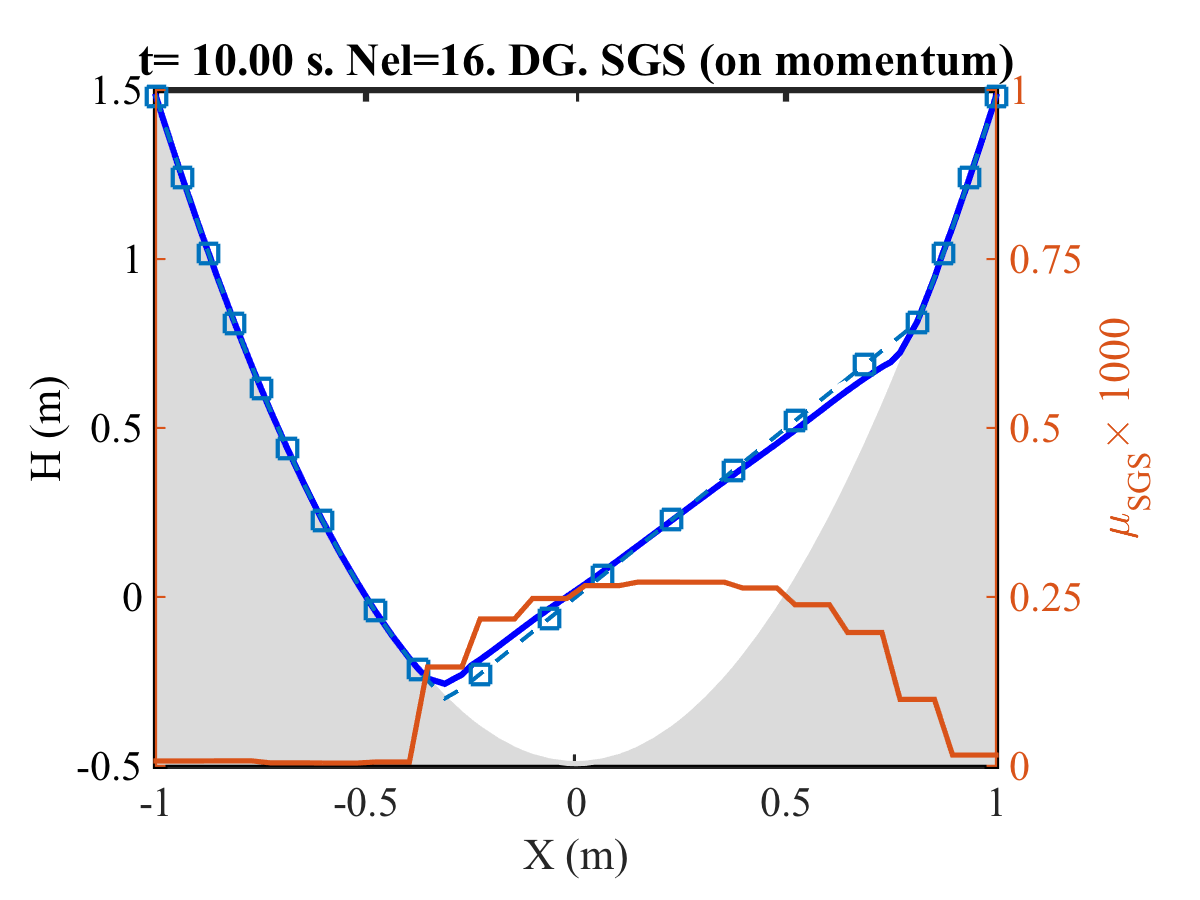}\\

\caption{Like Fig.\ \ref{1DparabolicBowlCGlowRes}, but using {\bf DG}.}
\label{1DparabolicBowlDGlowRes}
\end{figure}

To measure the convergence rate of the model, we compare the computed solutions against the analytic solution of the flow in a one-dimensional parabolic bowl \cite{gallardoEtAl2007,thacker1981}. 
The parabolic topography is defined as:
\[
H_b(x) = h_0\left(\frac{1}{a^2} x^2\right) - 0.5
\]
where $h_0=2$ m and $a=1$ in $\Omega=x=[-1,1]$ m. The initial velocity is $u=0\, {\rm m s^{-1}}$ and the water surface begins to oscillate due to gravity only.
The solution is computed using 
8, 16, 32, 64, and 128 elements of order 4. The CG solution using 128 elements is plotted in Fig.\ \ref{1DparabolicBowl128CG}, where we compare the inviscid solution against the viscous computation at $t=[2.5,\,5,\,10]$ s. 
The solution preserves its smoothness at all times as long as the resolution is sufficiently fine. 
At high resolution, it is evident that the CG approximation to this problem does not require stabilization. 
As the resolution is decreased by a factor of 4 (Fig.\ \ref{1DparabolicBowlCGlowRes}), 
the inviscid solution begins to develop oscillations although stability is not compromised. 
Unlike the inviscid one, the stabilized solution preserves the surface flatness although it moves slower 
due to a viscous source term in the momentum equation. 
This fact is also reflected in the computation of the normalized $L_2$ error norms plotted in 
 Fig.\ \ref{L2error1DparabolicBowl}. 
The analytic solution against which the error is calculated is defined for an inviscid flow. 
A difference between the viscous and inviscid computation is hence expected. 
The same observation applies to the CG and DG computations alike. The DG solutions to the same problem are 
plotted in Fig.\ \ref{1DparabolicBowl128DG} and \ref{1DparabolicBowlDGlowRes}.

\subsection{2D oscillation in a symmetrical paraboloid}
\label{2dparaboloidTest}
\begin{figure}
\centering
\includegraphics[width=0.49\textwidth]{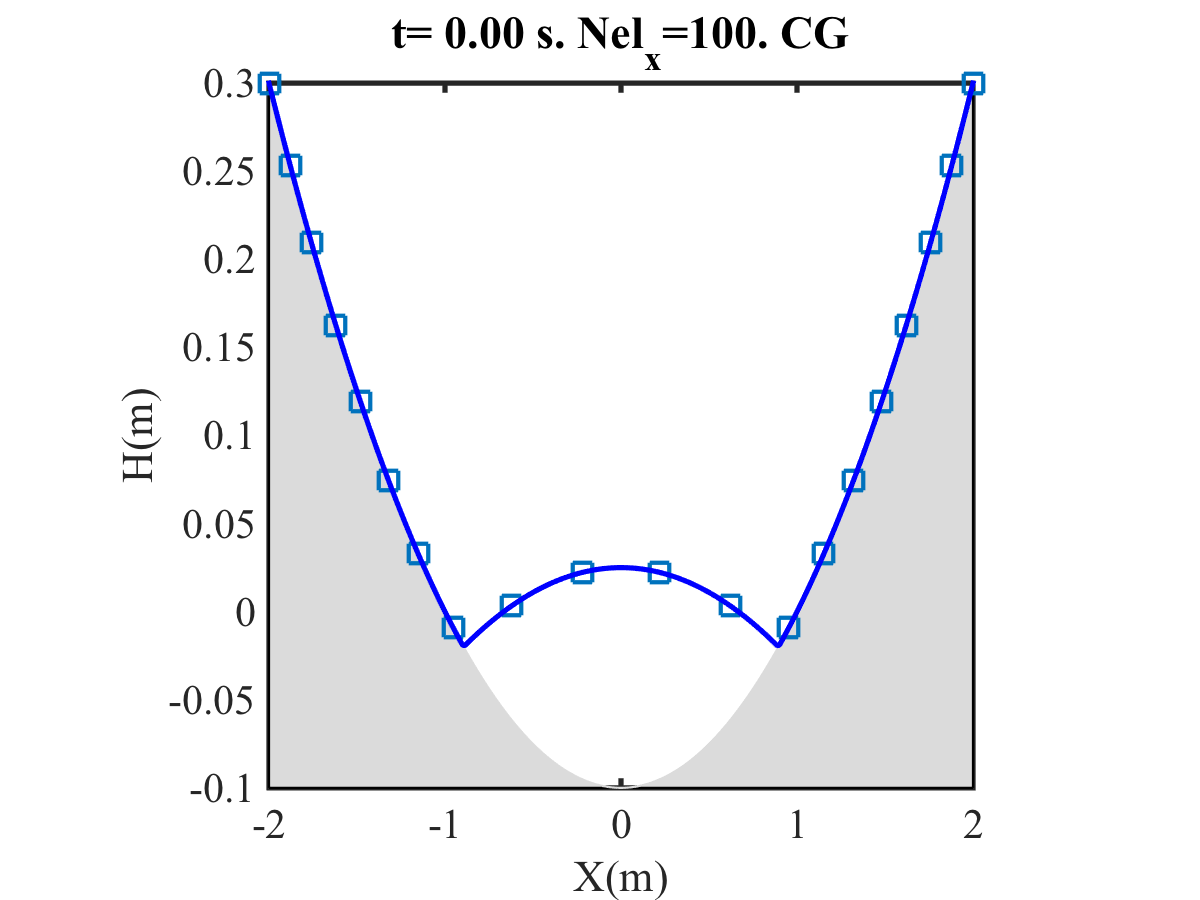}
\includegraphics[width=0.49\textwidth]{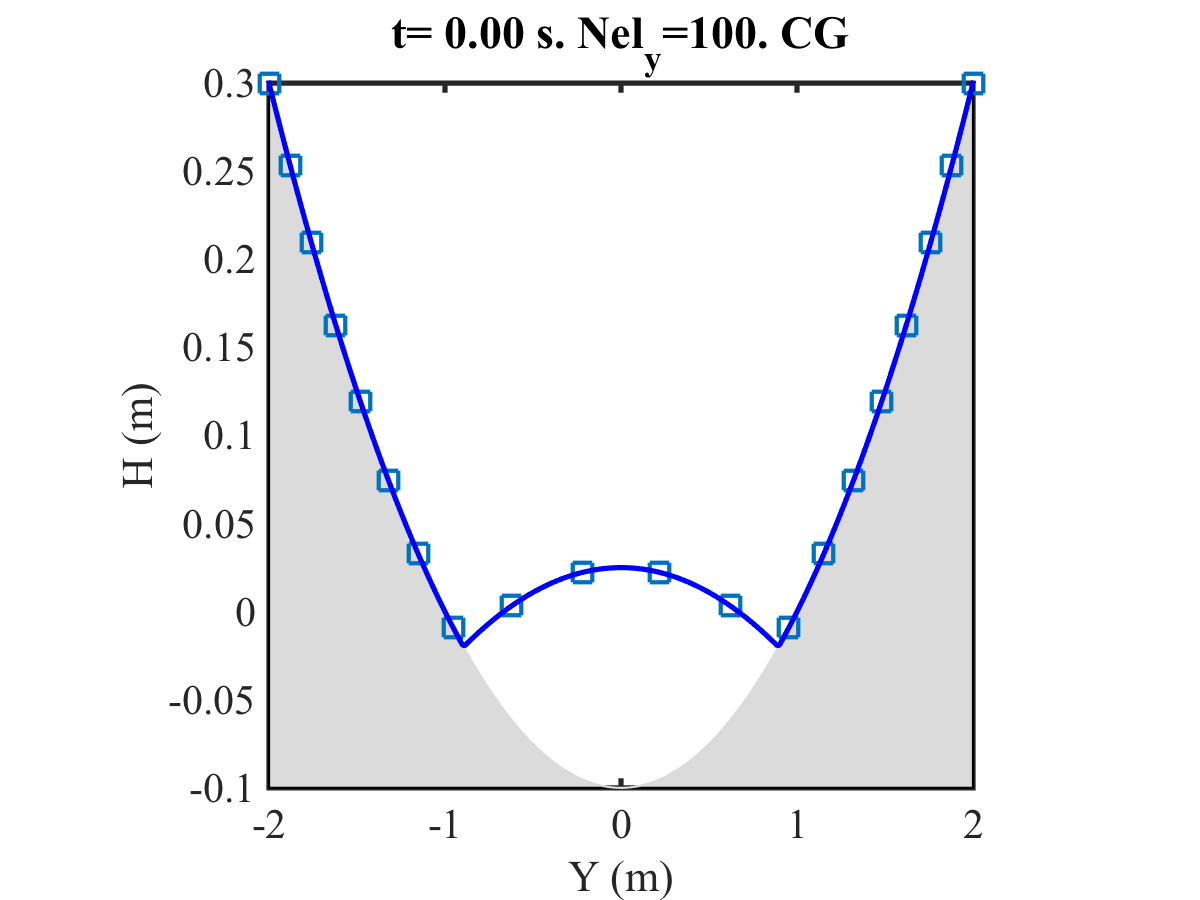}\\

\includegraphics[width=0.49\textwidth]{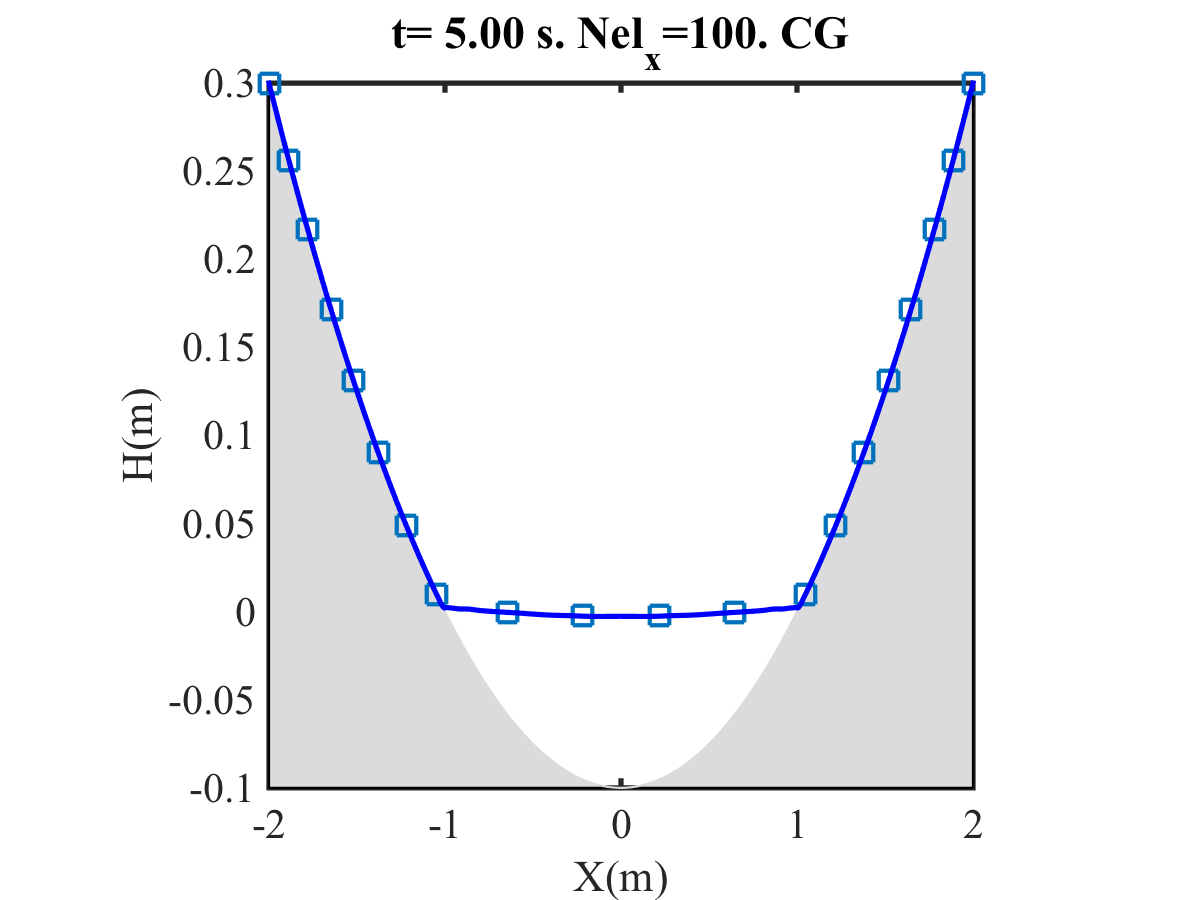}
\includegraphics[width=0.49\textwidth]{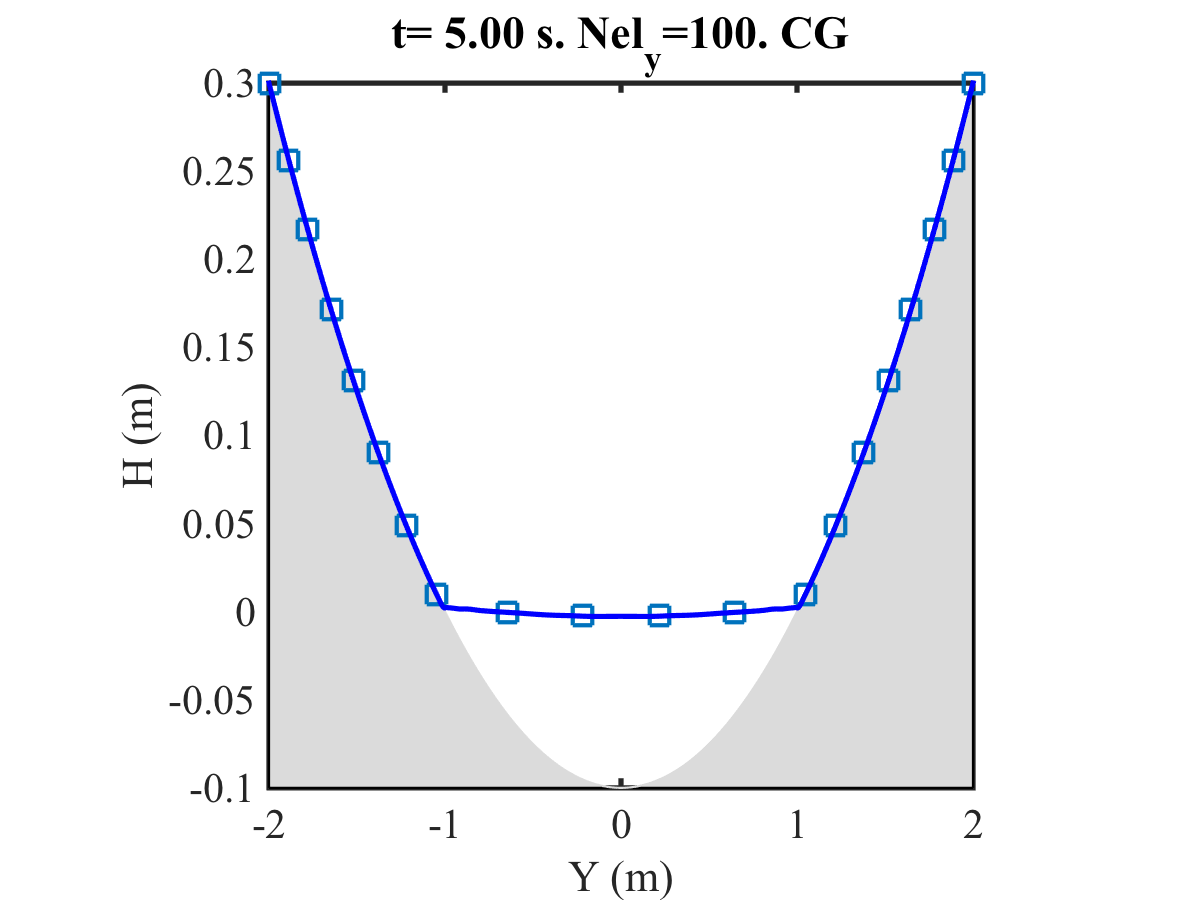}\\

\includegraphics[width=0.49\textwidth]{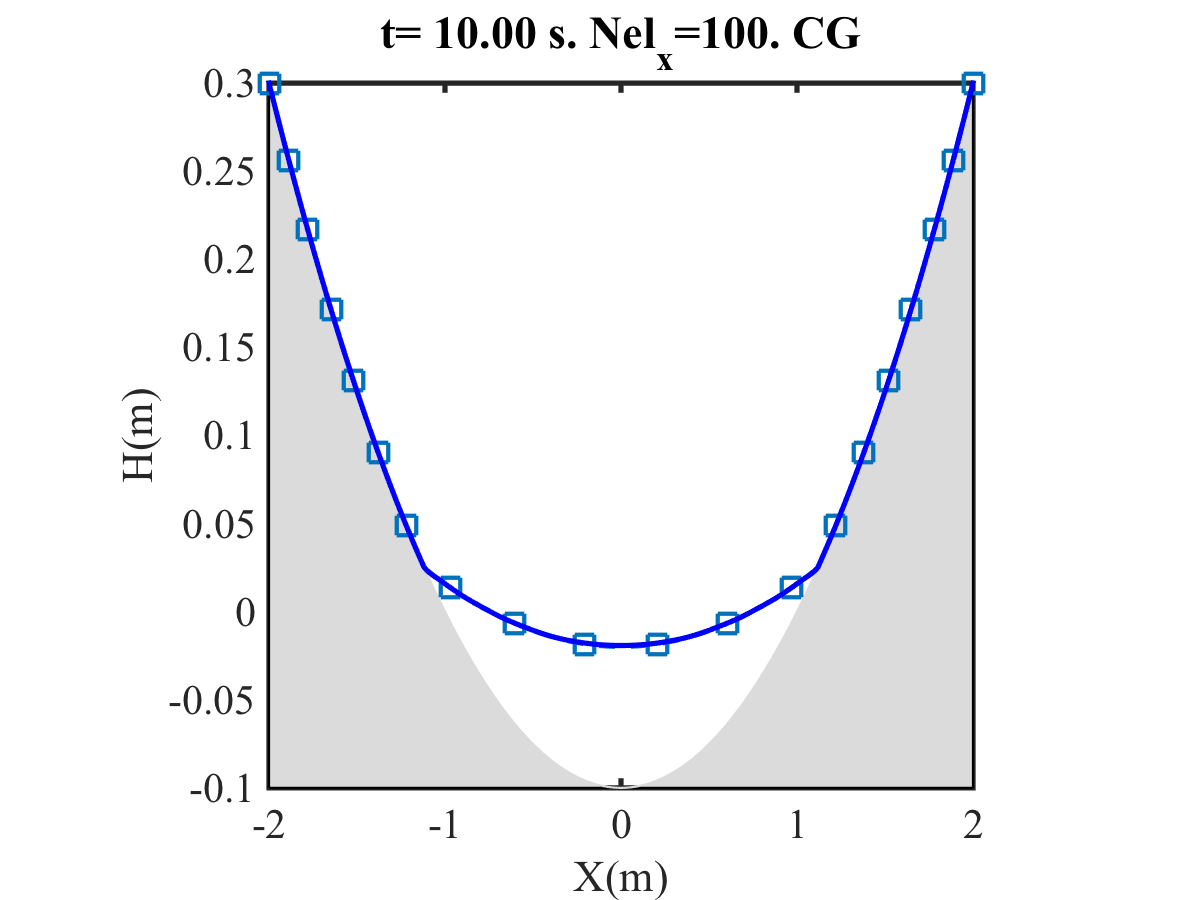}
\includegraphics[width=0.49\textwidth]{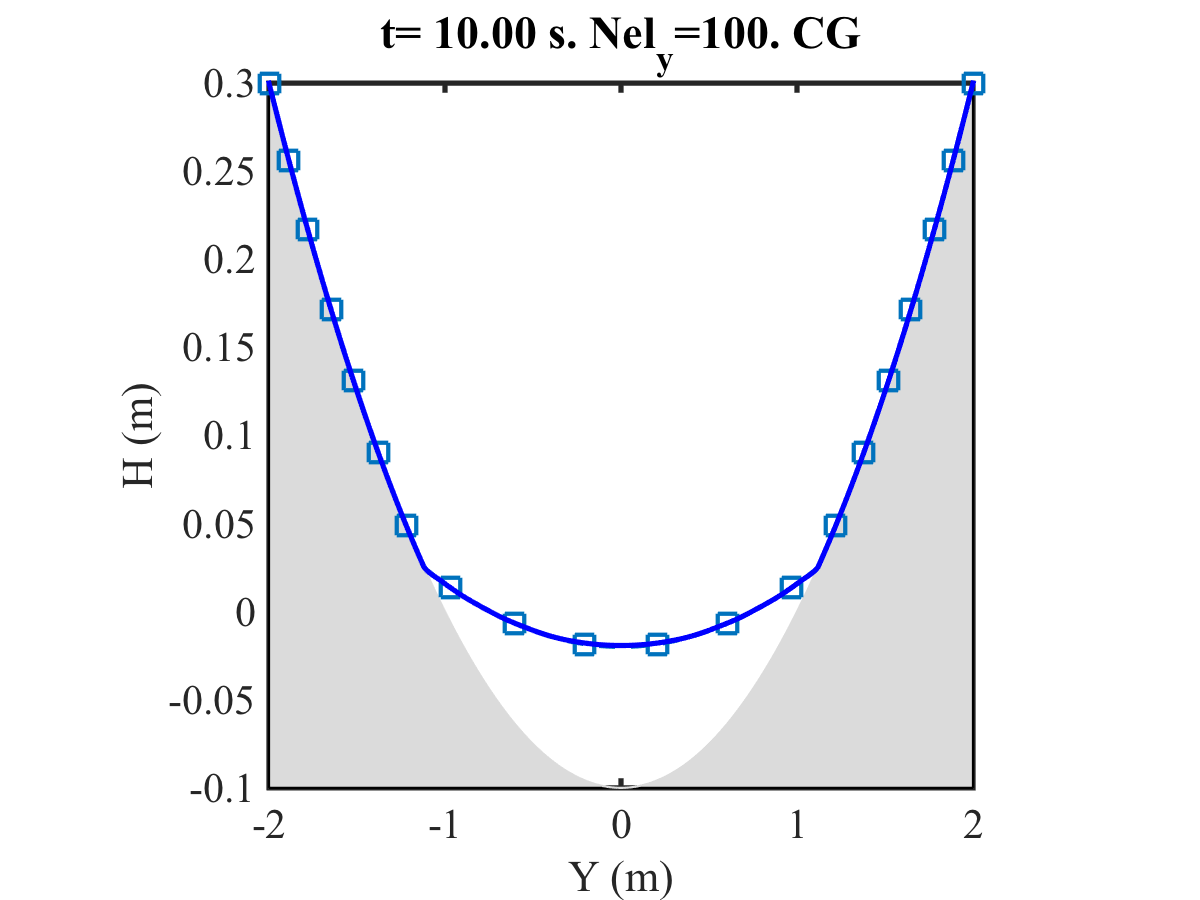}\\
\caption{2D oscillation in a paraboloid. Inviscid {\bf CG} solution. Left: $x-z$ view. Right: $y-z$ view. From top to bottom: $t=[0,\; 5,\; 10]\,s$. 
Solution with $100\times 100$ elements of order 4.}
\label{2DparabolicBowlCG}
\end{figure}

\begin{figure}
\centering
\includegraphics[width=0.49\textwidth]{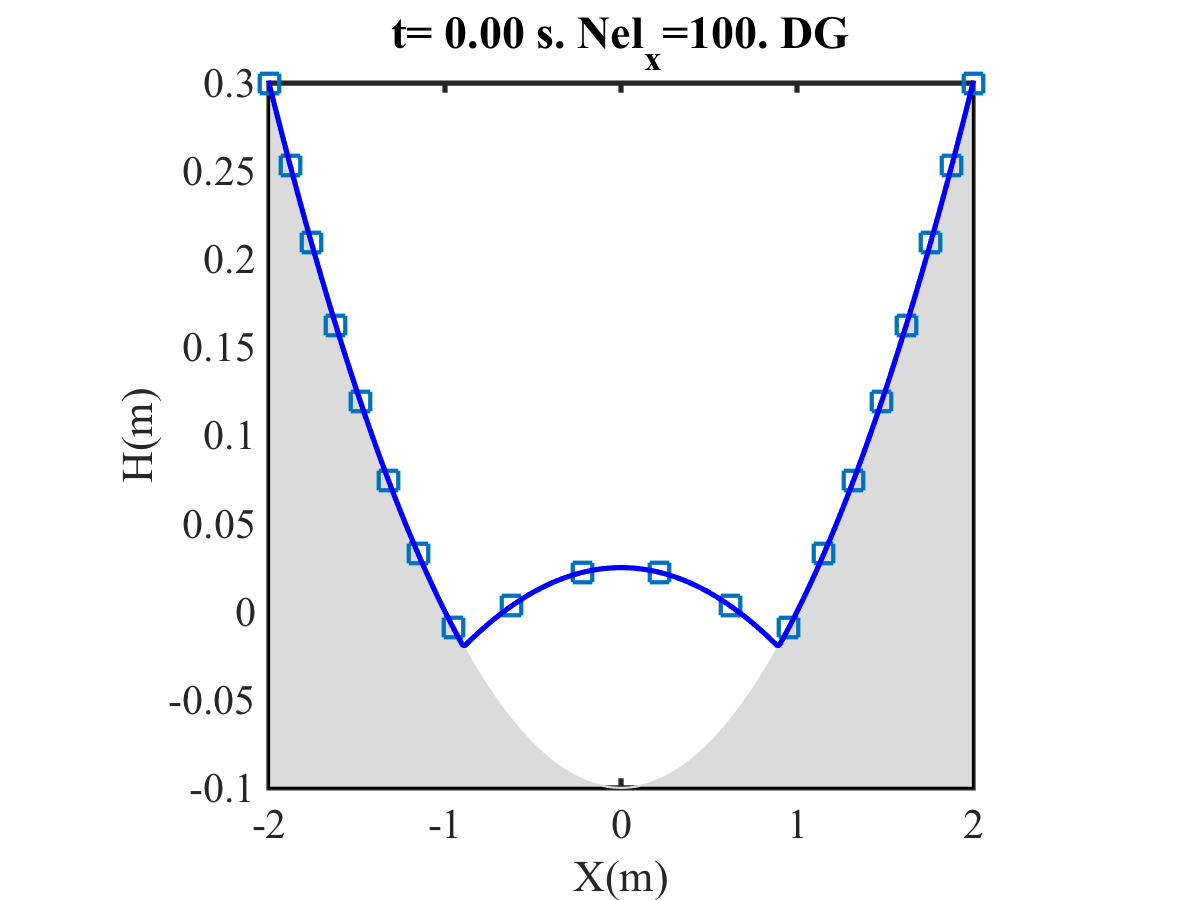}
\includegraphics[width=0.49\textwidth]{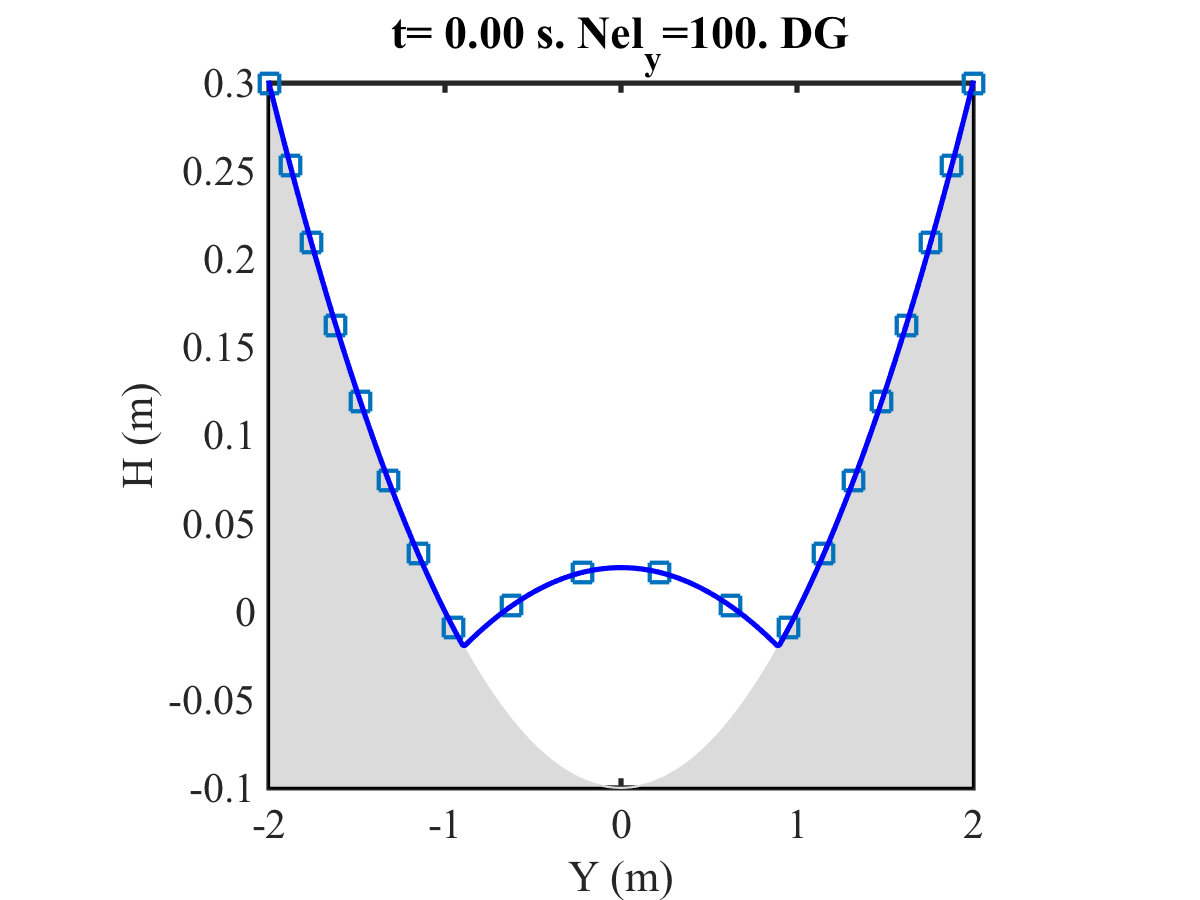}\\

\includegraphics[width=0.49\textwidth]{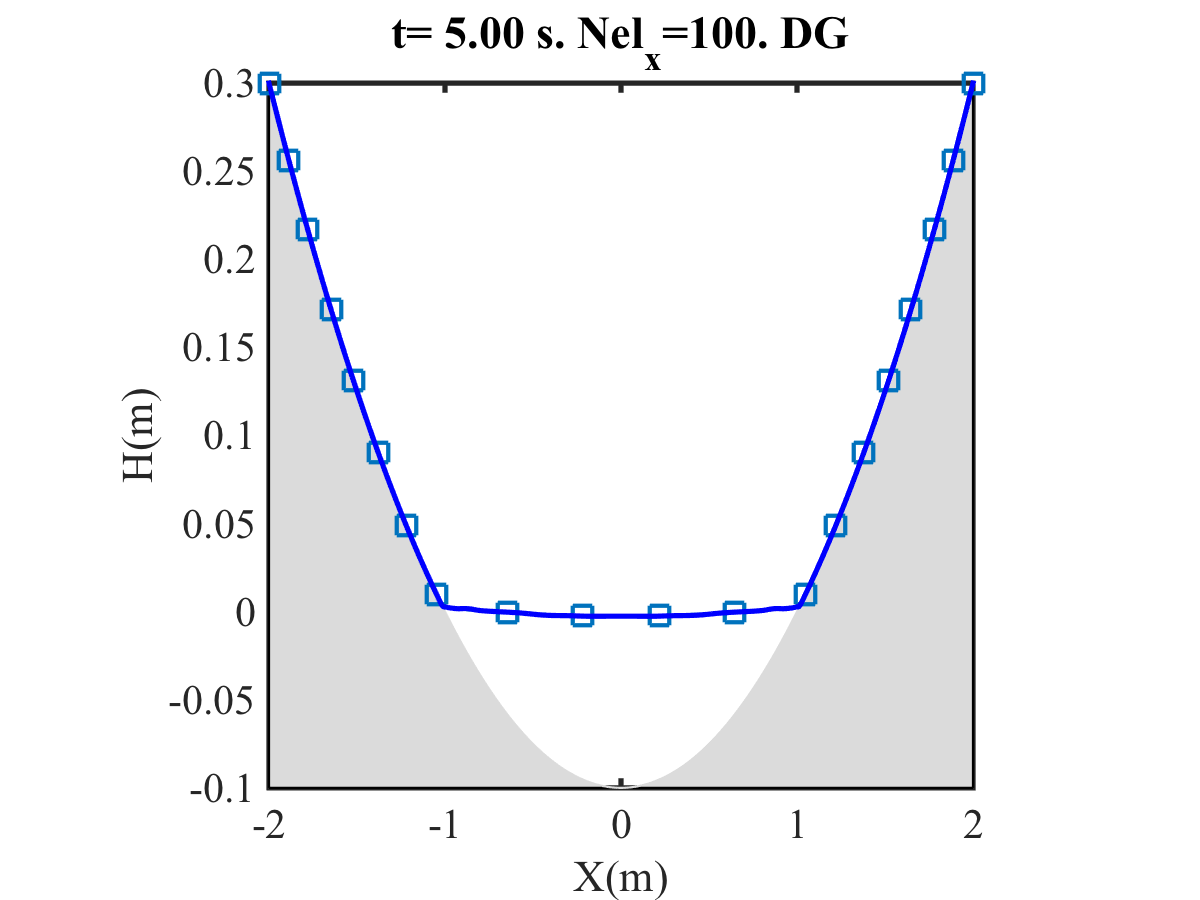}
\includegraphics[width=0.49\textwidth]{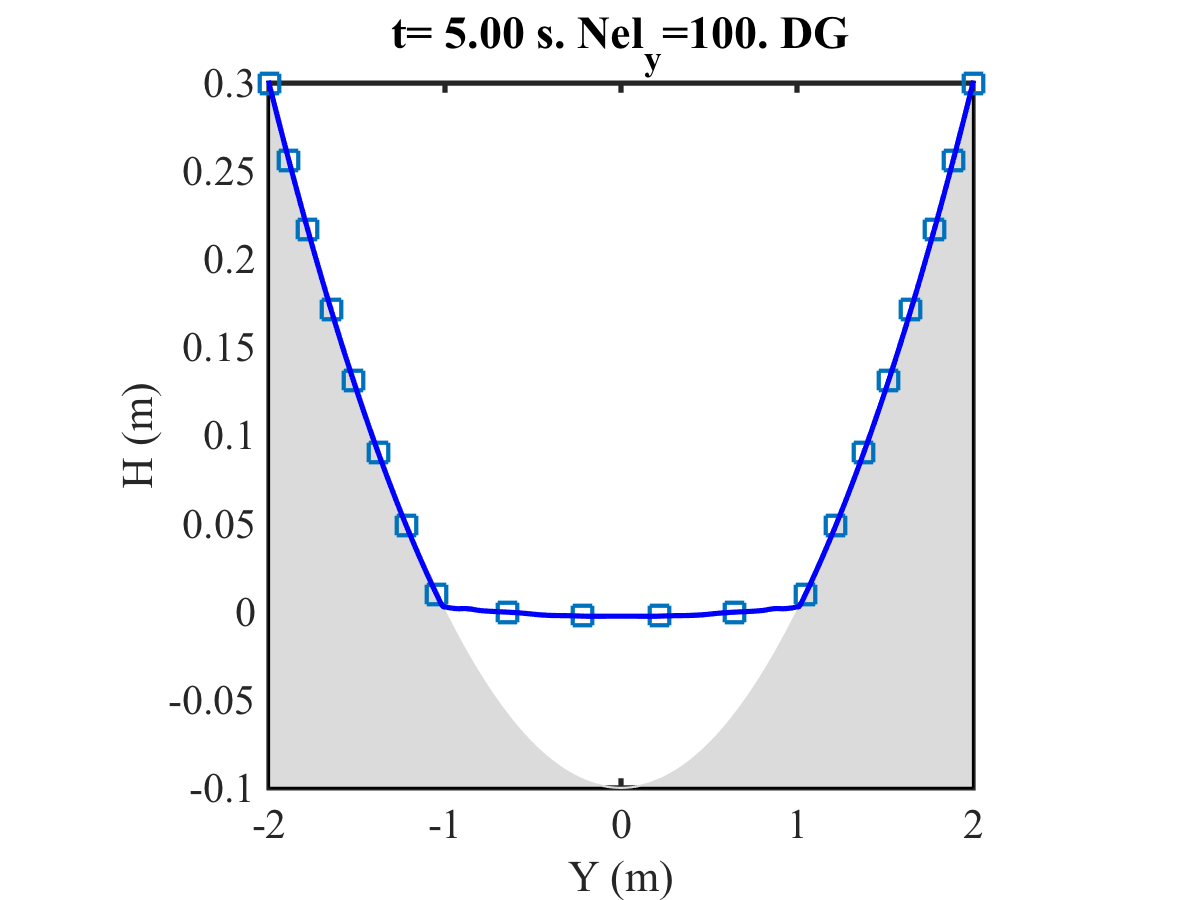}\\

\includegraphics[width=0.49\textwidth]{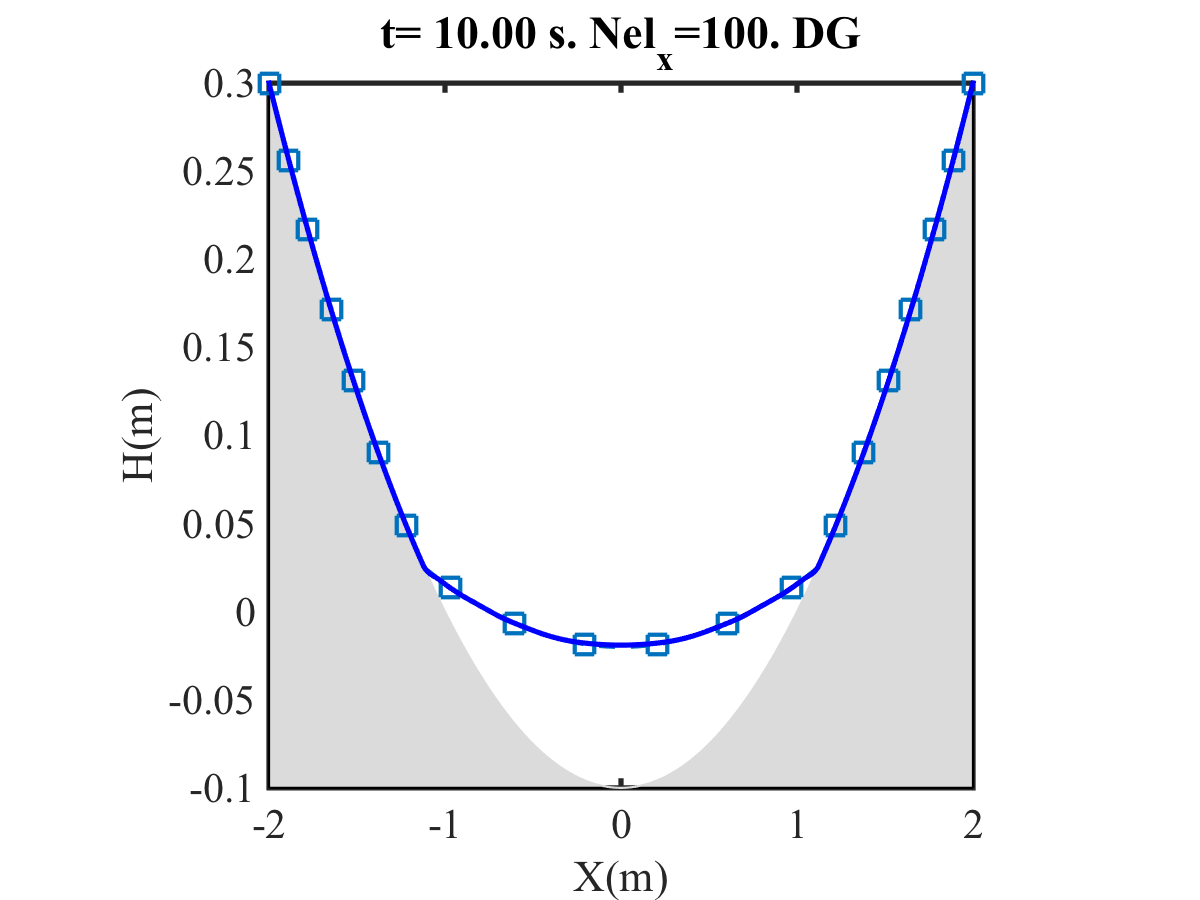}
\includegraphics[width=0.49\textwidth]{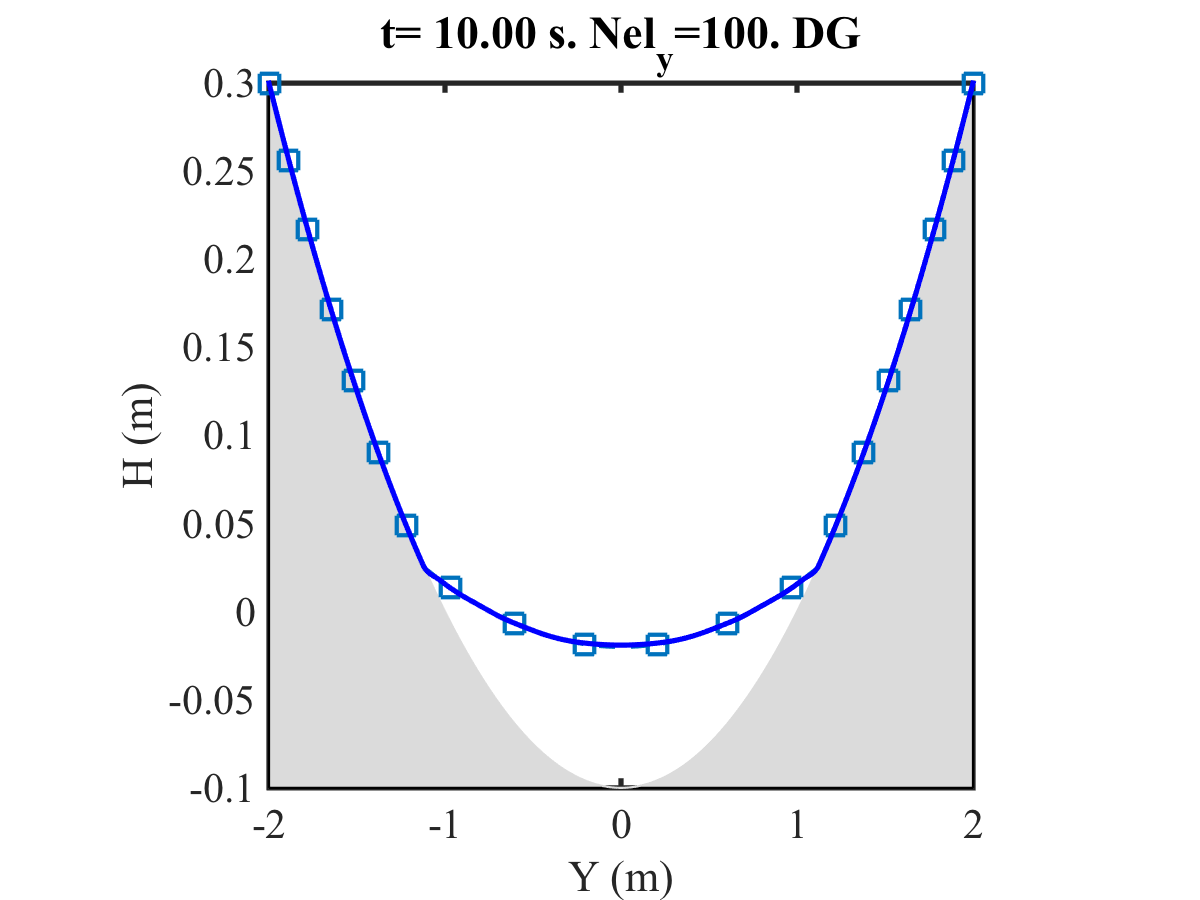}\\
\caption{2D oscillation in a paraboloid. Inviscid {\bf DG} solution. Left: $x-z$ view. Right: $y-z$ view. From top to bottom: $t=[0,\; 5,\; 10]\,s$. 
Solution with $100\times 100$ elements of order 4.}
\label{2DparabolicBowlDG}
\end{figure}

This test was defined by Thacker \cite{thacker1981} who computed the exact solution of the problem. The paraboloid of revolution is defined as :
\[
H_b(x,y) = h_0\left(1 - \frac{\sqrt{x^2 + y^2}}{a^2}\right) - 0.1
\]
in the domain $\Omega=[-2,2]\times [-2,2]\,{\rm m}^2$, with $h_0=0.2$ m.
The initial condition is a reversed paraboloid (See top row of Fig.\ \ref{2DparabolicBowlDG}).
The variation of the bed slope and the radial symmetry are a perfect test to assess the inundation scheme with two dimensional effects.
The inviscid CG and DG solutions are plotted in Figs.\ \ref{2DparabolicBowlCG} and \ref{2DparabolicBowlDG}, respectively, between $t=0$ and $t=10$ s. 
The radial symmetry of the solution is preserved throughout the simulation, and no spurious modes
can be observed in the proximity of the moving shoreline. Although it is not fully visible from the plots, 
we observed a small discrepancy between the exact and 
computed solutions that occurs at the very center of the bowl at the latest time.
More specifically, the numerical solution tends to move slightly slower 
than the exact solution along the bowl centerline. This occurs for CG and DG alike.

\subsection{2D flooding problem in a channel with three mounds}
\label{3islands}
\begin{figure}
\centering
\includegraphics[width=0.49\textwidth]{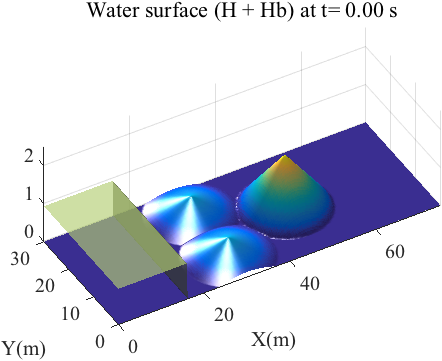}
\includegraphics[width=0.49\textwidth]{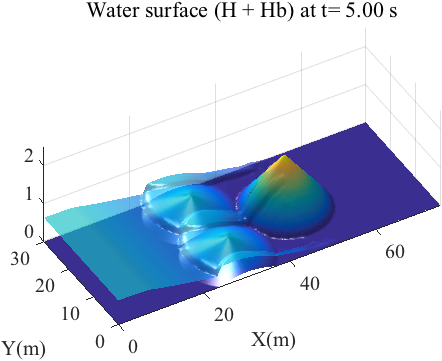}
\includegraphics[width=0.49\textwidth]{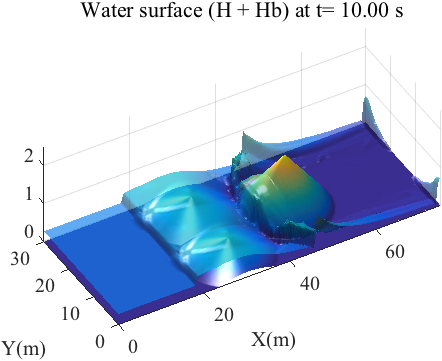}
\includegraphics[width=0.49\textwidth]{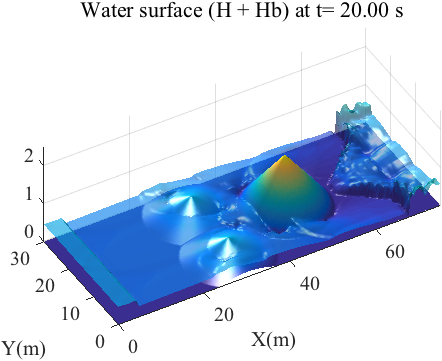}
\includegraphics[width=0.49\textwidth]{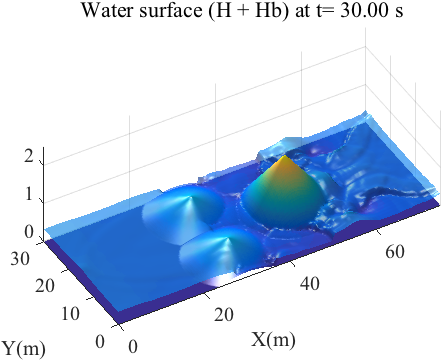}
\includegraphics[width=0.49\textwidth]{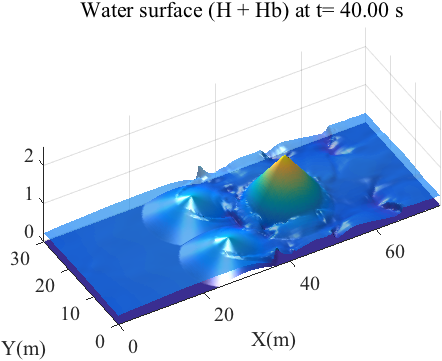}
\caption{Dam break problem in a closed channel with three mounds. Viscous {\bf CG} solution at $t=[0,\,5,\,10,\,20\,30\,40]$ s. 
Solution computed using  $15\times 6$ elements of order 4 for the domain $\Omega=75\times30\,{\rm m^2}$.}
\label{humpConfig1flow}
\end{figure}

This test \cite{kawaharaUmetsu1986,brufauEtAl2002} was used by Xing and Zhang \cite{xingZhang2013} and Gallardo et al.\ \cite{gallardoEtAl2007} to assess, on unstructured grids, the same positivity-preserving method 
that we are applying here \cite{xingZhangShu2010}. We reproduce their results here to verify its correct implementation in our code.
The three mounds in the channel $\Omega=75 \times 30\, {\rm m^2}$, are defined by the function $H_b({\bf x})=\max(0.0, m_1,m_2,m_3)$ where
\begin{subequations}
\[
m_1=1.0-0.10\sqrt{(x-30.0)^2 + (y-22.5)^2},
\]
\[
m_2=1.0-0.10\sqrt{(x-30.0)^2 + (y-7.50)^2},
\]
\[
m_3=2.8-0.28\sqrt{(x-47.5)^2 + (y-15.0)^2},
\]
\end{subequations}
The flooding is triggered by a dam break at $t=0$ s. 
The evolution of the flooding from $t=0$ to $t=40$ s is shown in Fig. \ref{humpConfig1flow}. 
We computed the solution using a relatively coarse grid made of $15\times 6$ elements of order 4. 
The results are in agreement with \cite{gallardoEtAl2007} and \cite{xingZhang2013}: 
the flow separation is well resolved and occurs after 10 seconds from the dam breaking. 
After 15 seconds, the first water front has reached the back wall and is fully reflecting back by 20 seconds. 
By 40 seconds, the flow has almost reached full rest.

\subsection{2D solitary wave runup and run-down on a circular island}
\label{OneislandTest}
A solitary wave runup on a circular island was studied in \cite{briggsEtAl1995} at the Waterways Experiment Station 
of the US Army Corps of Engineers. In this example, the initial wave is modeled via the following analytic definition by Synolakis \cite{synolakis1987}:
\begin{equation}
\label{synowave}
	\eta({\bf x}, 0) = \frac{A}{h_0} {\rm sech}^2 \left( \gamma(x - x_c) \right),
\end{equation}
where $A=0.064$ m is the wave amplitude, $x_c=2.5$ m, $h_0=0.32$ m is the initial still water level, and 
\[
\gamma = \sqrt{\frac{3A}{4h_0}}.
\]

The island is a cone given as
\[
H_b = 0.93 \left( 1 - \frac{r}{r_c} \right), \quad {\rm if}\; r\leq r_c,
\]
where $r=\sqrt{ (x - x_c)^2 + (y-y_c)^2}$, $r_c=3.6$ m, and is centered at $(x_c, y_c)=(12.5, 15)$ m. The cone is 
installed on a flat bathymetry. The fluid is confined within four solid walls.

\begin{figure}
\centering                                                                                                                                                 
\includegraphics[width=0.49\textwidth]{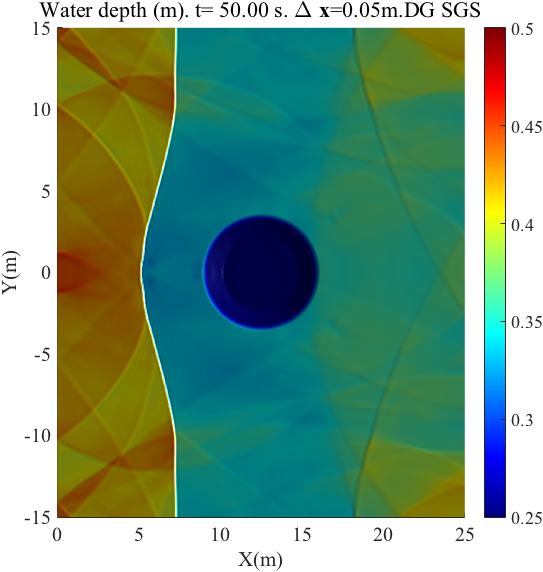}
\includegraphics[width=0.49\textwidth]{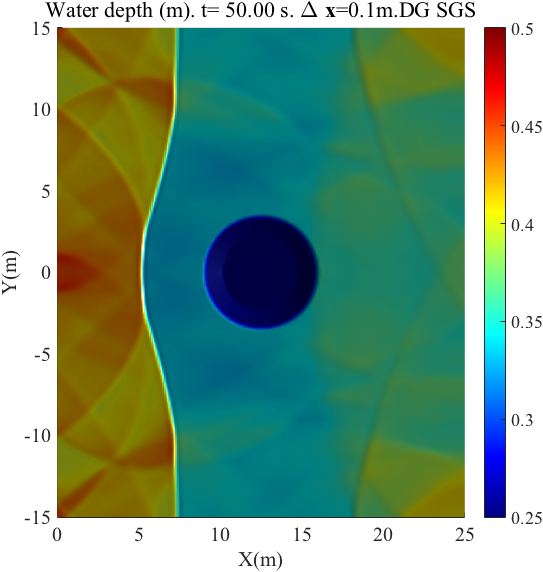}
\includegraphics[width=0.49\textwidth]{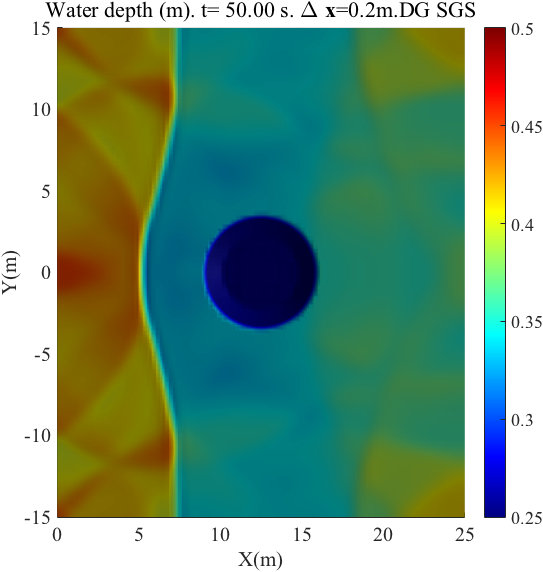}
\includegraphics[width=0.49\textwidth]{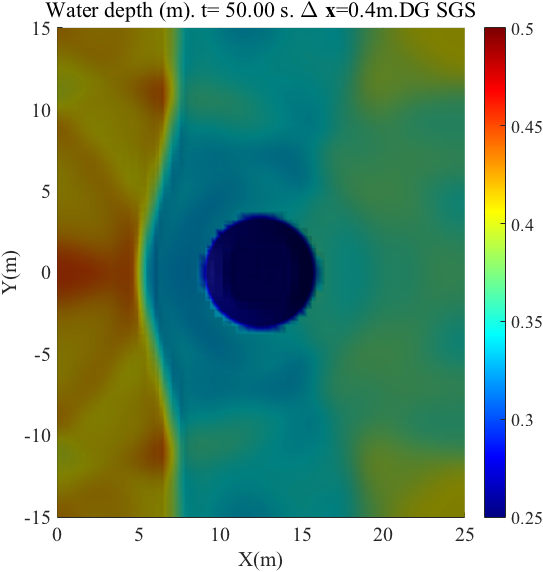}
\caption{{\bf DG} solution of water depth for the single-hill configuration. Results obtained with the four grid resolutions $\Delta{\bf x} = [0.05, 0.10, 0.20, 0.40]$ m (indicated in the figures) using $4^{th}$-order elements. The color bar is cut at 0.25 
to preserve the visibility of the smallest features. The dark blue coloring within the region of the cone corresponds to the water depth equal to the threshold water layer of 1e-3 m.}
\label{briggs4Resolutions}
\end{figure}

To understand how the proposed diffusion and numerical approximation depend on the grid, 
we ran the simulation at the four resolutions $\Delta {\bf x} \approx [0.05\,,0.10\,,0.20\,,0.40]$ m. 
Fig.\ \ref{briggs4Resolutions} shows that the stabilized DG solution is converging to the same solution and the main features 
of the interacting waves are reproduced almost equally across the four grids. Certainly, the 40 cm grid spacing is the most 
dissipative, although it is encouraging to see how the important features resolved at 5 cm are still well represented on the coarsest grid. 
The same observation applies to the CG solution (plot not shown).

\begin{figure}
\centering                                                                                                                                                 
\includegraphics[width=0.49\textwidth]{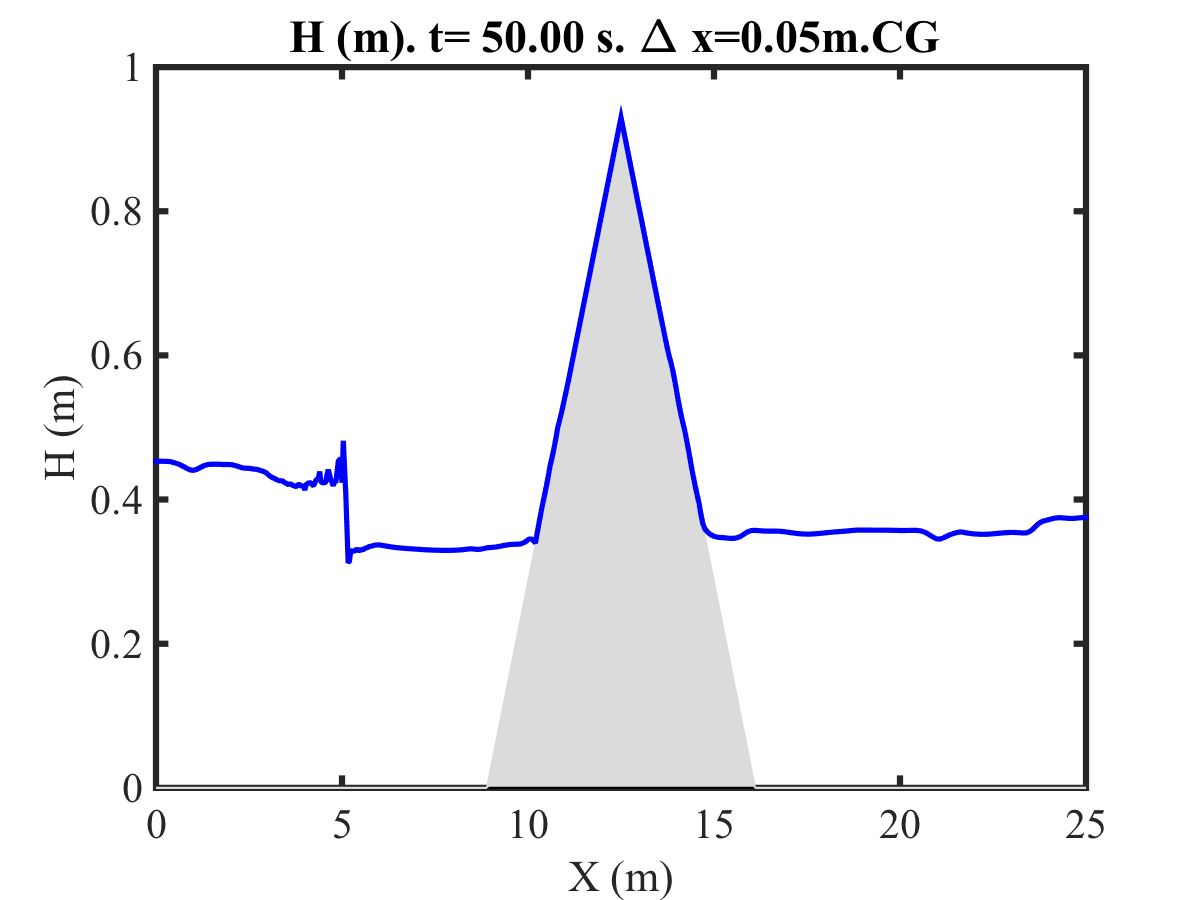}
\includegraphics[width=0.49\textwidth]{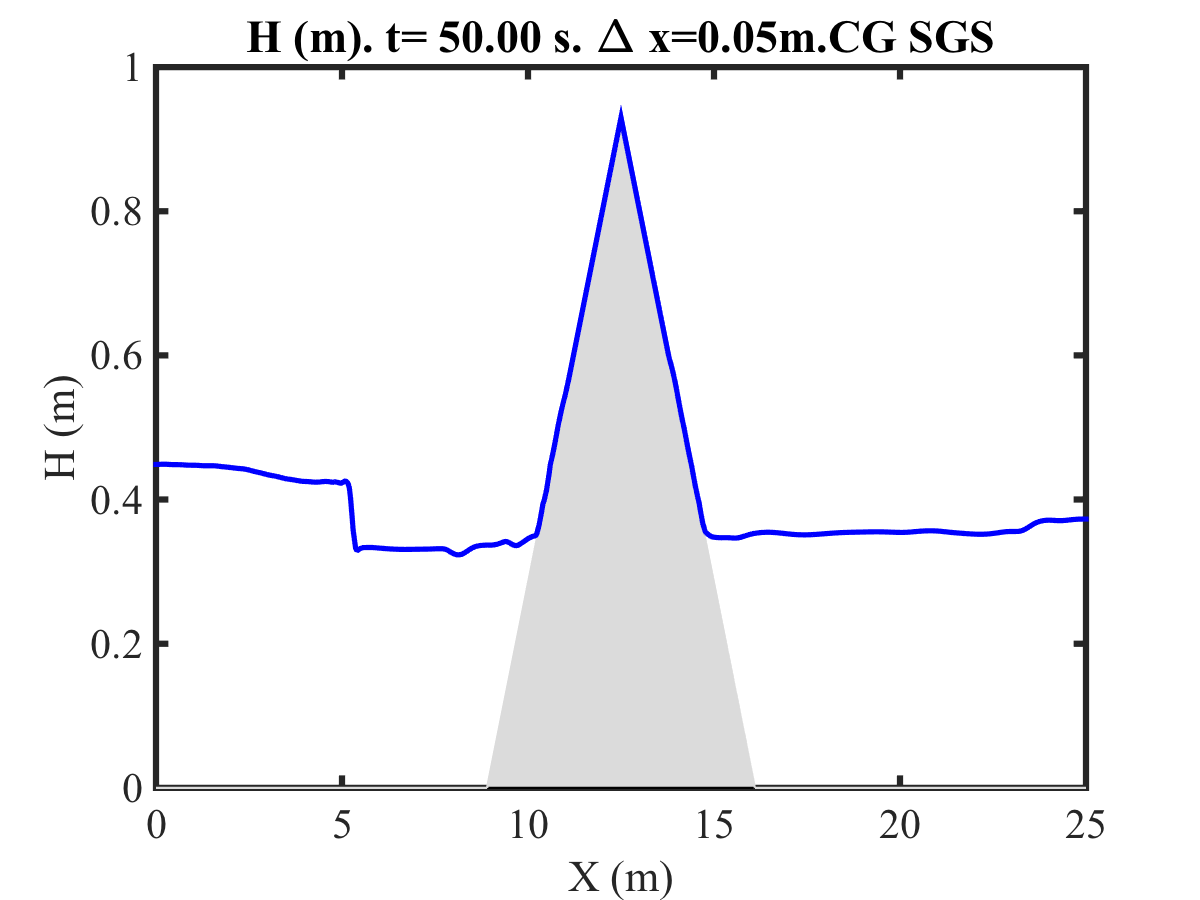}
\caption{$x-z$ view of the 2D {\bf CG} solution of water depth $(H)$ for the single-hill configuration. 
Inviscid (left) against viscous solution using SGS (right).
Solutions obtained using $4^{th}$-order elements. 
The dark blue coloring within the region of the cone corresponds to the water depth equal to the threshold water layer of 1e-3 m.}
\label{briggs4Resolutions1DsliceXZ_CG}
\end{figure}

\begin{figure}
\centering                                                                                                                                                 
\includegraphics[width=0.49\textwidth]{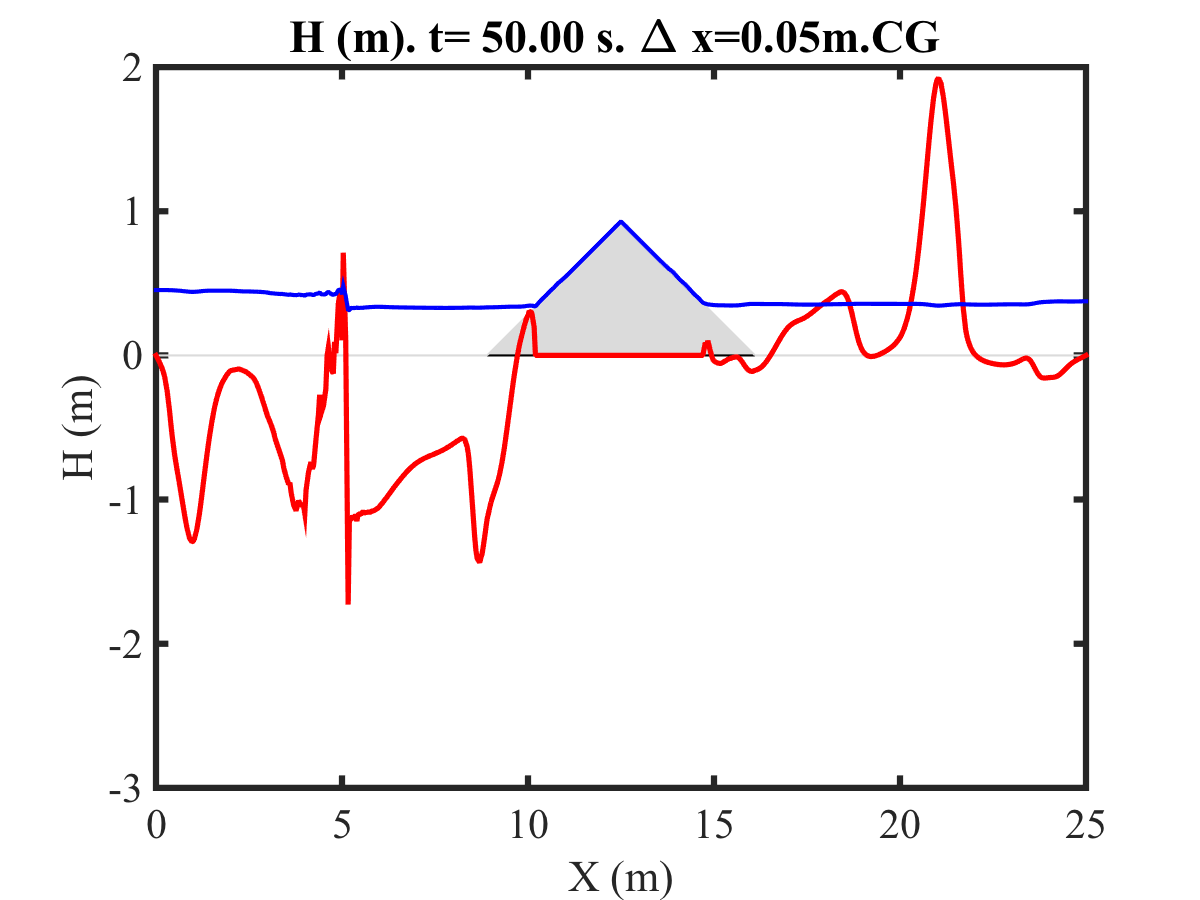}
\includegraphics[width=0.49\textwidth]{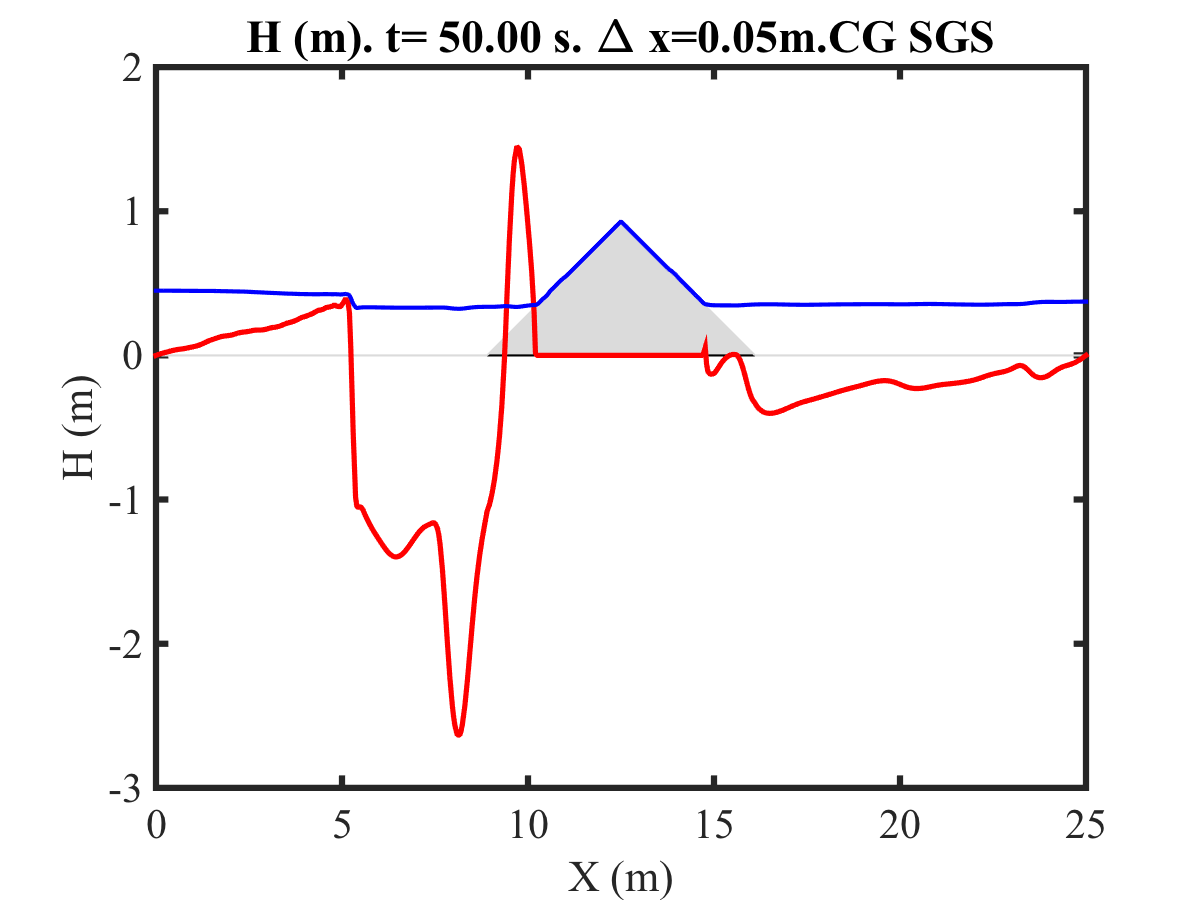}
\caption{Like Fig. \ref{briggs4Resolutions1DsliceXZ_CG}, but for $u$-velocity component. Inviscid (left) against viscous solution using SGS (right).}
\label{briggs4Resolutions1DsliceXZ_CG_UVELO}
\end{figure}

\begin{figure}
\centering                                                                                                                                                 
\includegraphics[width=0.49\textwidth]{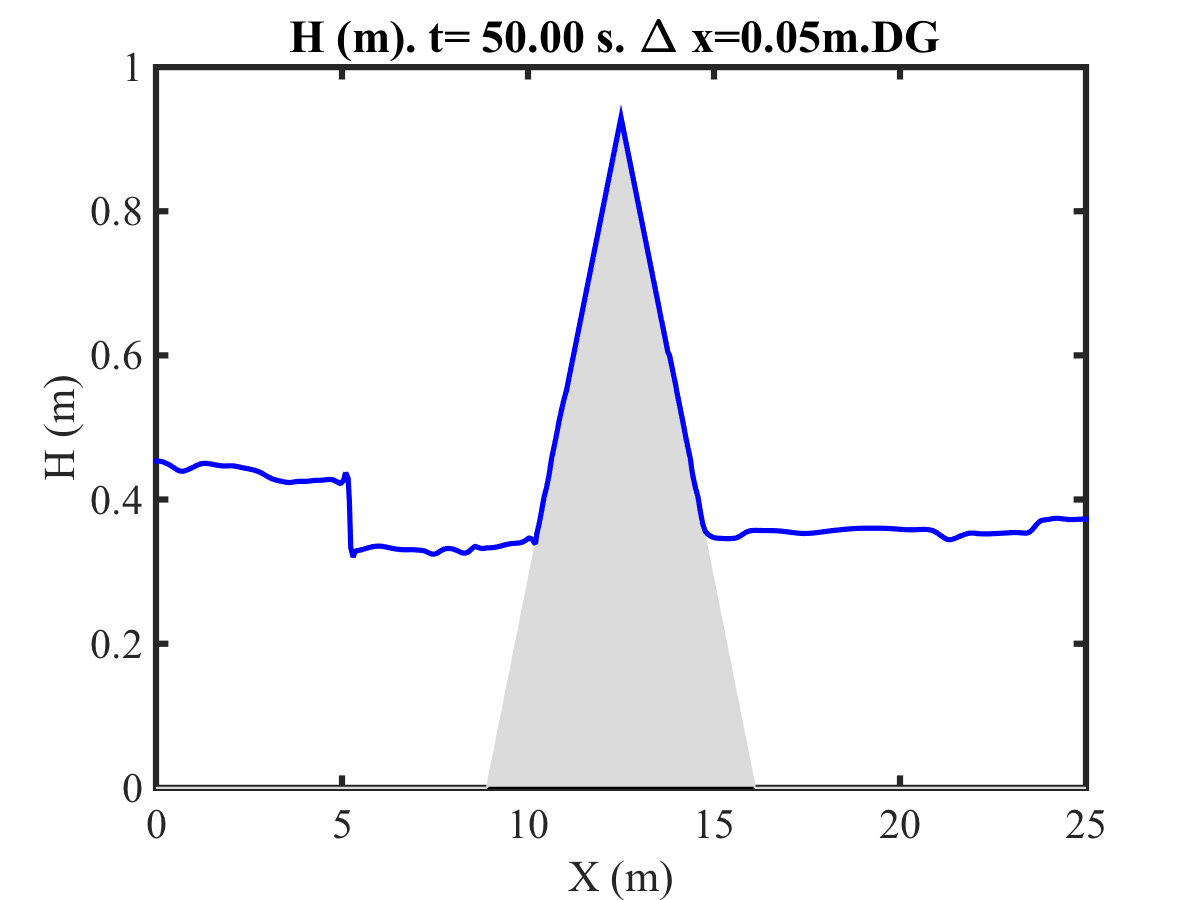}
\includegraphics[width=0.49\textwidth]{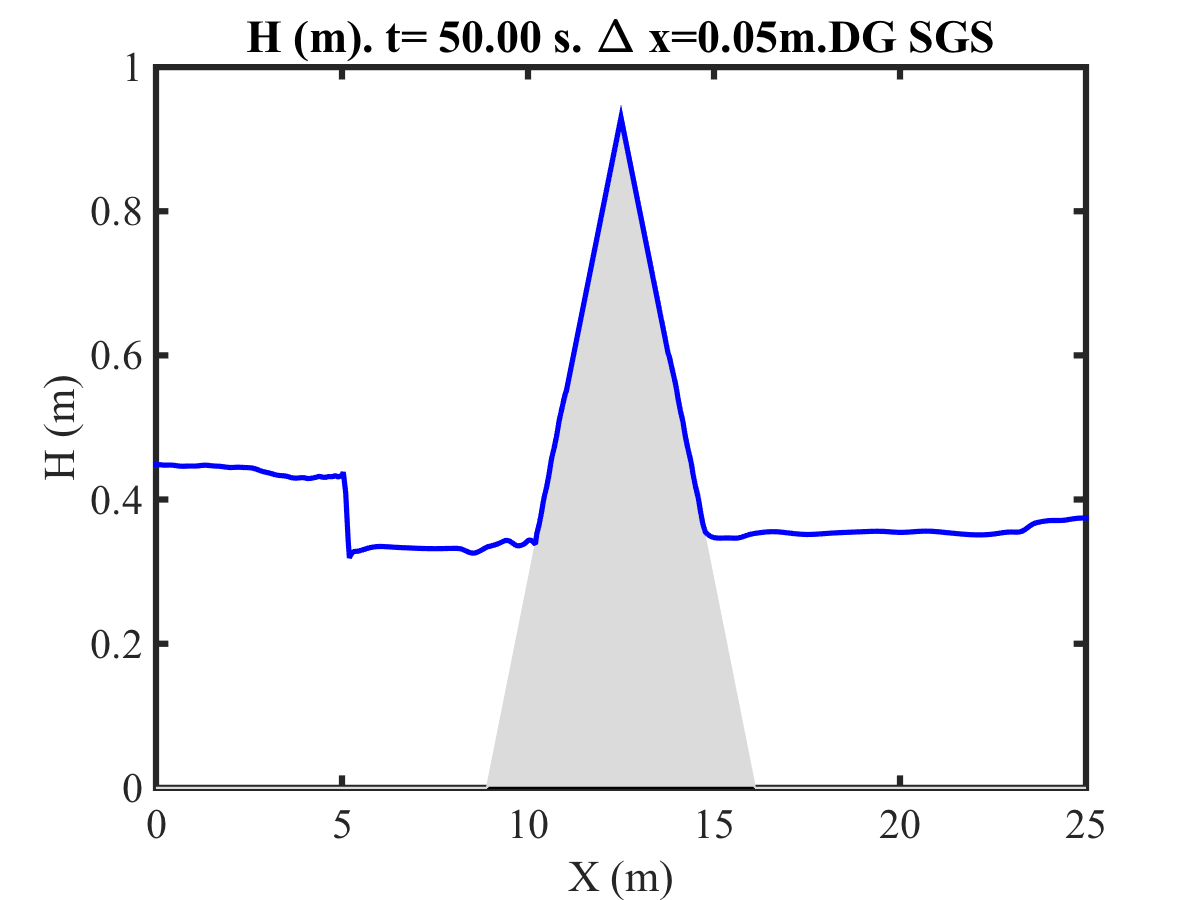}
\caption{$x-z$ view of the 2D {\bf DG} solution of water depth $(H)$ for the single-hill configuration. 
Inviscid (left) against viscous solution using SGS (right).
Solution obtained using $4^{th}$-order elements. 
The dark blue coloring within the region of the cone corresponds to the water depth equal to the threshold water layer of 1e-3 m. This plot shows the power of DG. Without SGS it still almost captures the bore sharply.}
\label{briggs4Resolutions1DsliceXZ_DG}
\end{figure}

\begin{figure}
\centering                                                                                                                                                 
\includegraphics[width=0.49\textwidth]{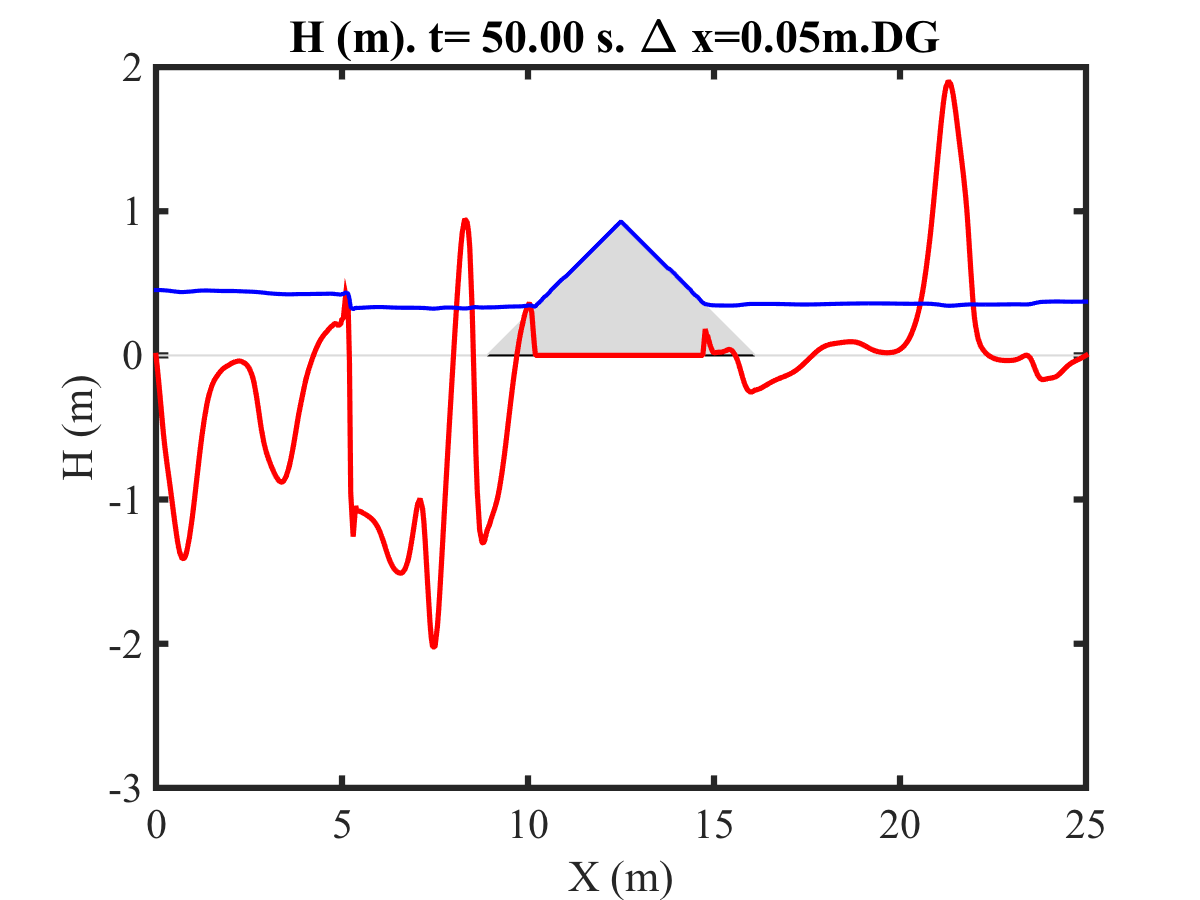}
\includegraphics[width=0.49\textwidth]{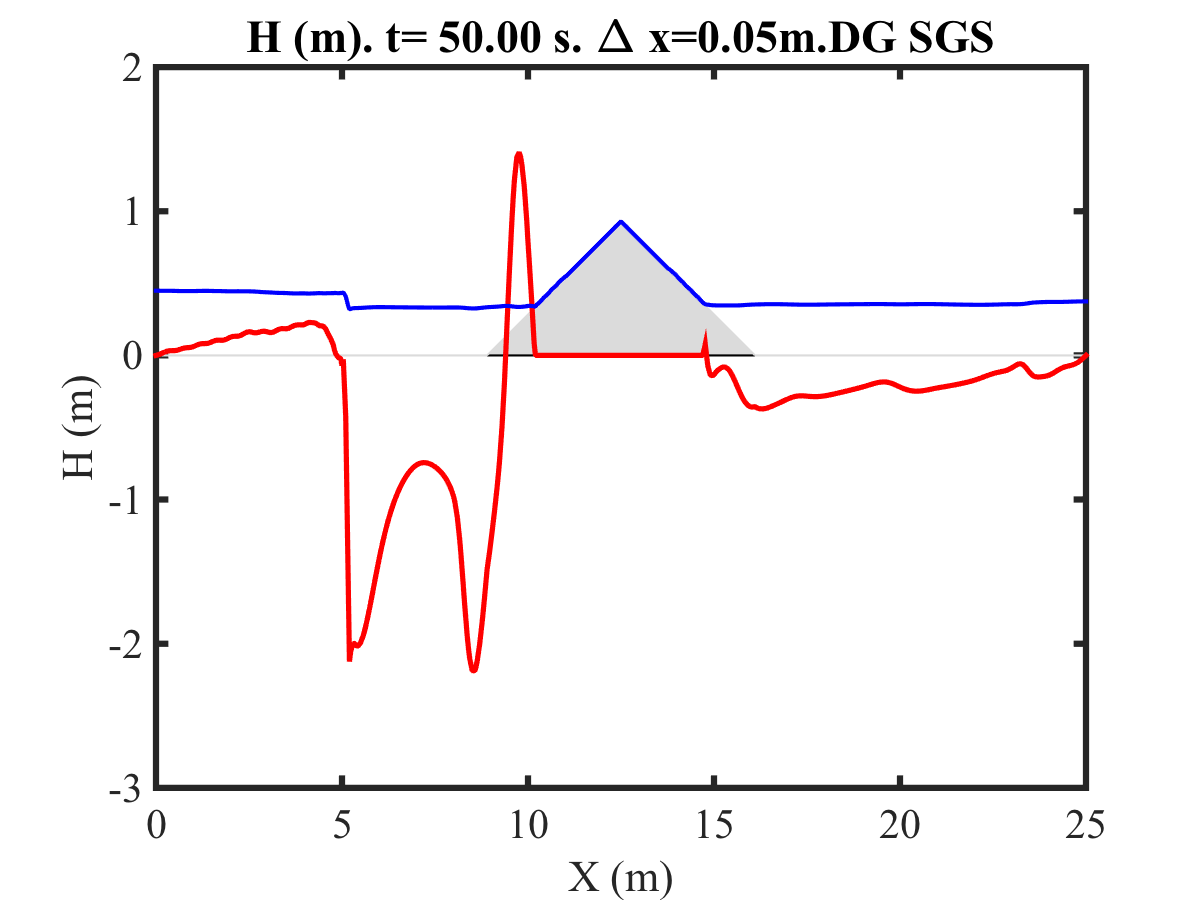}
\caption{Like Fig. \ref{briggs4Resolutions1DsliceXZ_DG}, but for $u$-velocity component. 
Inviscid (left) against viscous solution using SGS (right). 
A comparison of this figure with Fig.\ \ref{briggs4Resolutions1DsliceXZ_CG_UVELO}  
shows the same dispersion properties of the CG and DG solutions.}
\label{briggs4Resolutions1DsliceXZ_DG_UVELO}
\end{figure}

In Figs.\ \ref{briggs4Resolutions1DsliceXZ_CG}-\ref{briggs4Resolutions1DsliceXZ_DG_UVELO} we plot the projection of the 2D solution  
on the plane $y=0$. 
The spurious modes that characterize the water surface in the proximity of the sharpest wave front are fully removed by {\it Dyn-SGS}
(Fig.\ \ref{briggs4Resolutions1DsliceXZ_CG}) without weakening the front sharpness. We observe the same in the velocity field plotted in Fig.\ \ref{briggs4Resolutions1DsliceXZ_CG_UVELO}.
This is in full agreement with the application of {\it Dyn-SGS} to non-linear wave problems with strong discontinuities, as previously shown in \cite[Figs.\ 16, 17]{marrasEtAl2015b} for the solution of the Burgers' equation. 
We briefly mentioned above how DG already has built-in dissipation. 
This is clearly visible from the plot of Fig.\ \ref{briggs4Resolutions1DsliceXZ_DG}; the unstabilized DG solution shows no 
oscillations except for, at most, some minimal under- and over-shooting. This implies that its residual is so small that the effect of 
SGS eddy viscosity reduces to a minimum value. This explains why the inviscid and viscous DG solutions look similar.

It is important to notice the same space distribution of the CG and DG velocity waves shown in Figs.\ \ref{briggs4Resolutions1DsliceXZ_CG_UVELO} and \ref{briggs4Resolutions1DsliceXZ_DG_UVELO}.
This is a visual confirmation that the two solutions present the same dispersion properties, which is 
to be expected from two analogous numerical approximations to the same problem. This is indicative 
of a correct implementation of the unified CG/DG framework.

\begin{figure}
\centering                                                                                                                                                 
\includegraphics[width=0.9\textwidth]{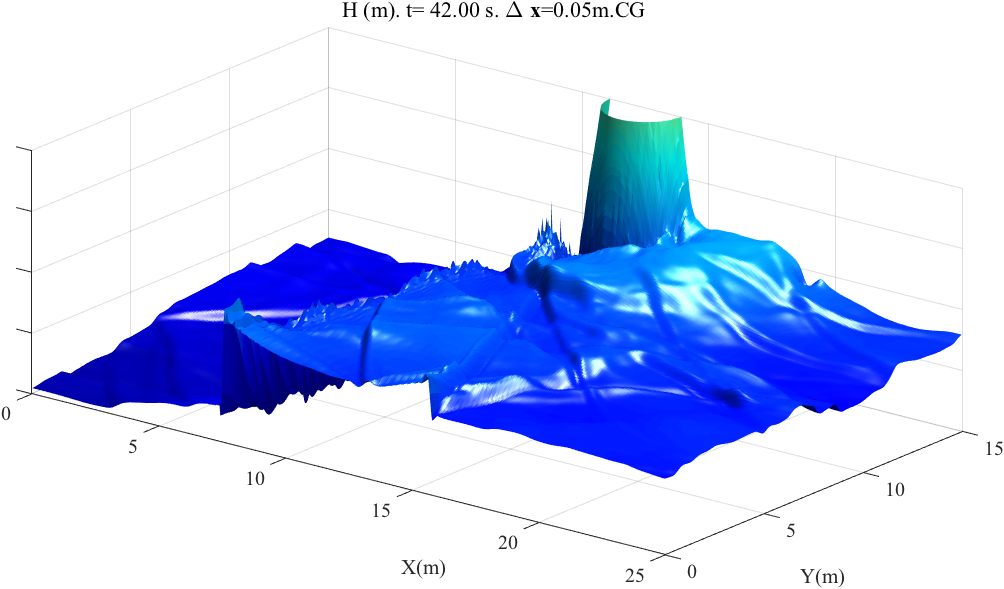}
\includegraphics[width=0.9\textwidth]{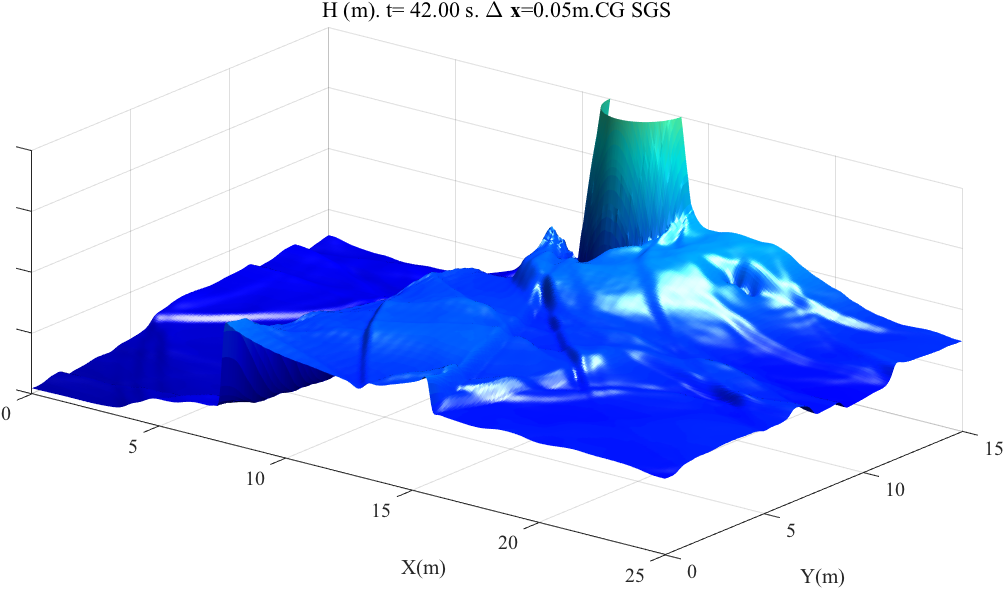}
\caption{{\bf CG} solution for the single-hill configuration. 
Instantaneous perspective view of the unstabilized (top) and stabilized (bottom) water surface for $\Delta{\bf x} = 5$ cm using $4^{th}$-order elements. The high frequency instabilities are removed by Dyn-SGS without compromising the overall sharpness of the interacting waves. Both solutions are characterized by the same wave features at all scales.}
\label{briggsCGstableVSunstableH3dview}
\end{figure}

\begin{figure}
\centering                                                                                                                                                 
\includegraphics[width=0.9\textwidth]{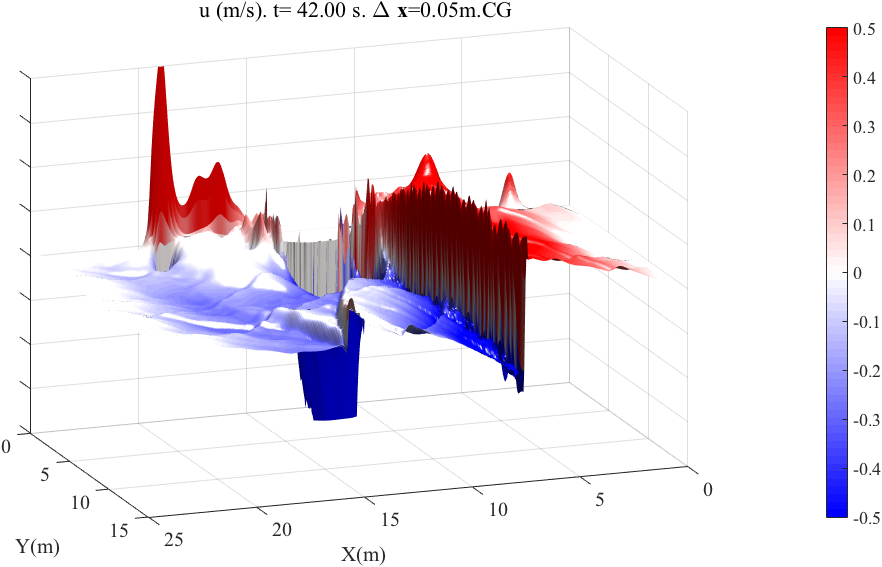}
\includegraphics[width=0.9\textwidth]{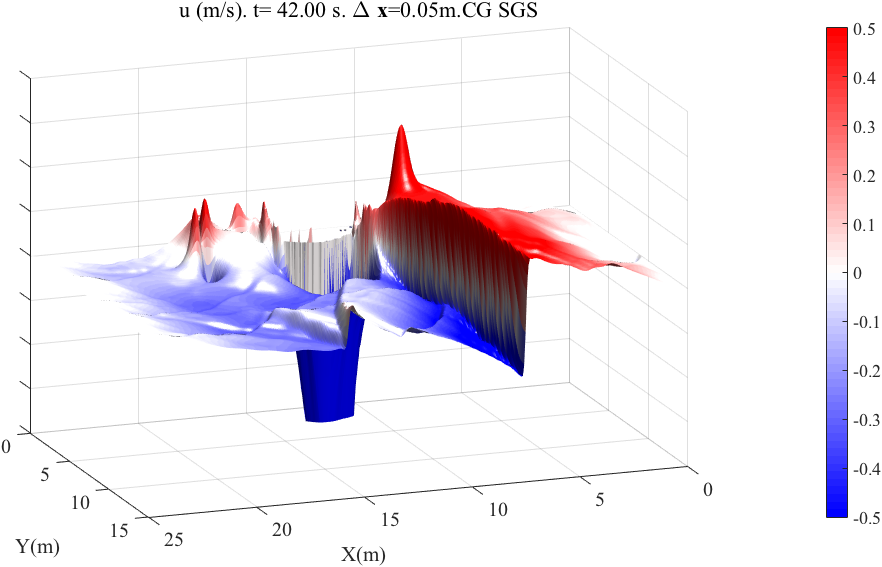}
\caption{Like Fig.\ \ref{briggsCGstableVSunstableH3dview}, but surface of horizontal velocity, $u$. 
For best view of the velocity surface, the view angle is different from the one of 
Fig.\ \ref{briggsCGstableVSunstableH3dview}.}
\label{briggsCGstableVSunstableUVELO3dview}
\end{figure}

The solution of the top plot in Fig.\ \ref{briggsCGstableVSunstableH3dview}, 
was computed with stabilization applied to the continuity ($\delta=1$ in Eq.\ (\ref{SWE})) and momentum equations.
When compared against the un-stabilized solution (bottom plot in the same figure), 
we notice that the features of the fronts of the interacting waves are fully preserved. 
Furthermore, the fronts are not excessively smeared out as the high frequency modes are removed.
The effect of momentum stabilization on the velocity field is shown in 
Fig.\ \ref{briggsCGstableVSunstableUVELO3dview}. 

%

%

\begin{figure}
\centering
\includegraphics[width=0.49\textwidth]{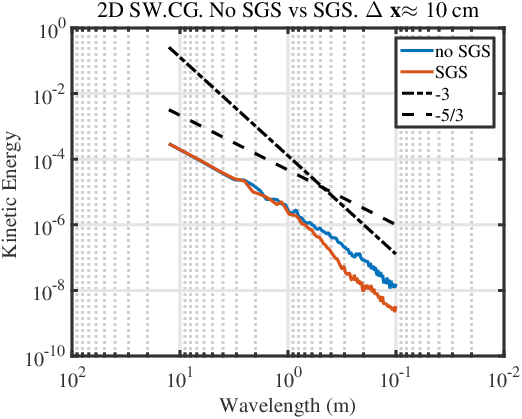}
\includegraphics[width=0.49\textwidth]{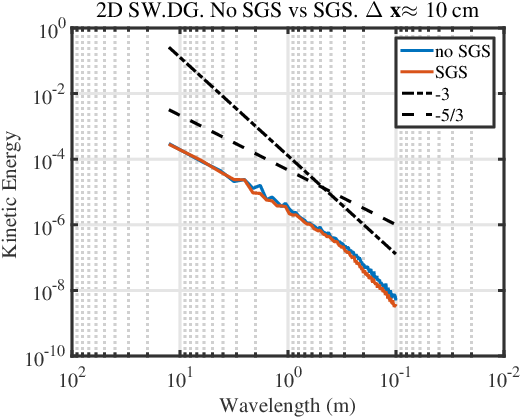}

\includegraphics[width=0.49\textwidth]{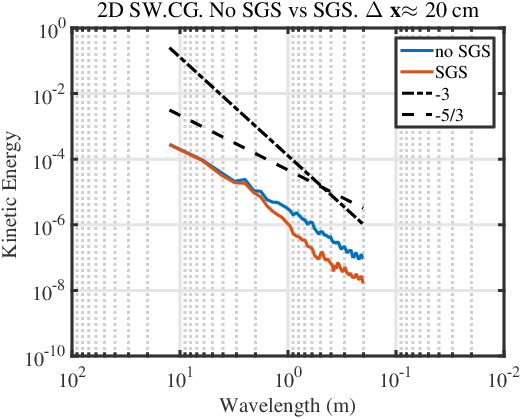}
\includegraphics[width=0.49\textwidth]{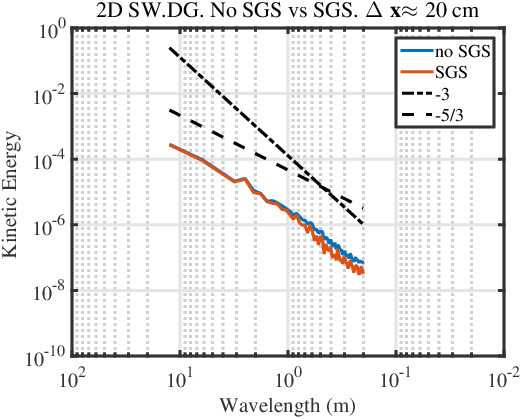}

\includegraphics[width=0.49\textwidth]{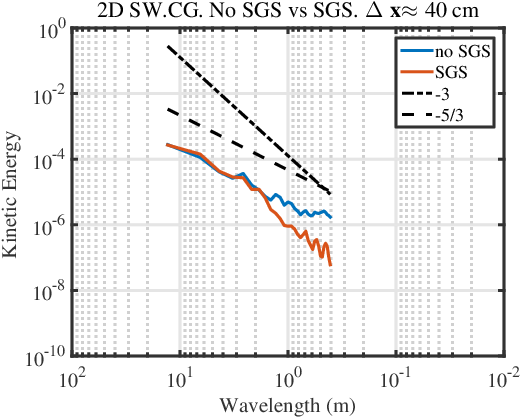}
\includegraphics[width=0.49\textwidth]{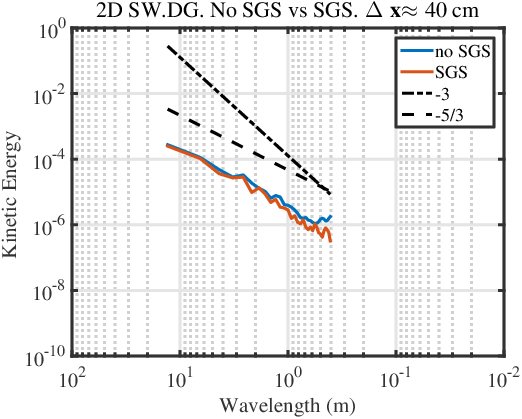}
\caption{Instantaneous 1D energy spectra of the single hill problem of Fig.\ \ref{briggs4Resolutions} at $t=50$ s.  
Left: CG with and without viscosity. Right: DG with and without viscosity. From top to bottom the resolution decreases.}
\label{briggs1995SpectraCGvsDG}
\end{figure}

The instantaneous energy spectra of the viscous and inviscid CG and DG solutions 
at $t=50$ s are plotted in Fig.\ \ref{briggs1995SpectraCGvsDG}. 
The difference between the CG and DG curves is striking.
The viscous and inviscid DG spectra overlap almost fully and 
show approximately the same decay across the entire spectrum, from a -5/3 slope in the 
inertial sub-range to a -3 slope in the dissipation wave numbers (refer to \cite{boffettaEcke2012} for a review on two-dimensional flows and their energetics).
This is only true as long as the resolution is not too coarse, especially so in the case of CG.
At very coarse resolutions ($\Delta {\bf x} \geq 0.4$ m), neither CG or DG can 
avoid energy from building up in the highest modes unless artificial diffusion is used. 
The inherent dissipation of DG is no longer sufficient.

\begin{sidewaysfigure}
\centerfloat
\includegraphics[trim=2cm 2cm 2cm 0cm, clip=true,width=0.2\textwidth]{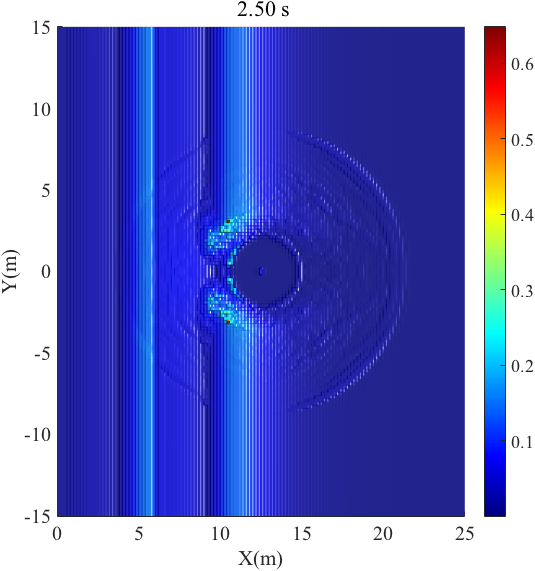}
\includegraphics[trim=2cm 2cm 2cm 0cm, clip=true,width=0.2\textwidth]{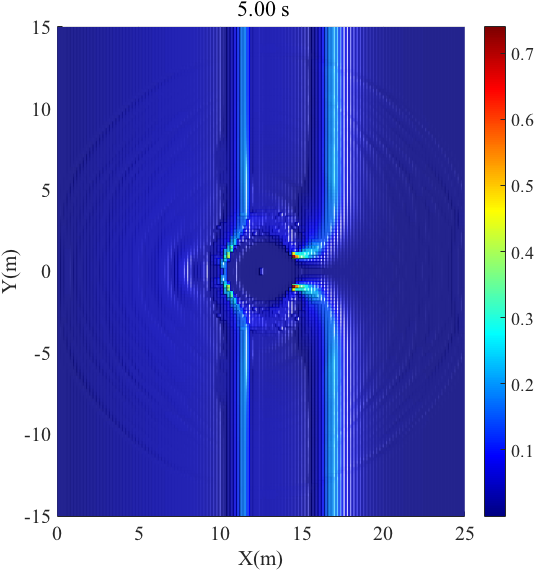}
\includegraphics[trim=2cm 2cm 2cm 0cm, clip=true,width=0.2\textwidth]{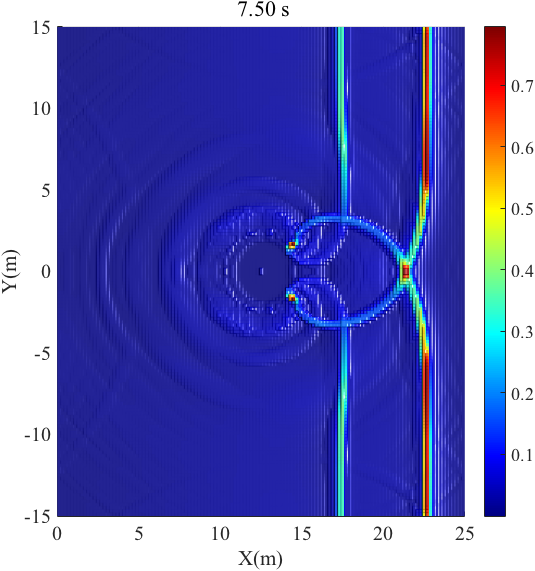}
\includegraphics[trim=2cm 2cm 2cm 0cm, clip=true,width=0.2\textwidth]{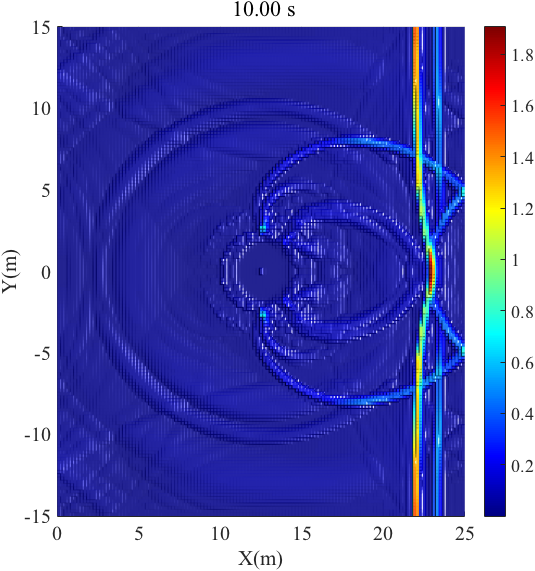}
\includegraphics[trim=2cm 2cm 2cm 0cm, clip=true,width=0.2\textwidth]{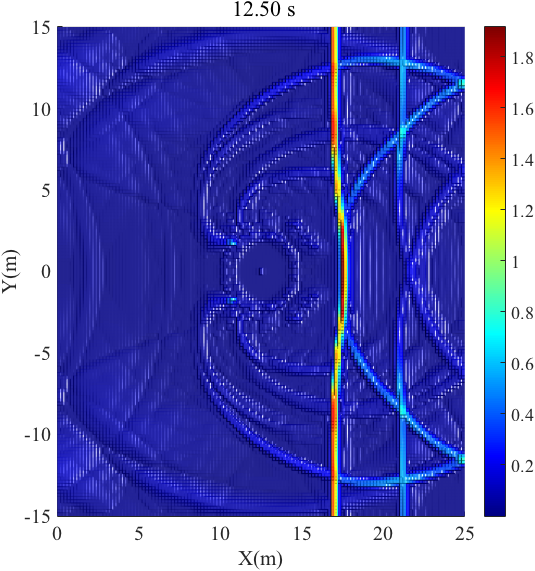}\\
\includegraphics[trim=2cm 2cm 2cm 0cm, clip=true,width=0.2\textwidth]{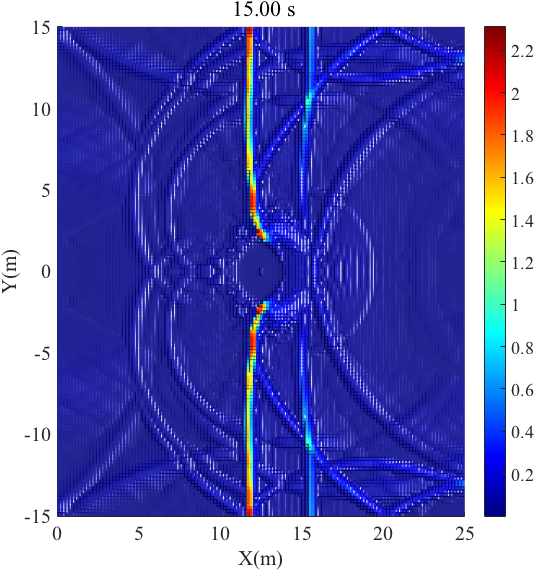}
\includegraphics[trim=2cm 2cm 2cm 0cm, clip=true,width=0.2\textwidth]{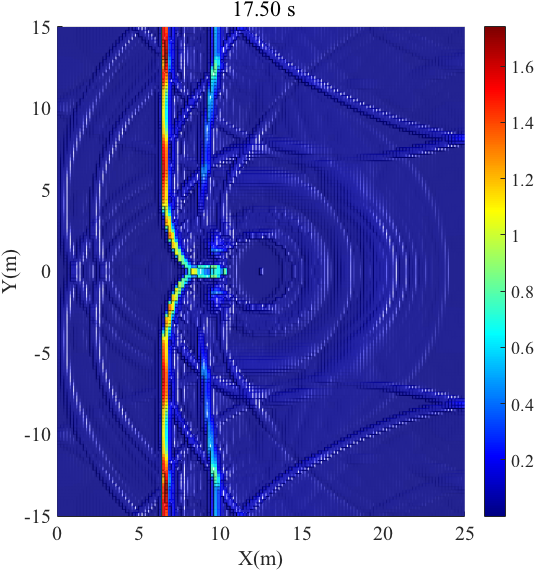}
\includegraphics[trim=2cm 2cm 2cm 0cm, clip=true,width=0.2\textwidth]{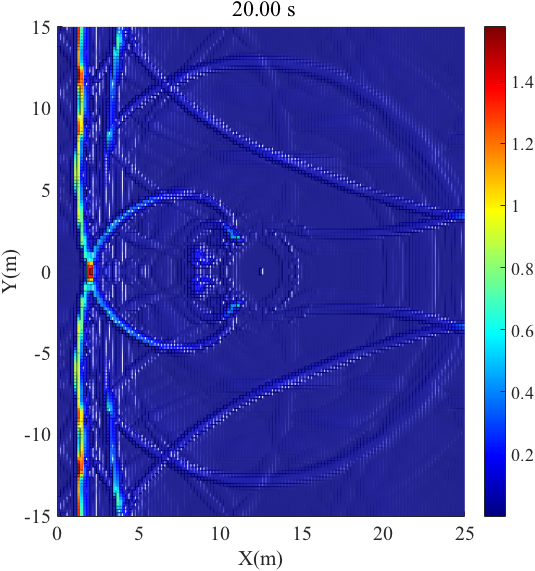}
\includegraphics[trim=2cm 2cm 2cm 0cm, clip=true,width=0.2\textwidth]{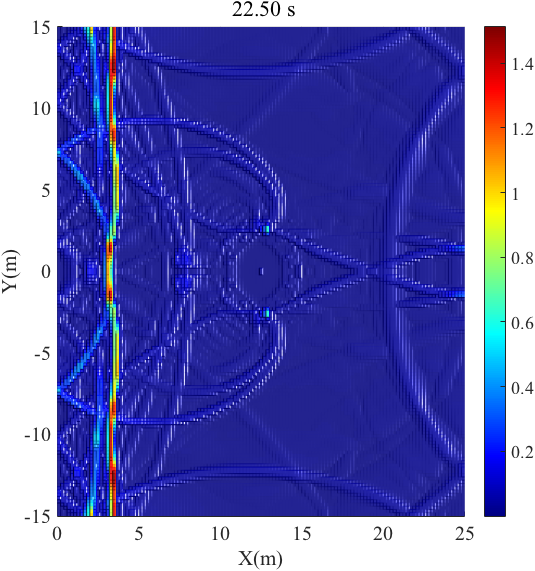}
\includegraphics[trim=2cm 2cm 2cm 0cm, clip=true,width=0.2\textwidth]{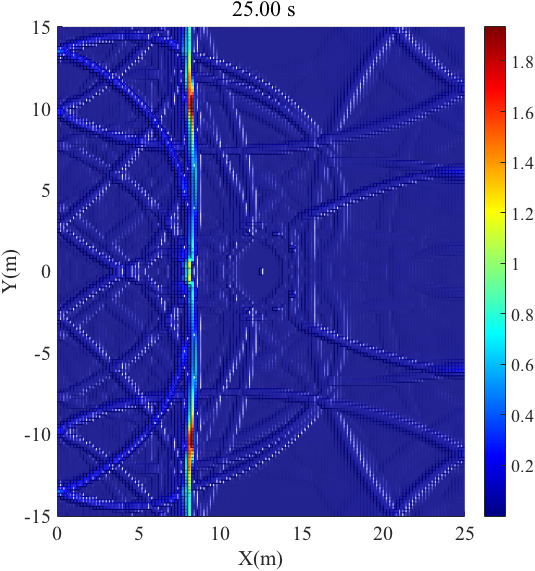}\\
\includegraphics[trim=2cm 2cm 2cm 0cm, clip=true,width=0.2\textwidth]{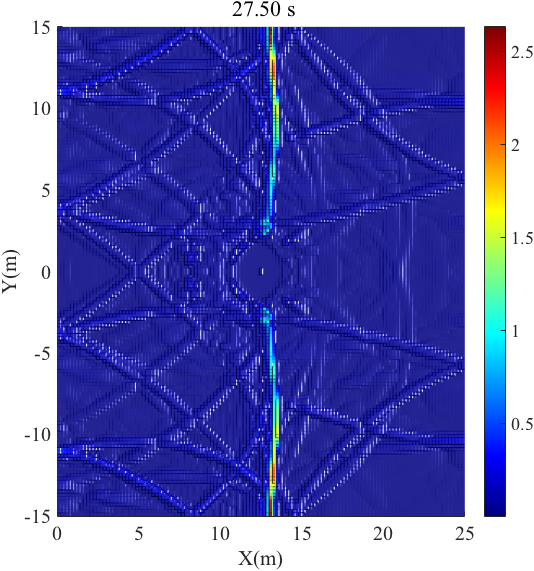}
\includegraphics[trim=2cm 2cm 2cm 0cm, clip=true,width=0.2\textwidth]{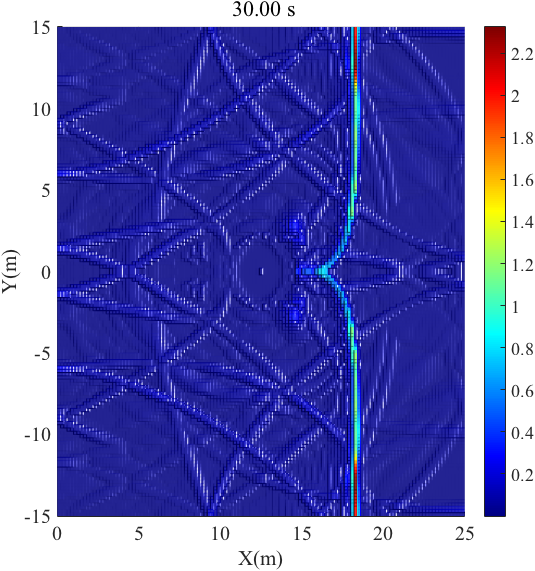}
\includegraphics[trim=2cm 2cm 2cm 0cm, clip=true,width=0.2\textwidth]{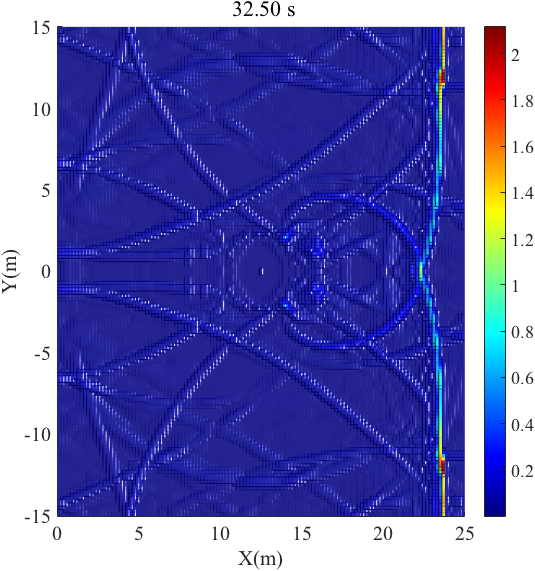}
\includegraphics[trim=2cm 2cm 2cm 0cm, clip=true,width=0.2\textwidth]{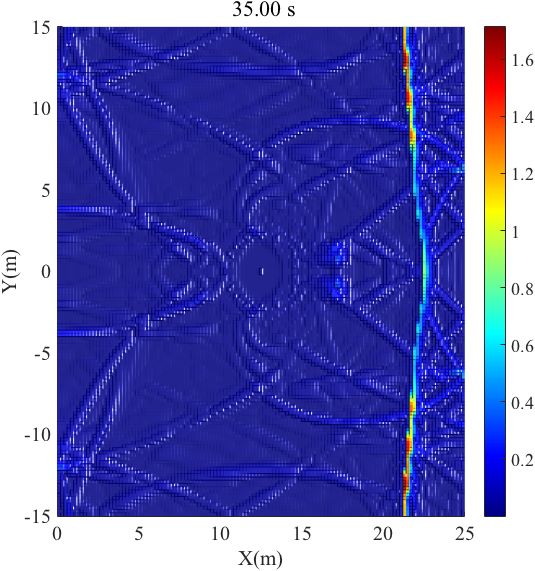}
\includegraphics[trim=2cm 2cm 2cm 0cm, clip=true,width=0.2\textwidth]{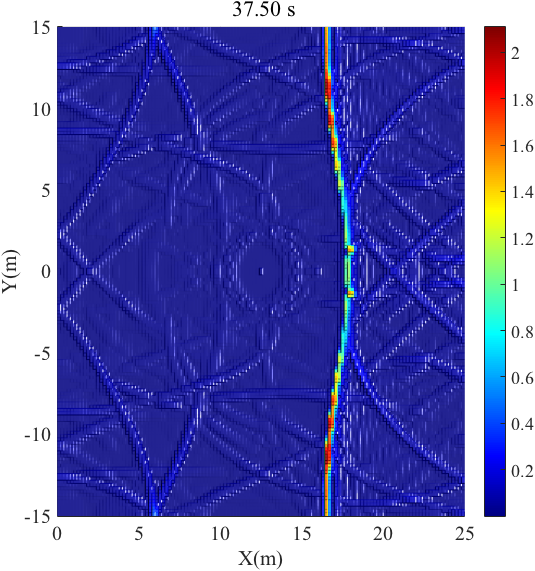}\\
\includegraphics[trim=2cm 2cm 2cm 0cm, clip=true,width=0.2\textwidth]{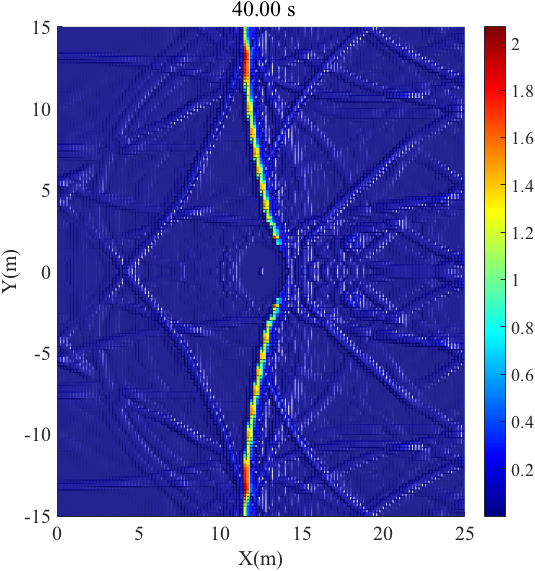}
\includegraphics[trim=2cm 2cm 2cm 0cm, clip=true,width=0.2\textwidth]{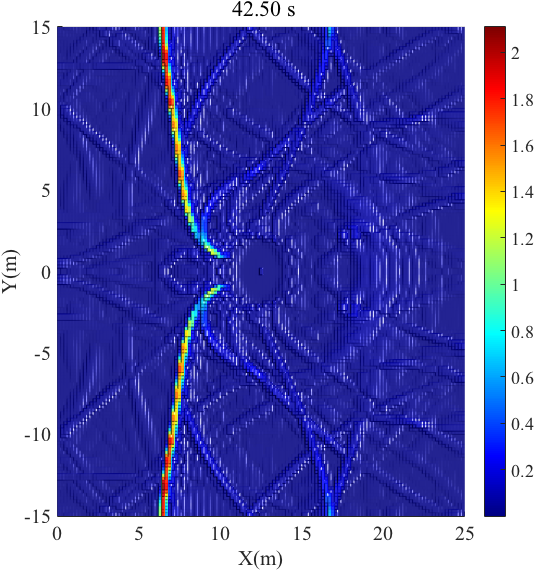}
\includegraphics[trim=2cm 2cm 2cm 0cm, clip=true,width=0.2\textwidth]{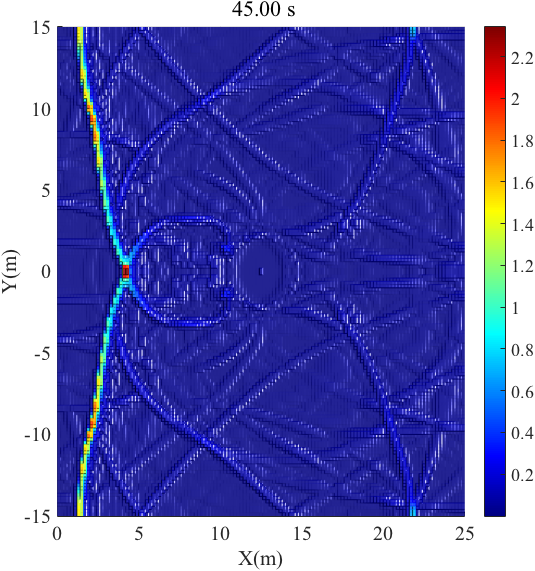}
\includegraphics[trim=2cm 2cm 2cm 0cm, clip=true,width=0.2\textwidth]{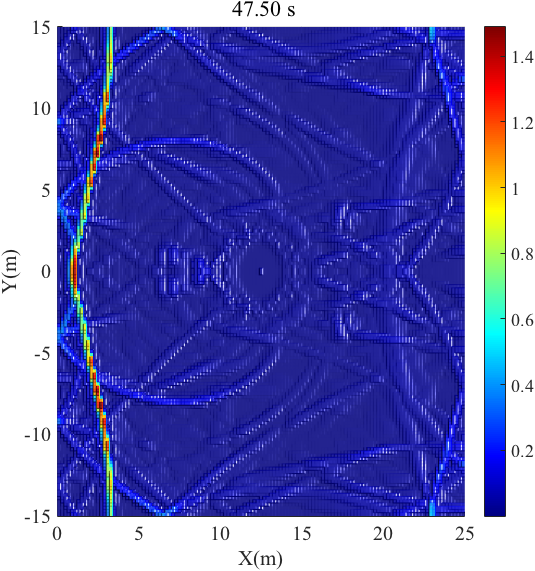}
\includegraphics[trim=2cm 2cm 2cm 0cm, clip=true,width=0.2\textwidth]{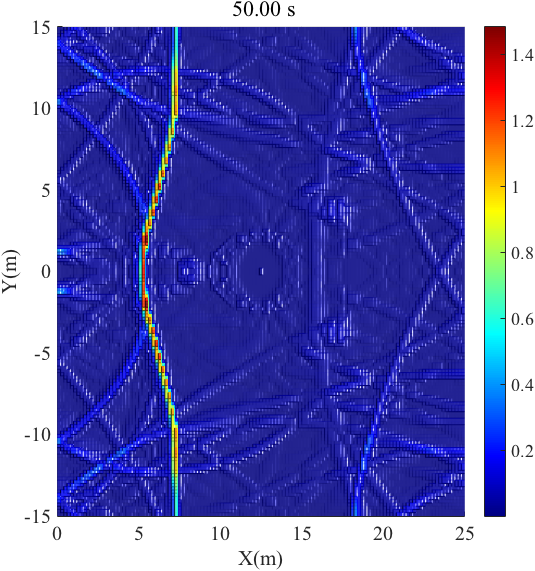}
\caption{Time evolution of $\mu_{SGS}$ from $t=0$ s to $t=50$ s. The colorscale ranges between 0 (blue) and 5 (red) m$^2/$s. 
The plotted domain is $\Omega=[0, 25]\times [-15, 15]$ m$^2$.}
\label{nusgs_timeevolution}
\end{sidewaysfigure}

\begin{sidewaysfigure}
\centerfloat
\includegraphics[trim=2cm 2cm 2.5cm 0cm, clip=true,width=0.2\textwidth]{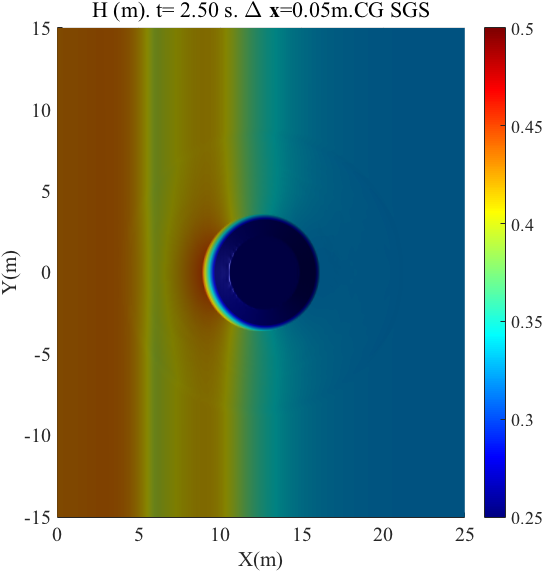}
\includegraphics[trim=2cm 2cm 2.5cm 0cm, clip=true,width=0.2\textwidth]{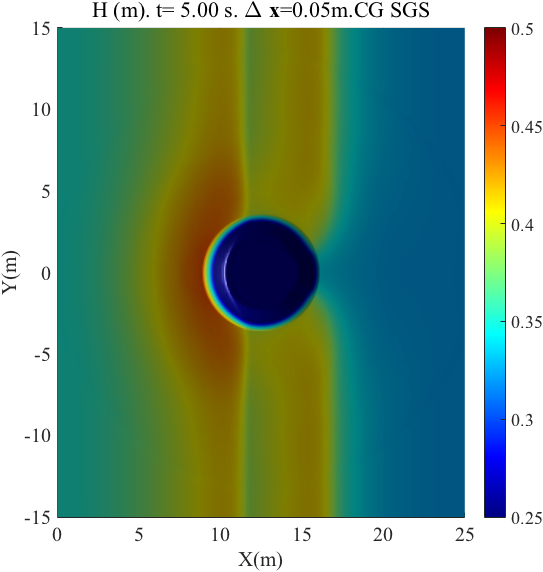}
\includegraphics[trim=2cm 2cm 2.5cm 0cm, clip=true,width=0.2\textwidth]{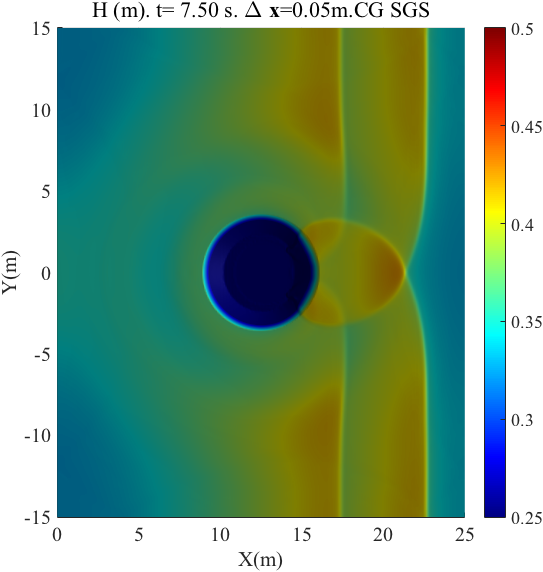}
\includegraphics[trim=2cm 2cm 2.5cm 0cm, clip=true,width=0.2\textwidth]{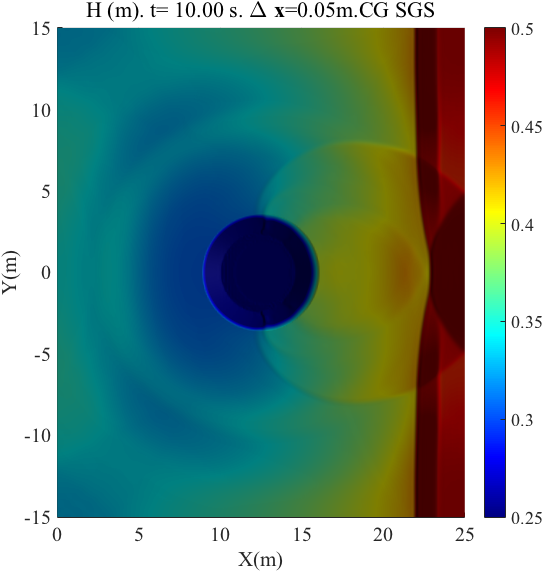}
\includegraphics[trim=2cm 2cm 2.5cm 0cm, clip=true,width=0.2\textwidth]{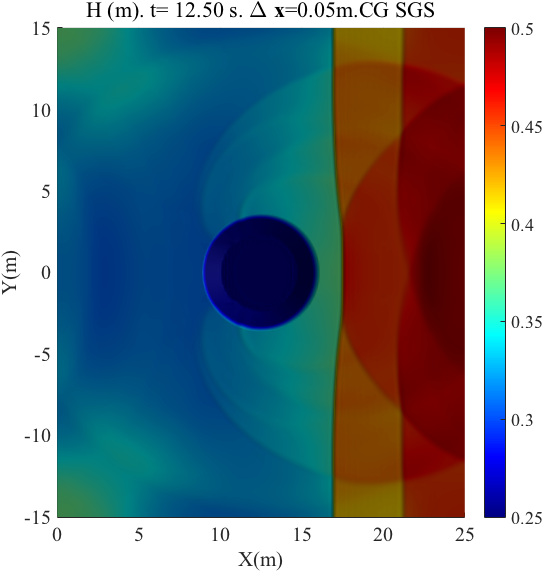}\\
\includegraphics[trim=2cm 2cm 2.5cm 0cm, clip=true,width=0.2\textwidth]{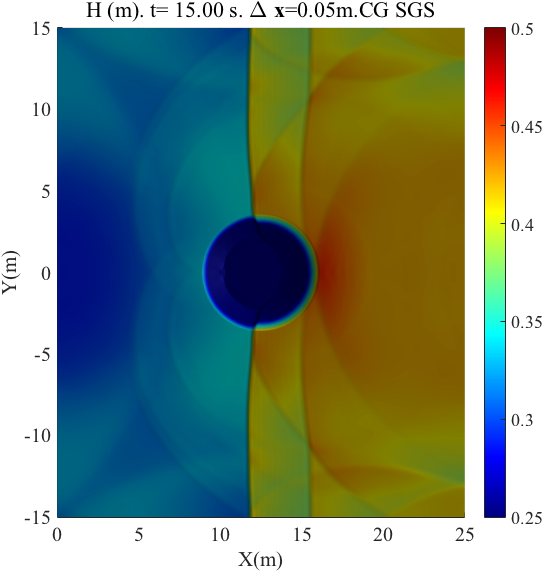}
\includegraphics[trim=2cm 2cm 2.5cm 0cm, clip=true,width=0.2\textwidth]{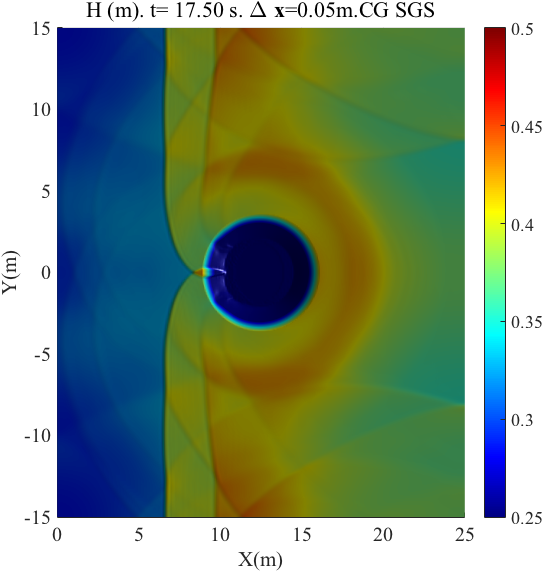}
\includegraphics[trim=2cm 2cm 2.5cm 0cm, clip=true,width=0.2\textwidth]{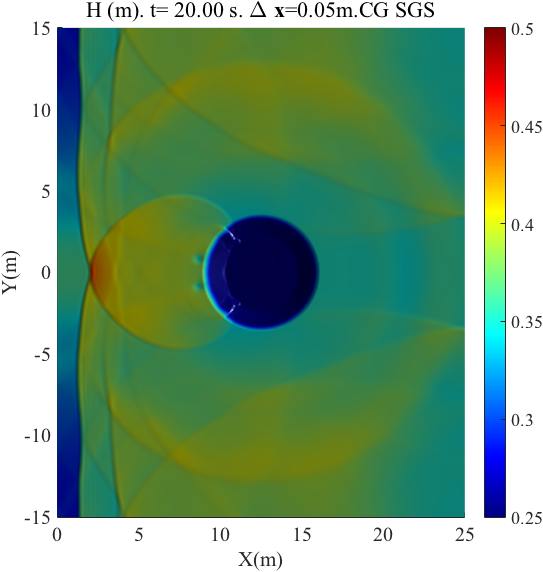}
\includegraphics[trim=2cm 2cm 2.5cm 0cm, clip=true,width=0.2\textwidth]{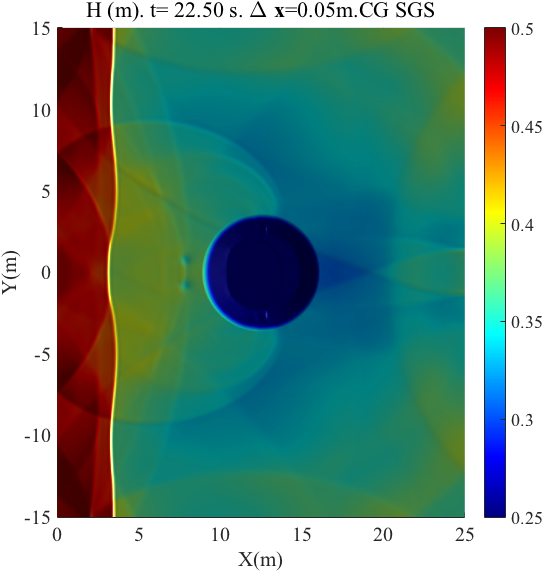}
\includegraphics[trim=2cm 2cm 2.5cm 0cm, clip=true,width=0.2\textwidth]{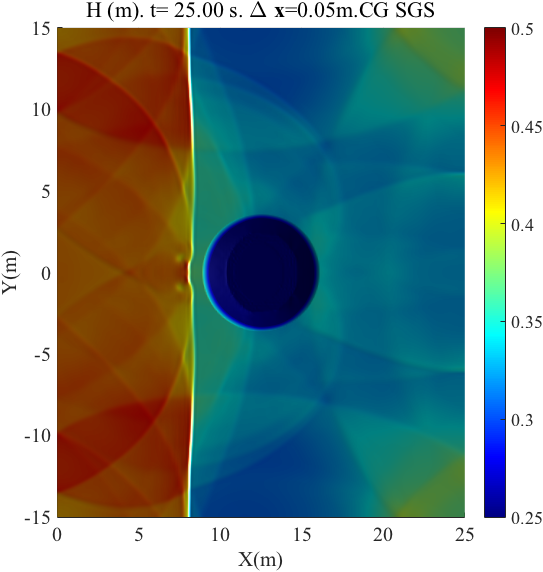}\\
\includegraphics[trim=2cm 2cm 2.5cm 0cm, clip=true,width=0.2\textwidth]{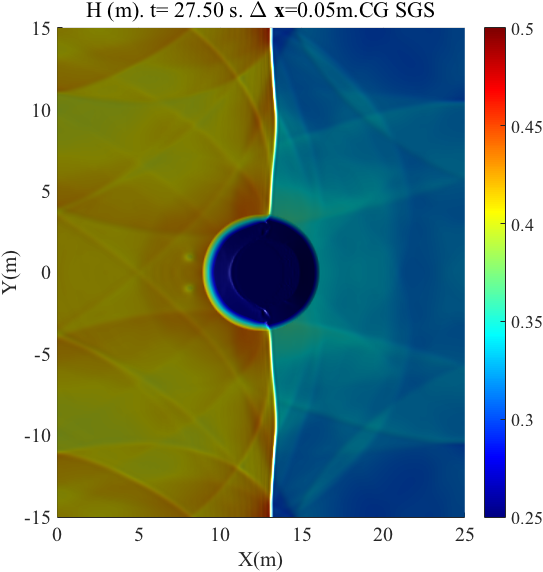}
\includegraphics[trim=2cm 2cm 2.5cm 0cm, clip=true,width=0.2\textwidth]{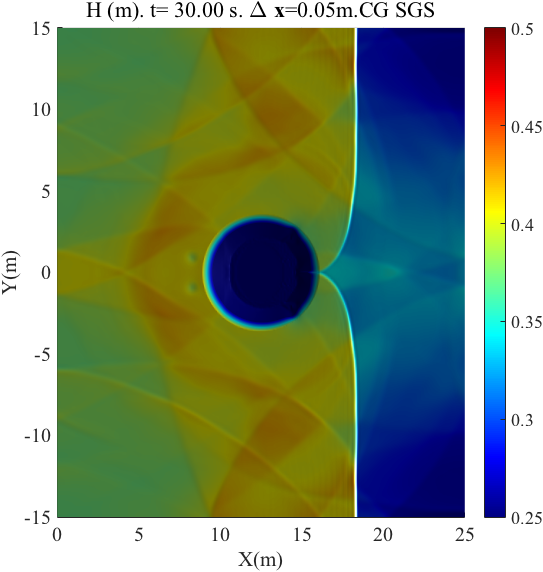}
\includegraphics[trim=2cm 2cm 2.5cm 0cm, clip=true,width=0.2\textwidth]{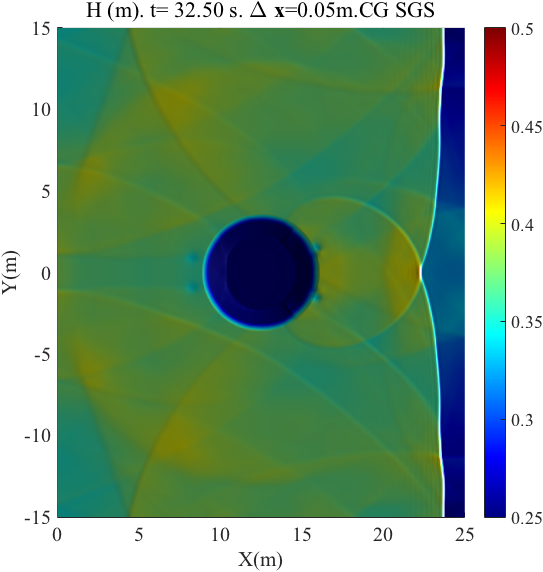}
\includegraphics[trim=2cm 2cm 2.5cm 0cm, clip=true,width=0.2\textwidth]{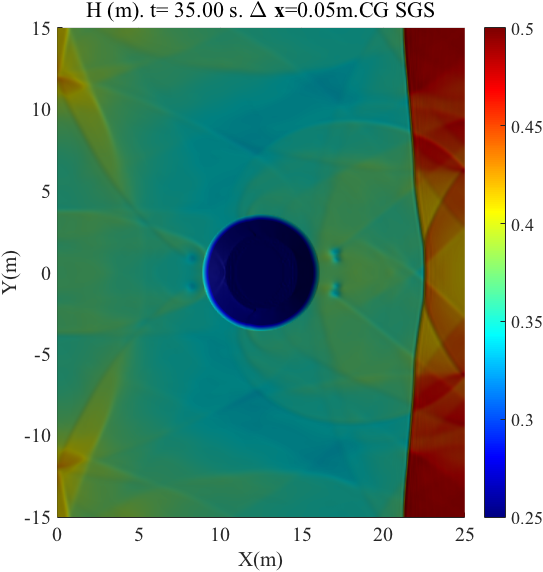}
\includegraphics[trim=2cm 2cm 2.5cm 0cm, clip=true,width=0.2\textwidth]{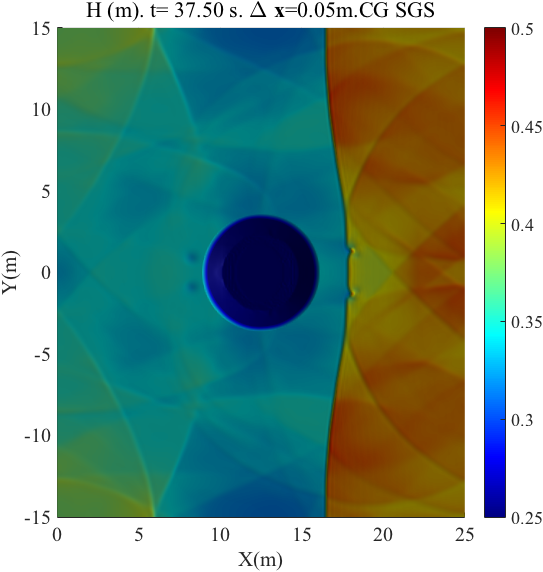}\\
\includegraphics[trim=2cm 2cm 2.5cm 0cm, clip=true,width=0.2\textwidth]{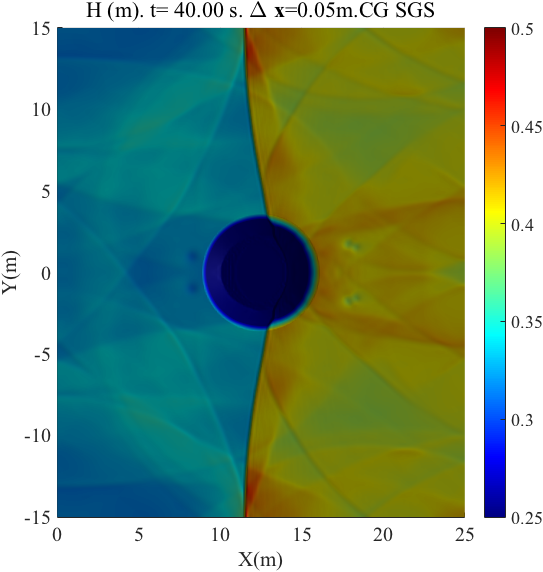}
\includegraphics[trim=2cm 2cm 2.5cm 0cm, clip=true,width=0.2\textwidth]{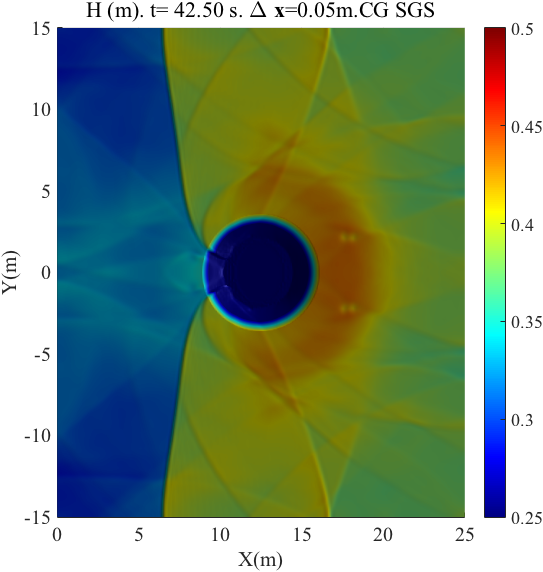}
\includegraphics[trim=2cm 2cm 2.5cm 0cm, clip=true,width=0.2\textwidth]{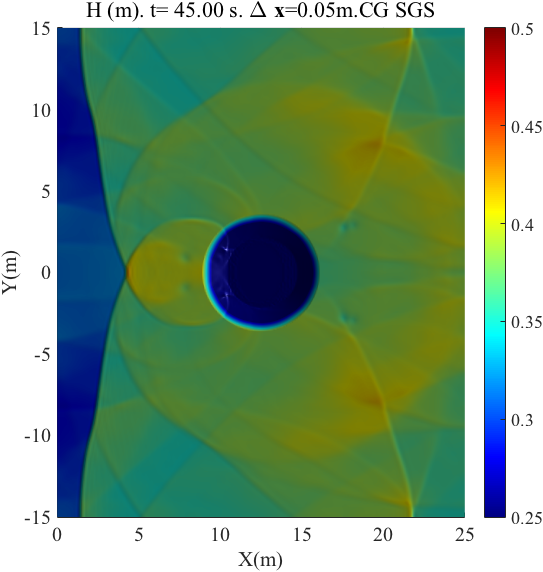}
\includegraphics[trim=2cm 2cm 2.5cm 0cm, clip=true,width=0.2\textwidth]{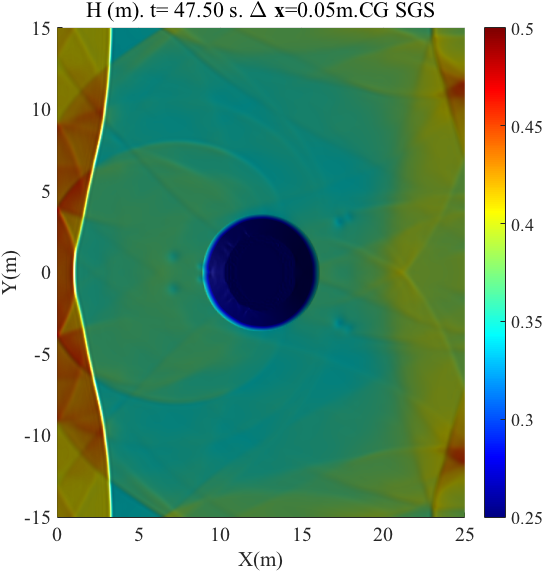}
\includegraphics[trim=2cm 2cm 2.5cm 0cm, clip=true,width=0.2\textwidth]{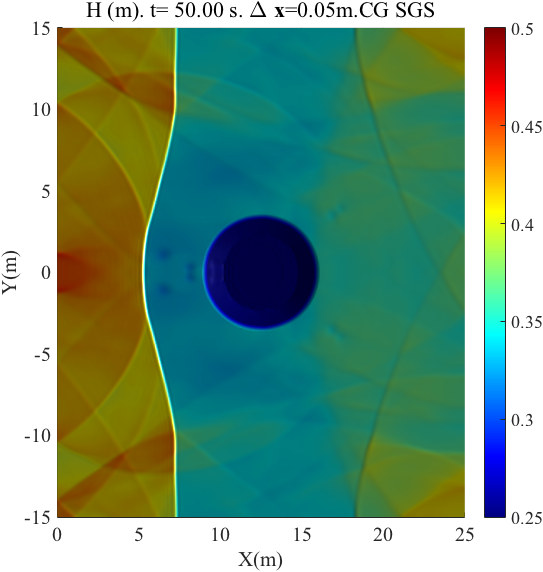}
\caption{Time evolution of $H$ from $t=0$ s to $t=50$ s. The color scale ranges between 0.25 (blue) and 0.5 (red) m. 
The plotted domain is $\Omega=[0, 25]\times [-15, 15]$ m$^2$.}
\label{H_timeevolution}
\end{sidewaysfigure}

We stated above that $\mu_{SGS}$ is only active where the equation residuals (i.e. gradients) are important. In the case of 
water waves, this occurs in the proximity of the wave fronts. 
In Fig.\ \ref{nusgs_timeevolution}, we plot $\mu_{SGS}$ to show its spatial structure and its evolution between $t=0$ and $t=50$ seconds.
This plots clearly show how diffusion is equally zero away from the fronts and only activates where really necessary. 
It may not be so obvious to achieve this by using an artificial diffusion that is not residual-based. 
To provide a visual correlation between $\mu_{SGS}$ and the wave features, in Fig.\ \ref{H_timeevolution} we plot the stabilized spectral element solution 
of the water surface at a grid resolution $\Delta x \approx 0.05$ m.

\section{Conclusions}
\label{conclusionsSCT}
We presented the numerical solution of the shallow water equations via continuous and discontinuous Galerkin (CG/DG)
methods to model problems involving inundation. 
Careful handling of the transition between dry and wet surfaces is particularly challenging for high-order methods.
The most adopted solution in these cases consists of lowering the approximation order in either the transition zone alone, or in the whole domain.
By extending to CG a simple slope limiter originally designed for DG \cite{xingZhangShu2010}, 
and by combining it with a thin water layer always present in dry regions, 
we showed that the nominal high-order accuracy of the underlying space approximations was preserved 
and that the solution smoothness in the proximity of moving boundaries was maintained in one- and two-dimensions.
The simplicity of this approach is effective for CG as well as DG. More important, we demonstrated its effectiveness using higher order elements.
This limiter, unless embedded with a positivity-preserving scheme for water height, 
does not prevent the high-order solution from 
triggering unphysical Gibbs oscillations in the proximity of strong gradients.
To overcome this problem, we presented a dynamically adaptive dissipation based on a residual-based sub-grid scale eddy viscosity model ({\it Dyn-SGS}).
By numerical examples, we demonstrated the following properties of this model when solving the shallow water equations:
$(i)$ It removes the Gibbs oscillations that form in the proximity of sharp wave fronts while preserving 
the overall accuracy and sharpness of the solution everywhere else. This is possible thanks to the residual-based definition of the 
dynamic diffusion coefficient.
$(ii)$ For coarse grids, it prevents energy from building up at small wave-numbers; this is very important to preserve 
numerical stability in the flow regimes we are interested in. $(iii)$ The model has no user tunable parameter, which is of great advantage 
when the model is to be used by an external user. 
It is important to underline that the natural, built-in dissipation of DG be may be large enough that the contribution of {\it Dyn-SGS} becomes irrelevant.
When this happens, the dynamic dissipation detects it from the residual, and hence limits its own strength.
Finally, a three stage, second order explicit-first-stage, singly diagonally implicit Runge-Kutta 
(ESDIRK) time integration scheme was implemented to overcome the small time-step restriction incurred by high-order Galerkin approximations.

\section{Acknowledgments}
The authors would like to acknowledge the contribution of Haley Lane, 
who implemented the one-dimensional version of the wetting and drying algorithm used in this work.
The authors would also like to acknowledge Karoline Hood who tested the 
correctness of the implicit solver in her NPS Master's thesis (\cite{karolineTHESIS2016}).
The authors are also thankful to Oliver Fringer and the Rogers for discussions regarding coastal flows, and 
to  Stephen R. Guimond for providing his MATLAB functions to compute the energy spectra.
FXG acknowledges the support of the ONR Computational Mathematics program, and FXG and EMC acknowledge the support of AFOSR Computational Mathematics.

\newpage

\end{document}